\documentclass[3p,sort & compress]{elsarticle}
\usepackage{amssymb}
\usepackage{latexsym}
\usepackage{amsfonts}
\usepackage{amsthm}
\usepackage{amsbsy}
\usepackage{amsmath}
\usepackage{graphics,graphicx}
\usepackage{subfigure}
\usepackage{multirow,multicol}
\usepackage{bm}
\usepackage{bbm}
\usepackage{color}
\usepackage{float}
\usepackage{calc}
\usepackage{caption}
\usepackage{dsfont}
\usepackage{ifthen}
\usepackage{mathrsfs}
\usepackage[colorlinks,linkcolor=blue,citecolor=blue]{hyperref}
\usepackage{epstopdf}
\usepackage{epsfig}
\usepackage{algorithm}
\usepackage{algorithmic}
\usepackage{pifont}
\usepackage{natbib}
\usepackage{overpic}
\usepackage{psfrag}
\usepackage{rotating}
\usepackage{stmaryrd}
\usepackage{verbatim}
\usepackage{setspace}
\usepackage{tikz}
\usetikzlibrary{fpu}

\newdefinition{remark}{Remark}

\newdefinition{method}{Method}
\newdefinition{example}[theorem]{Example}

\numberwithin{theorem}{section}

\newcommand{\orcid}[1]{\href{https://orcid.org/#1}{\includegraphics[width=8pt]{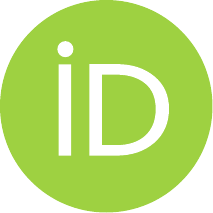}}}
\journal{Journal of \LaTeX\ Templates}
\begin{document}
\begin{frontmatter}
\title{Augmented physics informed extreme learning machine to solve the biharmonic equations via Fourier expansions}
\author[Ceyear,SDUInfo]{Xi'an Li\orcid{0000-0002-1509-9328}}\ead{lixian9131@163.com}
\author[acu]{Jinran Wu\orcid{0000-0002-2388-3614}} \ead{ryan.wu@acu.edu.au}
\author[qu]{Yujia Huang} \ead{rainfamilyh@gmail.com}
\author[qut]{Zhe Ding\orcid{0000-0002-5123-0135}\corref{cor1}}\ead{zhe.ding@hdr.qut.edu.au}
\author[Ceyear]{Xin Tai}\ead{taixin@ceyear.com}
\author[Ceyear]{Liang Liu}\ead{liuliang@ceyear.com}
\author[qu]{You-Gan Wang\orcid{0000-0003-0901-4671}}\ead{ygwanguq2012@gmail.com}
\cortext[cor1]{Corresponding author.}
\address[Ceyear]{Ceyear Technologies Co., Ltd, Qingdao 266555, China}
\address[SDUInfo]{School of Information Science and Engineering, Shandong University, Qingdao 266237, China}
\address[acu]{Australian Catholic University, Banyo 4014, Australia}
\address[qu]{The University of Queensland, St Lucia 4072, Australia}
\address[qut]{Queensland University of Technology, Brisbane 4001, Australia}

\begin{abstract} 
To address the sensitivity of parameters and limited precision for physics-informed extreme learning machines (PIELM) with common activation functions, such as sigmoid, tangent, and Gaussian, in solving high-order partial differential equations (PDEs) relevant to scientific computation and engineering applications, this work develops a Fourier-induced PIELM (FPIELM) method. This approach aims to approximate solutions for a class of fourth-order biharmonic equations with two boundary conditions on both unitized and non-unitized domains. By carefully calculating the differential and boundary operators of the biharmonic equation on discretized collections, the solution for this high-order equation is reformulated as a linear least squares minimization problem. We further evaluate the FPIELM with varying hidden nodes and scaling factors for uniform distribution initialization, and then determine the optimal range for these two hyperparameters. Numerical experiments and comparative analyses demonstrate that the proposed FPIELM method is more stable, robust, precise, and efficient than other PIELM approaches in solving biharmonic equations across both regular and irregular domains.
    
\end{abstract}
\end{frontmatter}

\section{Introduction} \label{sec:introduction}
Partial differential equations (PDEs), which have their roots in science and engineering, have been instrumental in advancing the fields of physics, chemistry, biology, image processing, and more. In general, the analytical solutions of PDEs are often elusive, necessitating the reliance on numerical methods and approximation techniques for their resolution. This paper focuses on the numerical solution of the following biharmonic equation:
\begin{equation}\label{eq:biharmonic}
\Delta^2 u(\bm{x}) = f(\bm{x}),   ~\text{in}~ \Omega,
\end{equation}
where $\Omega$ is a polygonal or polyhedral domain of dimension $d$ in Euclidean space, with a piecewise Lipschitz boundary satisfying an interior cone condition. Here, $f(\bm{x}) \in L^2(\Omega)$ is a prescribed function, and $\Delta $ denotes the standard Laplace operator. The expression of the biharmonic operator $\Delta^2 u(\bm{x})$ is given by:
\begin{equation*}
\Delta^2 u(\bm{x}) = \sum_{i=1}^{d}\frac{\partial^4 u}{\partial x_i^4} + \sum_{i=1}^{d}\sum_{j=1, j\neq i}^{d}\frac{\partial^4 u}{\partial x_i^2 \partial x_j^2}.
\end{equation*}
In practice, the solution of equation \eqref{eq:biharmonic} can be uniquely and precisely determined by the one of the following two boundary conditions:
\begin{itemize}
    \item Dirichlet boundary conditions:
    \begin{equation}\label{bd:Dirichlet}
    u(\bm{x}) = g(\bm{x}) ~~\text{and}~~ \frac{\partial u(\bm{x})}{\partial \vec{n}} = h(\bm{x}) ~~\text{on}~ \partial\Omega;
    \end{equation}
    \item Navier boundary conditions:
    \begin{equation}\label{bd:Navier}
    u(\bm{x}) = g(\bm{x}) ~~\text{and}~~\Delta u(\bm{x}) = k(\bm{x}) ~~\text{on}~ \partial\Omega,
    \end{equation}
\end{itemize}
where $\partial\Omega$ denotes the boundary of $\Omega$, and $\vec{n} = (n_1, n_2, \dots, n_d)$ is the outward normal vector on $\partial\Omega$. The functions $g$, $h$, and $k$ describe the boundary behavior of the solution and are assumed to be sufficiently smooth.

The biharmonic equation, a high-order PDE, frequently occurs in fields such as elasticity and fluid mechanics. Its application ranges from scattered data fitting with thin plate splines~\cite{wahba1990spline} to the simulation of fluid mechanics~\cite{greengard1998an, ferziger1996computational} and linear elasticity~\cite{christiansen1975integral, constanda1995the}. Solving it is challenging due to the fourth-order derivative terms, particularly in complex geometries. Various traditional numerical methods for solving \eqref{eq:biharmonic} can be broadly classified into direct (uncoupled) methods and coupled (splitting, mixed) methods.

In the direct approach, methodologies include finite difference methods(FDM)~\cite{gupta1979direct, altas1998multigrid, ben-artzi2008a, bialecki2012a}, finite volume techniques~\cite{chun-jiabi2004mortar, wang2004a, eymard2012finite}, and finite element methods (FEM) based on the variational formulation, such as non-conforming FEM~\cite{Baker1977Finite, Lascaux1975Some, Morley2016The} and conforming FEM~\cite{Zienkiewicz2005The, Ciarlet1978The}. In the coupled approach, the biharmonic equation is decomposed into two interrelated Poisson equations, solvable by methods like FDM, FEM, and mixed element methods~\cite{brezzi1991mixed, cheng2000some, davini2000an, lamichhane2011a, stein2019a}. Additional methods include the collocation method~\cite{mai-duy2009a, bialecki2010spectral, bialecki2020a} and the radial basis functions (RBF) method~\cite{maiduy2005an, adibi2007numerical, li2011a}.

Conventional methods, however, encounter challenges with irregular domains and high dimensions, often referred to as the curse of dimensionality. Machine learning, especially artificial neural networks (ANNs), has shown promise in approximating the solutions of PDEs, potentially overcoming the limitations of traditional numerical techniques. Various ANN approaches, such as extreme learning machines (ELM) and deep neural networks (DNN), have been developed to solve PDEs, including \eqref{eq:biharmonic}. Physics-informed neural network methods have also been proposed, minimizing loss and incorporating residuals and boundary conditions~\cite{guo2019a, huang2021ThinPlate, vahab2022physics}. Despite their potential, these methods may suffer from limited accuracy, uncertain convergence, and slow training. Recently, a physics-informed extreme learning machine (PIELM) approach was introduced, combining physical laws with ELM to address various linear or nonlinear PDEs~\cite{dwivedi2020physics, dwivedi2020solution, li2023extreme, wang2023extreme,joshi2024physics,liu2023bayesian}. This method retains the advantages of ELM, including mesh-free operation, rapid learning, and parallel architecture compatibility, while effectively handling complex PDEs.

The choice of activation function is critical for ELM networks. For regression and classification, the sigmoid $1/(1+e^{-x})$ and Gaussian $e^{-x^2}$ functions are often effective~\cite{chen2013modified, ding2015extreme}. For solving PDEs, experiments have shown that sigmoid, Gaussian, and hyperbolic tangent functions are ideal, producing high-precision solutions~\cite{dwivedi2020physics, fabiani2021numerical, calabro2021extreme, li2023extreme}. However, these functions may impede information transmission in high-order problems with large-scale inputs, as their high-order derivatives are complex. Fourier-type functions, such as $\sin(x)$, have been suggested as alternative to mitigate these issues~\cite{huang2006universal, huang2015local, rahimi2008uniform}. Moreover, initializing ELM’s internal weights and biases uniformly within $[-\delta, \delta]$ with $\delta > 0$ has been shown to enhance performance~\cite{dong2021local, dong2022computing, dong2021modified}.

This work addresses the biharmonic equations~\eqref{eq:biharmonic} with Dirichlet and Navier boundary conditions by combining the Fourier spectral theorem with PIELM, resulting in a Fourier-induced PIELM model (FPIELM). Various hidden node configurations and uniform initialization scale factors are tested to optimize hyperparameters. We further compare the performance of FPIELM against that of the PIELM model with sigmoid, tangent, and Gaussian activation functions. Numerical results indicate that FPIELM achieves exceptional performance and robustness for solving \eqref{eq:biharmonic} on both unitized and non-unitized domains.

The paper is organized as follows: Section \ref{sec:02} introduces the classical ELM framework. Section \ref{sec:03} presents the FPIELM algorithm construction. Numerical examples evaluating FPIELM's accuracy and efficiency are shown in Section \ref{sec:05}, and conclusions are drawn in Section \ref{sec:06}.

\section{Preliminary for ELM }\label{sec:02}

This section provides a concise overview of essential concepts and the mathematical formulation of the ELM model. ELM falls within the category of single hidden layer fully connected neural networks (SHFNN) in artificial neural networks and is designed to establish a mapping $\ell$: $\mathbb{R}^d \rightarrow \mathbb{R}^{k}$. A key characteristic of ELM is the random assignment of weights and biases in the input layer, while the output weights are determined analytically~\cite{schmidt1992feed, igelnik1995stochastic}. Mathematically, the ELM mapping $\ell(\bm{x})$ with input $\bm{x} \in \mathbb{R}^{d}$ and output $\bm{y} \in \mathbb{R}^k$ is defined as:
\begin{equation}\label{form2elm}
\bm{y} = (y_1, y_2, \cdots, y_k)~\textup{with}~y_{m} = \sum_{i=1}^{n} \bm{\beta}_{mi} \sigma\left(\sum_{j=1}^{d} w_{ij} x_j + b_i\right).
\end{equation}
For convenience, we denote $\bm{W} = \left(\bm{W}^{[1]}, \bm{W}^{[2]}, \cdots, \bm{W}^{[n]}\right)^T \in \mathbb{R}^{n \times d}$, where $\bm{W}^{[i]} = (w_{i1}, w_{i2}, \cdots, w_{id})^T \in \mathbb{R}^{d}$ represents the weight of the $i$-th hidden unit in ELM, and $\bm{b} = (b_1, b_2, \cdots, b_n)^T \in \mathbb{R}^{n}$ is the corresponding bias vector. The activation function $\sigma(\cdot)$ operates element-wise. The matrix $\bm{\beta} = \left(\bm{\beta}^{[1]}, \bm{\beta}^{[2]}, \cdots, \bm{\beta}^{[k]}\right)^T \in \mathbb{R}^{k \times n}$, where $\bm{\beta}^{[i]} = (\beta_{i1}, \beta_{i2}, \cdots, \beta_{in})^T \in \mathbb{R}^n$, represents the weights connecting the hidden units to the outputs. A schematic of the ELM with $n$ hidden units is shown in Fig.~\ref{fig2elm}.

\begin{figure}[htp]
	\centering
	\begin{tikzpicture}[scale=0.65]
	\node[circle] (x0) at (0, 3.5) {};
	\node[circle] (x1) at (0, 1.75) {};
	\node[circle] (x2) at (0, 0) {};
	\node[circle] (x3) at (0, -1.75) {};
	\node[circle] (x4) at (0, -3.5) {};
	
	\node[] (input) at (1.5, -5.5) {Input layer};
	
	\node[circle, fill=green!60,inner sep=5.5pt] (i0) at (2, 3.5) {};
	\node[circle, fill=green!60,inner sep=5.5pt] (i1) at (2, 1.75) {};	
	\node[circle, fill=green!60,inner sep=5.5pt] (i2) at (2, 0) {};	
	\node[circle, fill=green!60,inner sep=5.5pt] (i3) at (2, -1.75) {};	
	\node[circle, fill=green!60,inner sep=5.5pt] (i4) at (2, -3.5) {};	
	
	\draw[line width=1.0pt,->] (x0) -- node[above]{$x_0$} (i0);
	\draw[line width=1.0pt,->] (x1) -- node[above]{$x_1$} (i1);
	\draw[line width=1.0pt,->] (x2) -- node[above]{$x_2$} (i2);
	\draw[line width=1.0pt,->] (x3) -- node[above]{$\cdots$} (i3);
	\draw[line width=1.0pt,->] (x4) -- node[above]{$x_d$} (i4);
	
	\node[circle, fill=cyan!70,inner sep=1.5pt] (h10) at (6, 4.5) {\huge$\sigma$};
	\node[circle, fill=cyan!70,inner sep=1.5pt] (h11) at (6, 3) {\huge$\sigma$};
	\node[circle, fill=cyan!70,inner sep=1.5pt] (h12) at (6, 1.5) {\huge$\sigma$};
	\node[circle, fill=cyan!70,inner sep=1.5pt] (h13) at (6, 0) {\huge$\sigma$};
	\node[circle, fill=cyan!70,inner sep=1.5pt] (h14) at (6, -1.5) {\huge$\sigma$};
	\node[circle, fill=cyan!70,inner sep=1.5pt] (h15) at (6, -3) {\large$\cdots$};
	\node[circle, fill=cyan!70,inner sep=2.5pt] (h16) at (6, -4.5) {\huge$\sigma$};	
	\node[] (input) at (6, -5.5) {Hidden layer};
	
	\draw[line width=1.0pt,->] (i0) -- (h10);
	\draw[line width=1.0pt,->] (i0) -- (h11);
	\draw[line width=1.0pt,->] (i0) -- (h12);
	\draw[line width=1.0pt,->] (i0) -- (h13);
	\draw[line width=1.0pt,->] (i0) -- (h14);
	\draw[line width=1.0pt,->] (i0) -- (h15);
	\draw[line width=1.0pt,->] (i0) -- (h16);
	
	\draw[line width=1.0pt,->] (i1) -- (h10);
	\draw[line width=1.0pt,->] (i1) -- (h11);
	\draw[line width=1.0pt,->] (i1) -- (h12);
	\draw[line width=1.0pt,->] (i1) -- (h13);
	\draw[line width=1.0pt,->] (i1) -- (h14);
	\draw[line width=1.0pt,->] (i1) -- (h15);
	\draw[line width=1.0pt,->] (i1) -- (h16);
	
	\draw[line width=1.0pt,->] (i2) -- (h10);
	\draw[line width=1.0pt,->] (i2) -- (h11);
	\draw[line width=1.0pt,->] (i2) -- (h12);
	\draw[line width=1.0pt,->] (i2) -- (h13);
	\draw[line width=1.0pt,->] (i2) -- (h14);
	\draw[line width=1.0pt,->] (i2) -- (h15);
	\draw[line width=1.0pt,->] (i2) -- (h16);
	
	\draw[line width=1.0pt,->] (i3) -- (h10);
	\draw[line width=1.0pt,->] (i3) -- (h11);
	\draw[line width=1.0pt,->] (i3) -- (h12);
	\draw[line width=1.0pt,->] (i3) -- (h13);
	\draw[line width=1.0pt,->] (i3) -- (h14);
	\draw[line width=1.0pt,->] (i3) -- (h15);
	\draw[line width=1.0pt,->] (i3) -- (h16);;
	
	\draw[line width=1.0pt,->] (i4) -- (h10);
	\draw[line width=1.0pt,->] (i4) -- (h11);
	\draw[line width=1.0pt,->] (i4) -- (h12);
	\draw[line width=1.0pt,->] (i4) -- (h13);
	\draw[line width=1.0pt,->] (i4) -- (h14);
	\draw[line width=1.0pt,->] (i4) -- (h15);
	\draw[line width=1.0pt,->] (i4) -- (h16);
	
	\node[circle, fill=yellow!100,inner sep=5.5pt] (out) at (9, 0) {};	
	
	\draw[line width=1.0pt,->] (h10) -- (out);
	\draw[line width=1.0pt,->] (h11) -- (out);
	\draw[line width=1.0pt,->] (h12) -- (out);
	\draw[line width=1.0pt,->] (h13) -- (out);
	\draw[line width=1.0pt,->] (h14) -- (out);
	\draw[line width=1.0pt,->] (h15) -- (out);
	\draw[line width=1.0pt,->] (h16) -- (out);
	
	\node[circle] (y0) at (11, 0) {};
	\node[] (output) at (9.5, -5.5) {Out layer};
	
	\draw[line width=1.0pt,->] (out) -- node[above]{$y$}(y0);	
	\end{tikzpicture}
	\caption{ \small Basic structure of ELM.}
	\label{fig2elm}
\end{figure}

Incorporating physical laws into the ELM model introduces a novel approach called PIELM. PIELM optimizes its loss function, which integrates information from governing terms as well as enforced boundary or initial conditions, by analytically calculating the weights of the outer layer. Unlike PINN, which are commonly used for solving PDEs, PIELM addresses several limitations associated with PINN, such as restricted accuracy, uncertain convergence, and tremendous computational demand that can lead to slow training speeds. Similar to previous ELM applications in classification and fitting tasks, PIELM is data-adaptive, fast convergent, and outperforms the gradient descent-based methods. Additionally, PIELM preserves the core advantages of ELM as a universal approximator for any continuous function;~\citet{huang2006extreme} rigorously proved that for any infinitely differentiable activation function, SLFNs with $N$ hidden nodes can exactly learn $N$ distinct samples. This implies that if a learning error ($\epsilon$) is permissible, SLFNs require fewer than $N$ hidden nodes.

\section{Fourier induced PIELM architecture for biharmonic equations}\label{sec:03}

In this section, the unified architecture of FPIELM is proposed to solve the biharmonic equation \eqref{eq:biharmonic} with boundaries \eqref{bd:Dirichlet} or \eqref{bd:Navier} by embracing the vanilla PIELM method with Fourier theorem.

\subsection{Vanilla PIELM framework to solve biharmonic equations}\label{algor2CELM}

This section provides a detailed introduction to the PIELM method for solving the biharmonic equation \eqref{eq:biharmonic} with boundaries \eqref{bd:Dirichlet} or \eqref{bd:Navier}. To obtain an ideal solution to the biharmonic equation, one can seek the solution in trial function space spanned by the ELM model with a high-order smooth activation function. Then, an ansatz expressed by ELM is designed to approximate the solutions of \eqref{eq:biharmonic}, it is:  
\begin{equation*}
  u(\bm{x}; \bm{W},\bm{b})= \sigma(\bm{x}^T\cdot\bm{W}^T+\bm{b}^T)\cdot \bm{\beta}.
\end{equation*}
Substituting $u(\bm{x}; \bm{W},\bm{b})$  into the equation of \eqref{eq:biharmonic} for $u$, we can obtain the following equation
\begin{equation}\label{elm2solu}
\displaystyle \Delta^2\left(\sum_{i=1}^{N} \beta_i\sigma\big{(}V_i(\bm{x})\big{)}\right) = f(\bm{x})~~\text{for}~~\bm{x}\in \Omega,
\end{equation}
with $V_i(\bm{x})=\bm{x}^T\cdot\bm{W}^{[i]}+b_i$ being the linear output for $i_{th}$ hidden unit of ELM and $N$ is the number of hidden unit. Additionally,  $u(\bm{x}; \bm{W},\bm{b})$ still holds on the prescribed boundaries \eqref{bd:Dirichlet} or \eqref{bd:Navier},  we then have the following formulation
\begin{equation}\label{elm2ubd_Dirchlet}
\begin{cases}
\displaystyle \sum_{i=1}^{N} \beta_i\sigma\big{(}V_i(\bm{x})\big{)} = g(\bm{x}),\\
\displaystyle \partial \bigg{(} \sum_{i=1}^{N} \beta_i\sigma\big{(}V_i(\bm{x})\big{)} \bigg{)}\bigg{/}\partial \vec{n}= h(\bm{x}),
\end{cases}
~~\text{for}~~\bm{x}\in \partial \Omega.
\end{equation}
or
\begin{equation}\label{elm2ubd_Navier}
\begin{cases}
\displaystyle \sum_{i=1}^{N} \beta_i\sigma\big{(}V_i(\bm{x})\big{)} = g(\bm{x}),\\
\displaystyle \Delta \bigg{(} \sum_{i=1}^{N} \beta_i\sigma\big{(}V_i(\bm{x})\big{)} \bigg{)}= k(\bm{x}),
\end{cases}
~~\text{for}~~\bm{x}\in \partial \Omega.
\end{equation}

For solving the biharmonic equations by the aforementioned PIELM method, a collection of collocation points is sampled from both the interior and boundaries of the domain. At these collocation points, the governed function \eqref{elm2solu} and the \eqref{elm2ubd_Dirchlet} or \eqref{elm2ubd_Navier} exhibit coerciveness.
Let $Q$ denote the number of random collocation points in the interior of $\Omega$ and denote by $\mathcal{S}_I=\{\bm{x}_I^q\}_{q=1}^Q$, and $P$ denote the number of random collocation points on hyperface of $\partial\Omega$ and denote by $\mathcal{S}_B=\{\bm{x}_B^p\}_{p=1}^P$. Feeding the above collocation points into the pipeline of PIELM, we have the following residual terms for Dirichlet boundaries:
\begin{equation}\label{residual}
   \begin{cases}
	   \mathcal{R}_u(\bm{x}_I^q) = \displaystyle \Delta^2\left(\sum_{i=1}^{N} \beta_i\sigma\big{(}V_i(\bm{x}_I^q)\big{)}\right) - f(\bm{x}_I^q),\\
	   \mathcal{L}_g(\bm{x}_B^p) = \displaystyle \sum_{i=1}^{N} \beta_i\sigma\big{(}V_i(\bm{x}_B^p)\big{)} - h(\bm{x}_B^p)\\
	    \mathcal{L}_h(\bm{x}_B^p) =\displaystyle \partial \bigg{(} \sum_{i=1}^{N} \beta_i\sigma\big{(}V_i(\bm{x}_B^p)\big{)} \bigg{)}\bigg{/}\partial \vec{n}- h(\bm{x}_B^p)
   \end{cases}
\end{equation}
and the last term in \eqref{residual} will be replaced by 
\begin{equation}
     \mathcal{L}_k(\bm{x}_B^p) = \displaystyle \Delta \bigg{(} \sum_{i=1}^{N} \beta_i\sigma\big{(}V_i(\bm{x}_B^p)\big{)} \bigg{)} - k(\bm{x}_B^p),
\end{equation}
for Navier boundaries.

To approximate  the functions $u$ well, we need to adjust the output weights of the PIELM model to make each residual of \eqref{residual} to near-zero or minimize the error of the following system of linear equations:
\begin{equation}\label{opt2linear_sys}
\mathbf{H}\bm{\beta}=\mathbf{S} \Longleftrightarrow \left[
\begin{matrix}
\mathbf{H}_{f}\\
\mathbf{H}_{g}\\
\mathbf{H}_{h} 
\end{matrix}
\right]\bm{\beta}=
\left[
\begin{matrix}
\mathbf{S}_f\\ 
\mathbf{S}_g\\
\mathbf{S}_h 
\end{matrix}
\right]~~\text{or}~~
\left[
\begin{matrix}
\mathbf{H}_{f}\\
\mathbf{H}_{g}\\
\mathbf{H}_{k} 
\end{matrix}
\right]\bm{\beta}=
\left[
\begin{matrix}
\mathbf{S}_f\\ 
\mathbf{S}_g\\
\mathbf{S}_k 
\end{matrix}
\right].
\end{equation}
Here 
$\displaystyle \mathbf{H}_{f}=\mathbf{H}_{B}+\sum_{i=1}^{d}\sum_{j=1,j\neq i}^{d}\mathbf{H}_{ij}$ with
\begin{equation*}
\mathbf{H}_{B}=\left[
\begin{matrix}
\sigma''''\left(V_1(\bm{x}^1_I)\right) & \sigma''''\left(V_2(\bm{x}^1_I)\right)& \cdots &\sigma''''\left(V_N(\bm{x}^1_I)\right)\\
\sigma''''\left(V_1(\bm{x}^2_I)\right) & \sigma''''\left(V_2(\bm{x}^2_I)\right)& \cdots &\sigma''''\left(V_N(\bm{x}^2_I)\right)\\
\vdots&\vdots&\vdots&\vdots\\
\sigma''''\left(V_1(\bm{x}^Q_I)\right) & \sigma''''\left(V_2(\bm{x}^Q_I)\right)& \cdots &\sigma''''\left(V_N(\bm{x}^Q_I)\right)
\end{matrix}
\right]\bigodot
\left[
\begin{matrix}
\displaystyle \sum_{i=1}^{d}w^4_{i1}&\displaystyle \sum_{i=1}^{d}w^4_{i2} & \cdots &\displaystyle\sum_{i=1}^{d}w^4_{iN}
\end{matrix}
\right]
\end{equation*}
and
\begin{equation*}
\mathbf{H}_{ij}=\left[
\begin{matrix}
\sigma''''\left(V_1(\bm{x}^1_I)\right) & \sigma''''\left(V_2(\bm{x}^1_I)\right)& \cdots &\sigma''''\left(V_N(\bm{x}^1_I)\right)\\
\sigma''''\left(V_1(\bm{x}^2_I)\right) & \sigma''''\left(V_2(\bm{x}^2_I)\right)& \cdots &\sigma''''\left(V_N(\bm{x}^2_I)\right)\\
\vdots&\vdots&\vdots&\vdots\\
\sigma''''\left(V_1(\bm{x}^Q_I)\right) & \sigma''''\left(V_2(\bm{x}^Q_I)\right)& \cdots &\sigma''''\left(V_N(\bm{x}^Q_I)\right)
\end{matrix}
\right]\bigodot
\bigg{[}
\begin{matrix}
\displaystyle w^2_{j1}*w^2_{i1}&\displaystyle w^2_{j2}*w^2_{i2} & \cdots &\displaystyle w^2_{jN}*w^2_{iN}
\end{matrix}
\bigg{]}.
\end{equation*}
Additionally, the matrices $\mathbf{H}_{g}$, $\mathbf{H}_{h}$ and $\mathbf{H}_{k}$  are given by
\begin{equation*}
\mathbf{H}_{g} = \left[
\begin{matrix}
\sigma\left(V_1(\bm{x}^1_B)\right) & \sigma\left(V_2(\bm{x}_B^1)\right)& \cdots&\sigma\left(V_N(\bm{x}_B^1)\right)\\
\sigma\left(V_1(\bm{x}^2_B)\right) & \sigma\left(V_2(\bm{x}_B^2)\right)& \cdots&\sigma\left(V_N(\bm{x}_B^2)\right)\\
\vdots&\vdots&\vdots&\vdots\\
\sigma\left(V_1(\bm{x}^P_B)\right) & \sigma\left(V_2(\bm{x}^P_B)\right)& \cdots&\sigma\left(V_N(\bm{x}^P_B)\right)
\end{matrix}
\right],
\end{equation*}
$\displaystyle\mathbf{H}_{h} =\mathbf{H}_1\cdot n_1+ \mathbf{H}_2\cdot n_2+\cdots+\mathbf{H}_d\cdot n_d$ with 
\begin{equation*}
\mathbf{H}_{i} = \left[
\begin{matrix}
\sigma'\left(V_1(\bm{x}^1_B)\right) & \sigma'\left(V_2(\bm{x}_B^1)\right)& \cdots&\sigma'\left(V_N(\bm{x}_B^1)\right)\\
\sigma'\left(V_1(\bm{x}^2_B)\right) & \sigma'\left(V_2(\bm{x}_B^2)\right)& \cdots&\sigma'\left(V_N(\bm{x}_B^2)\right)\\
\vdots&\vdots&\vdots&\vdots\\
\sigma'\left(V_1(\bm{x}^P_B)\right) & \sigma'\left(V_2(\bm{x}^P_B)\right)& \cdots&\sigma'\left(V_N(\bm{x}^P_B)\right)
\end{matrix}
\right]\bigodot
\bigg{[}
\begin{matrix}
\displaystyle w_{i1}&\displaystyle w_{i2} & \ldots  &\displaystyle w_{iN}
\end{matrix}
\bigg{]},~~ i=1,2,\cdots, d,
\end{equation*}
and 
\begin{equation*}
\mathbf{H}_{k} = \left[
\begin{matrix}
\sigma''\left(V_1(\bm{x}^1_B)\right) & \sigma''\left(V_2(\bm{x}_B^1)\right)& \cdots&\sigma''\left(V_N(\bm{x}_B^1)\right)\\
\sigma''\left(V_1(\bm{x}^2_B)\right) & \sigma''\left(V_2(\bm{x}_B^2)\right)& \cdots&\sigma''\left(V_N(\bm{x}_B^2)\right)\\
\vdots&\vdots&\vdots&\vdots\\
\sigma''\left(V_1(\bm{x}^P_B)\right) & \sigma''\left(V_2(\bm{x}^P_B)\right)& \cdots&\sigma''\left(V_N(\bm{x}^P_B)\right)
\end{matrix}
\right]\bigodot
\left[
\begin{matrix}
\displaystyle \sum_{i=1}^{d}w^2_{i1}&\displaystyle \sum_{i=1}^{d}w^2_{i2} & \cdots &\displaystyle\sum_{i=1}^{d}w^2_{iN}
\end{matrix}
\right].
\end{equation*}
In the above expressions,  $\sigma''$ and $\sigma''''$ are the second and fourth-order derivatives of activation function $\sigma$, respectively, and $\bigodot$ stands for the Hadamard product operator. In addition, $\mathbf{S}_f=\big{(}f(\bm{x}_I^1),f(\bm{x}_I^2),\cdots, f(\bm{x}_I^Q)\big{)}^T$,  $\mathbf{S}_g=\big{(}g(\bm{x}_B^1),g(\bm{x}_B^2),\cdots, g(\bm{x}_B^P)\big{)}^T$, $\mathbf{S}_h=\big{(}h(\bm{x}_B^1),h(\bm{x}_B^2),\cdots, h(\bm{x}_B^P)\big{)}^T$ and $\mathbf{S}_k=\big{(}k(\bm{x}_B^1),k(\bm{x}_B^2),\cdots, k(\bm{x}_B^P)\big{)}^T$.

When the weights and biases of the hidden layer for the PIELM model are randomly initialized and fixed, the matrix $\mathbf{H}$ is determined for given input data, then the output weights $\bm{\beta}$ of PIELM model are obtained by solving the linear systems~\eqref{opt2linear_sys}. The system equation $\mathbf{H}\bm{\beta} = \mathbf{S}$ is solvable and meets the following several cases:
\begin{itemize}
    \item If the coefficient matrix $\mathbf{H}$ is square and invertible, then 
    \begin{equation}
      \bm{\beta} = \mathbf{H}^{-1}\mathbf{S}.
    \end{equation}
    \item  If the coefficient matrix $\mathbf{H}$ is rectangular and full rank, then 
    \begin{equation}
      \bm{\beta} = \mathbf{H}^{\dagger}\mathbf{S},
    \end{equation}
    where $\mathbf{H}^{\dagger}$ is the pseudo inverse of $\mathbf{H}$ defned by $\mathbf{H}^{\dagger}=(\mathbf{H}^T\mathbf{H})^{-1}\mathbf{H}^T$.
    \item If the coefficient matrix $\mathbf{H}$ is  singular, then the Tikhonov regularization\cite{adriazola2022role} can be utilized to to overcome the ill-posedness of \eqref{opt2linear_sys}, and the $\bm{\beta}$ is given by
    \begin{equation}	
    	\bm{\beta} = (\mathbf{H}^T\mathbf{H}+\lambda \mathbf{I})^{-1}\mathbf{H}\mathbf{S},
    \end{equation}
    where $\mathbf{I}$ is the identity matrix, and $\lambda>0$ is the ridge parameter.
\end{itemize}

According to the above discussions, the detailed procedures involved in the PIELM algorithm to address the biharmonic equation \eqref{eq:biharmonic} within finite-dimensional spaces are described in the following.

 \begin{algorithm}[H]
    \caption{PIELM algorithm for solving biharmonic equation~\eqref{eq:biharmonic}.}
    \begin{itemize}
        \item [1)] Generating the points set $\mathcal{S}^k$ includes interior points $S_{I}=\{\bm{x}^q_I\}_{i=1}^{Q}$ with $\bm{x}_I^q\in\mathbb{R}^d$ and boundary points $S_{B}=\{\bm{x}^p_B\}_{j=1}^{P}$ with $\bm{x}_B^j\in\mathbb{R}^d$. Here, we draw the random points $\bm{x}_I^q $ and $\bm{x}_B^p$ from $\mathbb{R}^d$ with positive probability density $\nu$, such as uniform distribution.
        \item [2)] Selecting a shallow neural network and randomly initialize weights and biases of the hidden layer, then fixing the parameters of the hidden layer.
        \item [3)] Embedding the input data into the implicit feature space by a linear transformation, then activating nonlinearly the results and casting them into the output layer.
        \item [4)] Combining linearly the all output features of the hidden layer and output this results.
        \item [5)] Computing the derivation about the input variable for the output of the FELM model and constructing the linear system~\eqref{opt2linear_sys}.
        \item [6)] Learning the output layer weights by the least-squares method.
    \end{itemize}
    \label{algor:FMPINN}
\end{algorithm}

\subsection{The option of activation function}

In terms of function approximation, the hidden nodes governed by a given activation function in the ELM model can be considered a set of basis functions, and then the ELM model produces its output through a linear aggregation of these fundamental functions. Similar to other artificial neural network (ANN) architectures, the choice of activation function is critical for the ELM model’s performance. Since the biharmonic equations \eqref{eq:biharmonic} studied here is a high-order PDE, activation functions with strong regularity are preferred. Previous studies, including~\citet{fabiani2021numerical}, \citet{calabro2021extreme}, and \citet{ma2022novel}, have investigated the influence of sigmoidal functions and radial basis functions in the ELM model for solving both linear and nonlinear PDEs. These studies indicate that ELM models with sigmoid, Gaussian, and tangent activation functions achieve commendable numerical accuracy.

\begin{remark}
	(\textbf{Lipschitz continuity})~If a function $\sigma$ is continuous (i.e., $\sigma\in C^1$) and satisfies the following boundedness condition:
	\begin{equation*}
	|\sigma(x)|<1~~~~\textup{and}~~~~|\sigma'(x)|<1
	\end{equation*}
	for any $x\in\mathbb{R}$. Then, we have
	\begin{equation*}
	|\sigma(x)-\sigma(y)|\leqslant |x-y|~~~~\textup{and}~~~~|\sigma'(x)-\sigma'(y)|<|x-y|
	\end{equation*}
	for any $x,y\in\mathbb{R}$. 
\end{remark}

\begin{figure}[H]
    \centering
    \subfigure[Sigmoid fucntion]{
        \label{Sigmoid} 
        \includegraphics[scale=0.22]{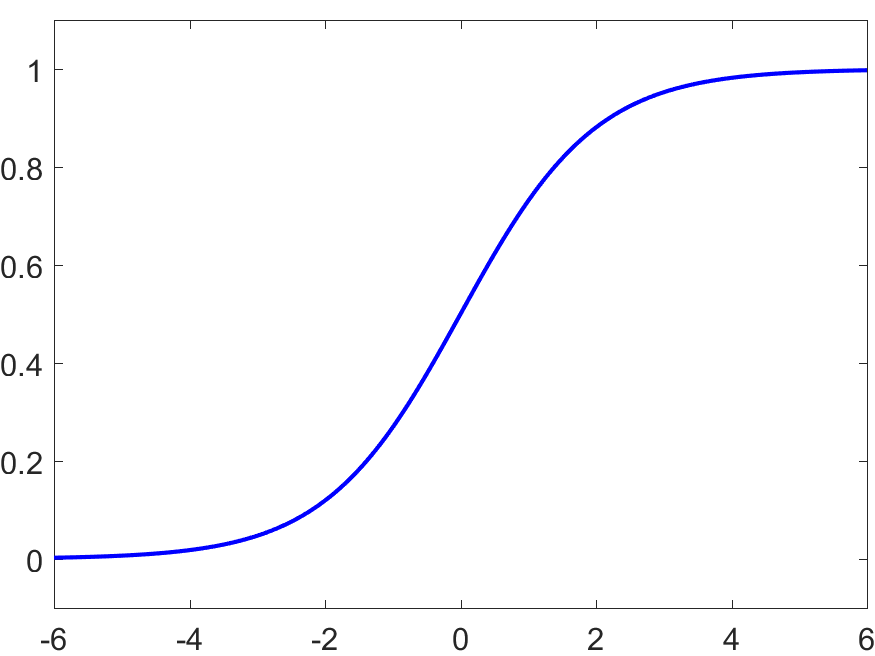}
    }
    \subfigure[1st-order derivative]{
        \label{1stDeri2Sigmoid}
        \includegraphics[scale=0.22]{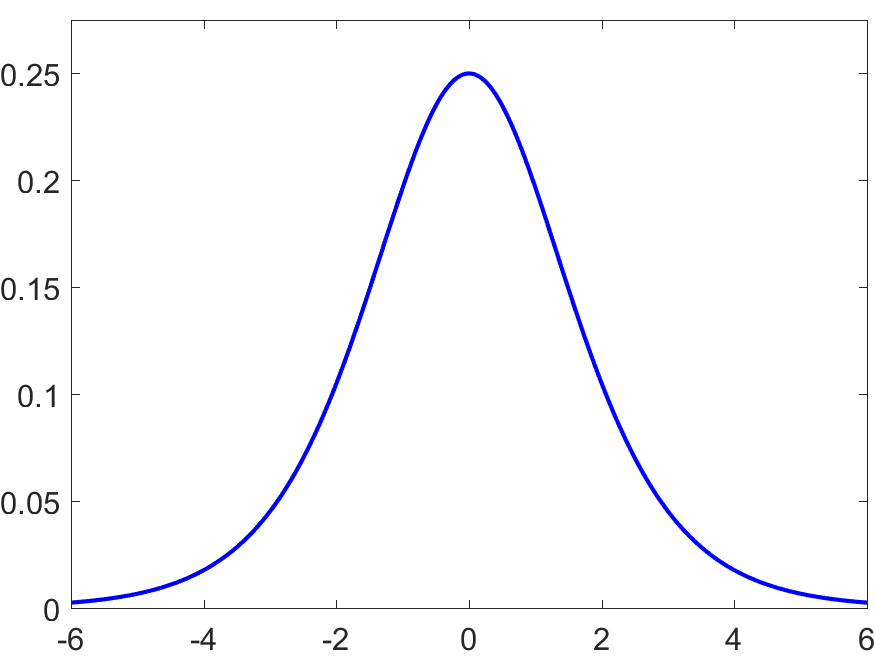}
    }
    \subfigure[2nd-order derivative]{
    \label{2ndDeri2Sigmoid}
    \includegraphics[scale=0.22]{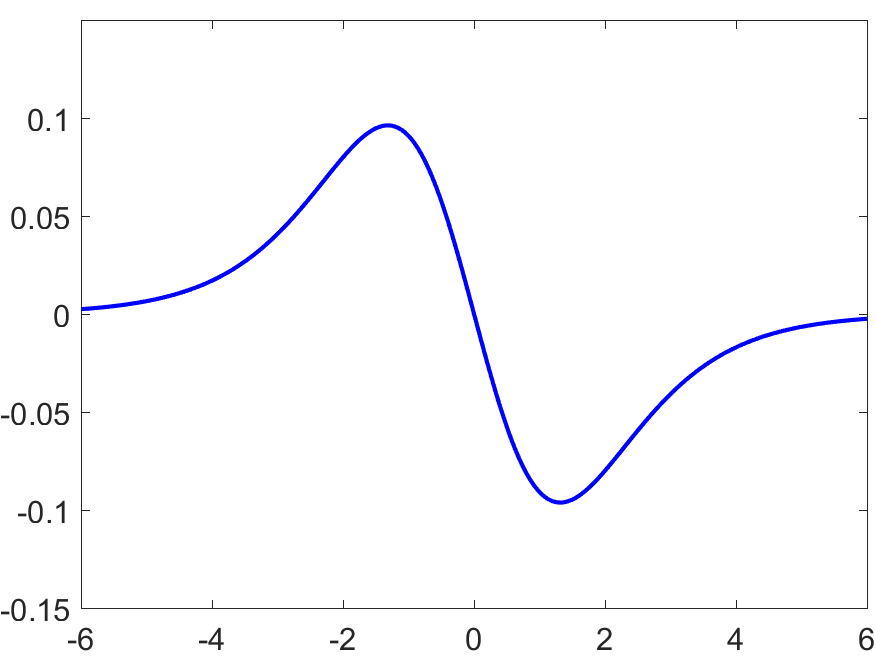}
    }
    \subfigure[4th-order derivative]{
    \label{4thDeri2Sigmoid}
    \includegraphics[scale=0.22]{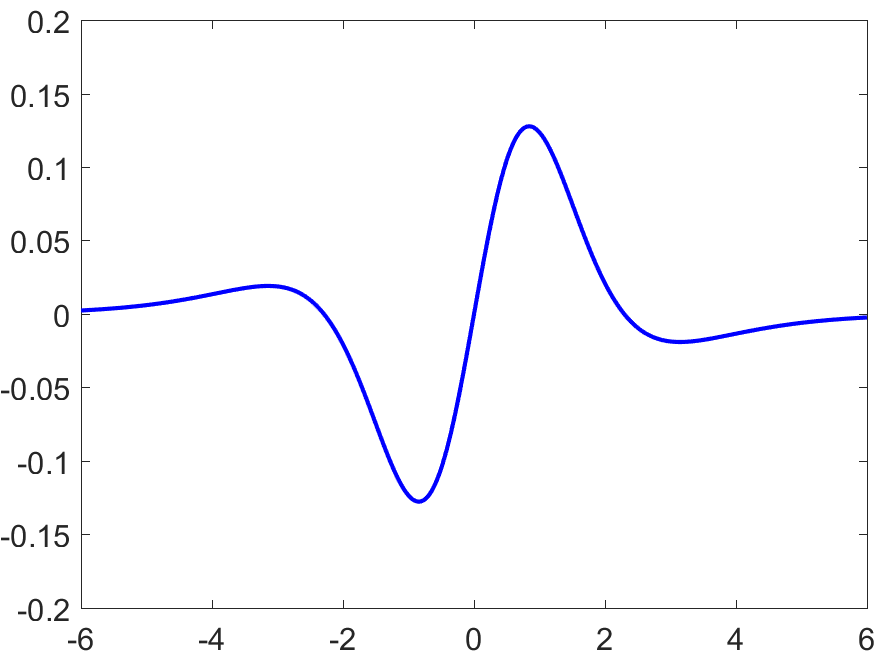}
    }
    \caption{\small The sigmoid function and its 1st-order, 2nd-order and 4th-order derivatives, respectively.}
    \label{sigmoid_act_func}
\end{figure}

\begin{figure}[H]
    \centering
    \subfigure[Gaussian fucntion]{
        \label{Gaussian} 
        \includegraphics[scale=0.22]{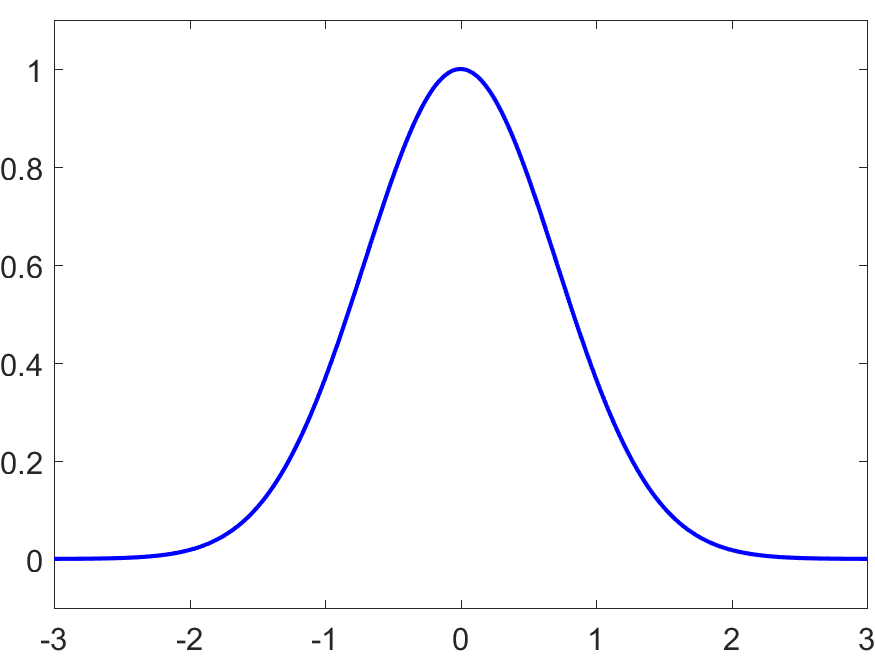}
    }
    \subfigure[1st-order derivative]{
        \label{1stDeri2Gaussian}
        \includegraphics[scale=0.22]{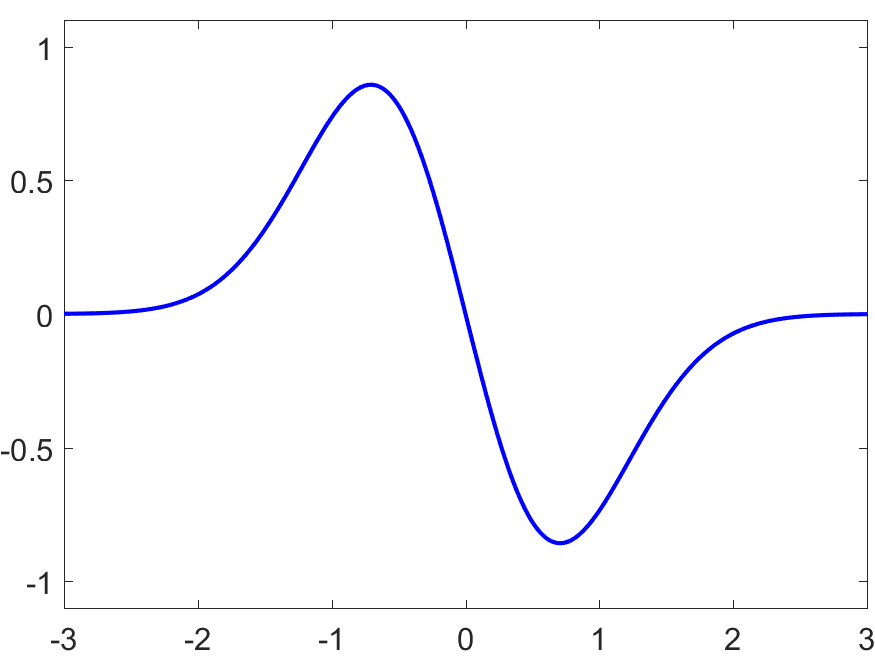}
    }
    \subfigure[2nd-order derivative]{
    \label{2ndDeri2Gaussian}
    \includegraphics[scale=0.22]{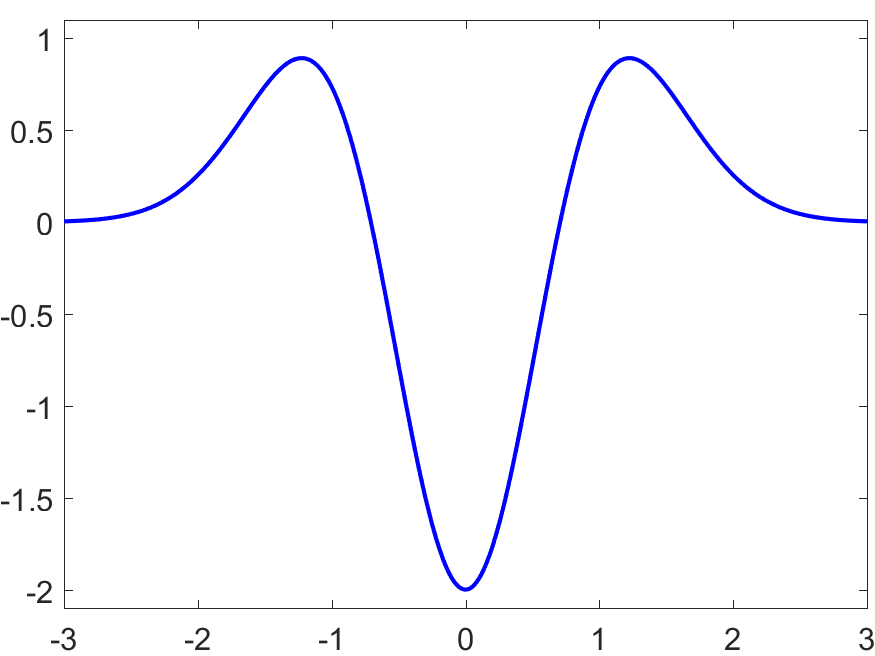}
    }
    \subfigure[4th-order derivative]{
    \label{4thDeri2Gaussian}
    \includegraphics[scale=0.22]{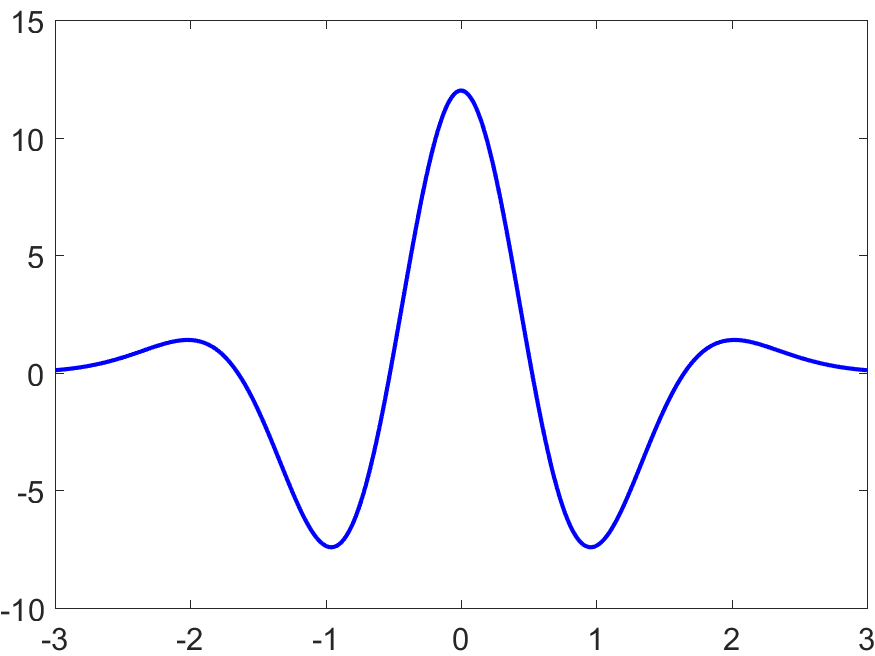}
    }
    \caption{\small The Gaussian function and its 1st-order, 2nd-order and 4th-order derivatives, respectively.}
    \label{gauss_act_func}
\end{figure}

\begin{figure}[H]
    \centering
    \subfigure[Tanh fucntion]{
        \label{Tanh} 
        \includegraphics[scale=0.22]{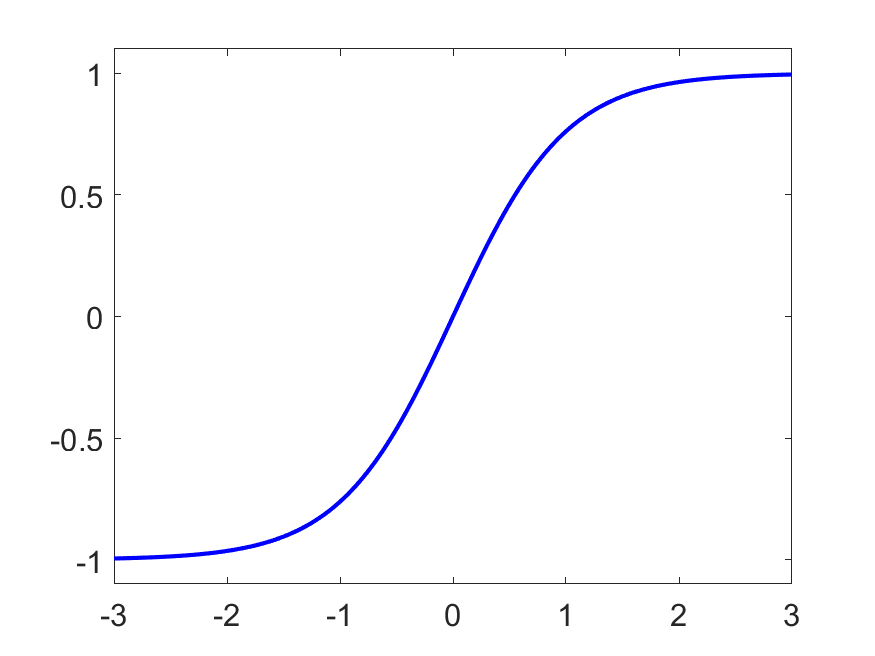}
    }
    \subfigure[1st-order derivative]{
        \label{1stDeri2tanh}
        \includegraphics[scale=0.22]{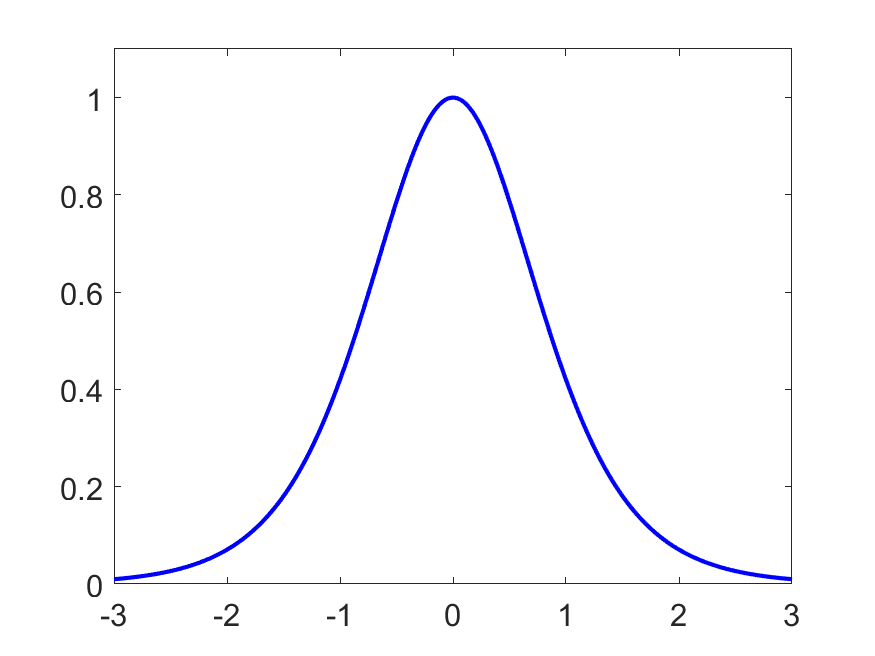}
    }
    \subfigure[2nd-order derivative]{
    \label{2ndDeri2tanh}
    \includegraphics[scale=0.22]{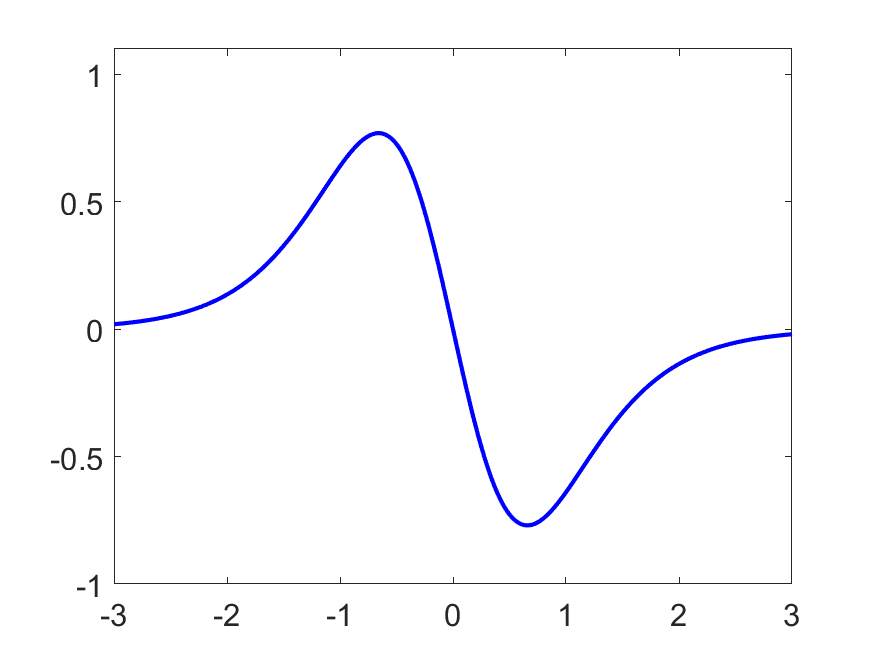}
    }
    \subfigure[4th-order derivative]{
    \label{4thDeri2tanh}
    \includegraphics[scale=0.22]{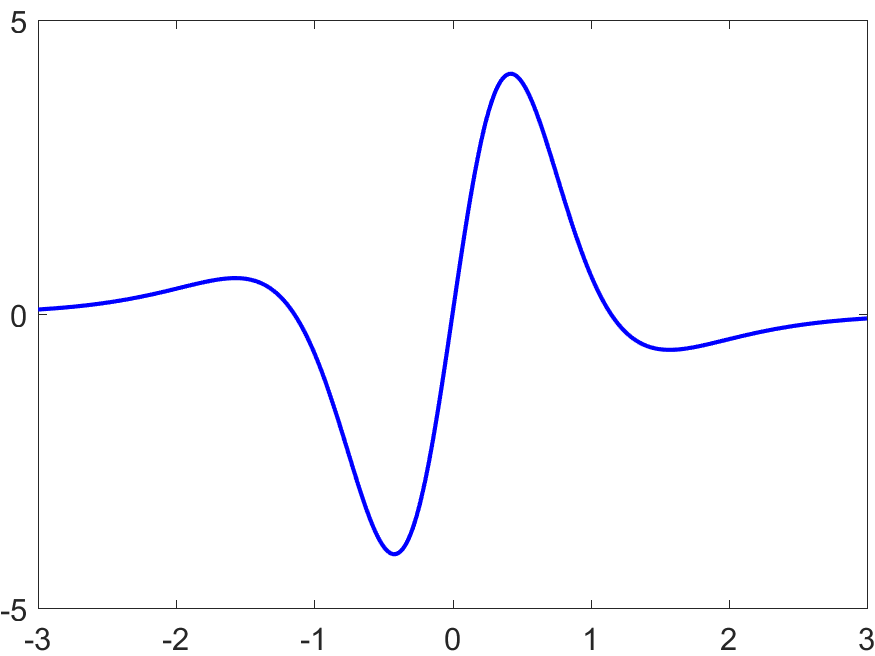}
    }
    \caption{\small The tanh function and its 1st-order, 2nd-order and 4th-order derivatives, respectively.}
    \label{tanh_act_func}
\end{figure}

The sigmoid, gaussian, and tangent functions all fulfill the Lipschitz continuity and exhibit robust regularity, they ensure the boundedness of the output for hidden units in ELM. In terms of the high-order PDEs, the high-order derivate of activation is required. Then, the first, second, and fourth-order derivatives for the sigmoidal function $s(x)=1/(1+e^{-x})$, Gaussian function $\sigma(x)=e^{-x^2}$ and hyperbolic tangent function $tanh(x) = (e^x-e^{-x})/(e^x+e^{-x})$ are given by
\begin{equation}\label{multi_deri2sigmoid}
\left\{
    \begin{aligned}
        s'(x) &= s(x) - s^2(x)\\
        s''(x) &= s(x)-3s^2(x)+2s^3(x)\\
        s''''(x) &= s(x) - 15s^2(x)+50s^3(x)-60s^4(x)+24s^5(x)
    \end{aligned}\right.,
\end{equation}
\begin{equation}\label{multi_deri2gaussian}
\left\{
\begin{aligned}
\Big{(}\sigma(x)\Big{)}' &= -2xe^{-x^2}\\
\Big{(}\sigma(x)\Big{)}''&= -2e^{-x^2}+4x^2e^{-x^2}\\
\Big{(}\sigma(x)\Big{)}''''&= 12e^{-x^2}-48x^2e^{-x^2}+16x^4e^{-x^2}
\end{aligned}\right.,
\end{equation}
and
\begin{equation}\label{multi_deri2tan}
\left\{
\begin{aligned}
\tanh'(x) &= 1-\tanh^2(x)\\
\tanh''(x) &= -2\tanh(x)+2\tanh^3(x)\\
\tanh''''(x) &= \Big{(}16\tanh(x)-24\tanh^3(x)\Big{)}\Big{(}1-\tanh^2(x)\Big{)}
\end{aligned}\right.,
\end{equation}
respectively. Figs.~\ref{sigmoid_act_func} -- \ref{tanh_act_func} illustrate the curves of the sigmoid, Gaussian, and tangent functions and their corresponding first-order, second-order and fourth-order derivatives, respectively. 

This visualization shows that the derivatives of the sigmoid, Gaussian, and tanh functions are bounded and tend to saturate for narrow-range inputs. However, their performance may degrade with large-range inputs. Inspired by the Fourier-based spectral method, a global basis function can generate responses, along with its derivatives, for both narrow- and large-range inputs. Furthermore, approximations based on such basis functions can handle data across various scales. The sine function serves as an ideal global Fourier basis function (with sine and cosine being fundamentally equivalent). In this study, we adopt the sine function as the activation function in the PIELM method, referring to this modified approach as FPIELM.

\section{Numerical experiments}\label{sec:05}

\subsection{Experimental setup}

In this section, we apply the proposed FPIELM method to solve the two- and three-dimensional biharmonic equation \eqref{eq:biharmonic} in both unitized and non-unitized domains. Additionally, PIELM models using Sigmoid, Gaussian, and Tanh as activation functions are established as baselines for solving \eqref{eq:biharmonic} and are referred to as SPIELM, GPIELM, and TPIELM, respectively. Following the results in~\citet{dong2022computing}, the weight matrices and bias vectors for the hidden layers of each model are initialized uniformly within the interval $[-\delta, \delta]$.

A criterion is provided to evaluate the performance of the four aforementioned methods, given by
\begin{equation*}
\text{REL} = \sqrt{\frac{\sum_{i=1}^{N'}|\tilde{u}(\bm{x}^i)-u^*(\bm{x}^i)|^2}{\sum_{i=1}^{N'}|u^*(\bm{x}^i)|^2}},
\end{equation*}
where $\tilde{u}(\bm{x}^i)$ and $u^*(\bm{x}^i)$ denote the approximate solution obtained by numerical methods and the exact solution, respectively, for the collection set $\{\bm{x}^i\}_{i=1}^{N'}$, with $N'$ being the number of collection points.

All experiments are implemented in Numpy (version 1.21.0) on a workstation equipped with 256 GB RAM and a single NVIDIA GeForce RTX 4090 Ti GPU with 24 GB of memory.

\subsection{Numerical examples for Dirichlet boundary}
\begin{example}\label{2D_Dirichlet_E1}
	To start, we aim to approximate the solution to the biharmonic equation \eqref{eq:biharmonic} with Dirichlet boundary in the regular domain $\Omega=[x_1^{min}, x_1^{max}]\times[x_2^{min}, x_2^{max}]$. An exact solution is given by
	\begin{equation*}
	u(x_1,x_2)=[(x_1-x_1^{min})(x_1^{max}-x_1)]^2[(x_2-x_2^{min})(x_2^{max}-x_2)]^2
	\end{equation*}
	such that 
	\begin{equation*}
	\begin{aligned}
	f(x_1, x_2) &= 24(x_2-x_2^{min})^2(x_2^{max}-x_2)^2+24(x_1-x_1^{min})^2(x_1^{max}-x_1)^2 \\&+2[2(x_1^{min}+x^{max}_1-2x_1)^2-4(x_1-x_1^{min})(x_1^{max}-x_1)]\\ &\times[2(x_2^{min}+x_2^{max}-2x_2)^2-4(x_2-x_2^{min})(x_2^{max}-x_2)].
	\end{aligned}
	\end{equation*}
	The boundary functions $g(x_1,x_2)$ and $h(x_1,x_2)$ can be obtained easily by direct computation.
\end{example}

First, we approximate the solution of \eqref{eq:biharmonic} using the four methods on 16,384 equidistant grid points distributed across regular domains: $[-1,1]\times[-1,1]$, $[0,5]\times[0,5]$, $[-4,6]\times[-3,7]$, and $[5,15]\times[0,10]$. In addition to the grid points, the training set includes 10,000 random points within $\Omega$ and 4,000 random points on $\partial\Omega$. The number of hidden nodes for each of the four PIELM methods is set to 1,000.

\begin{figure}[H]
	\centering
	\subfigure[Analytical solution]{
		\label{2D_Dirichlet_E1:Exact} 
		\includegraphics[scale=0.325]{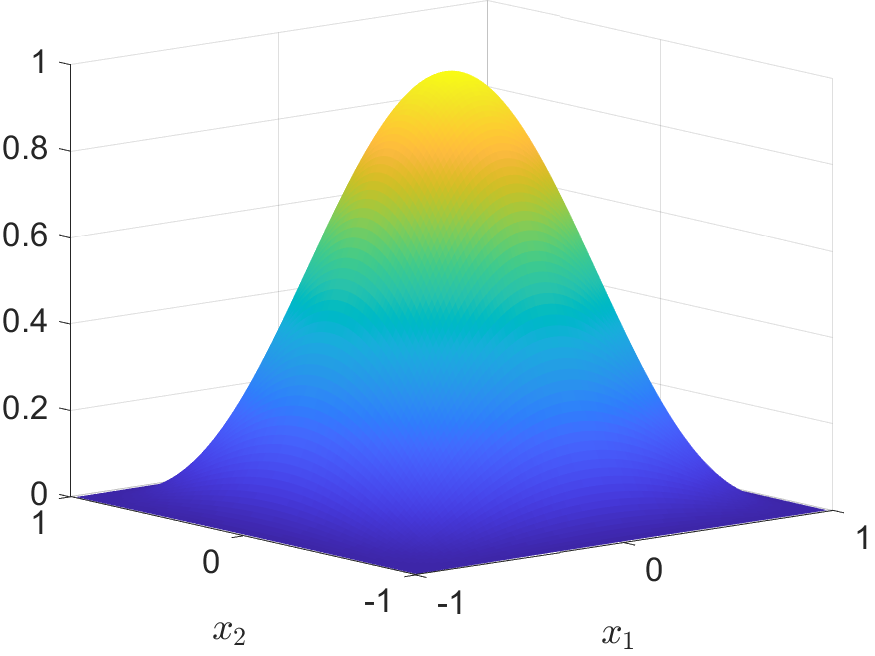}
	}
	\subfigure[SPIELM]{
		\label{2D_Dirichlet_E1:Perr2Simoid}
		\includegraphics[scale=0.325]{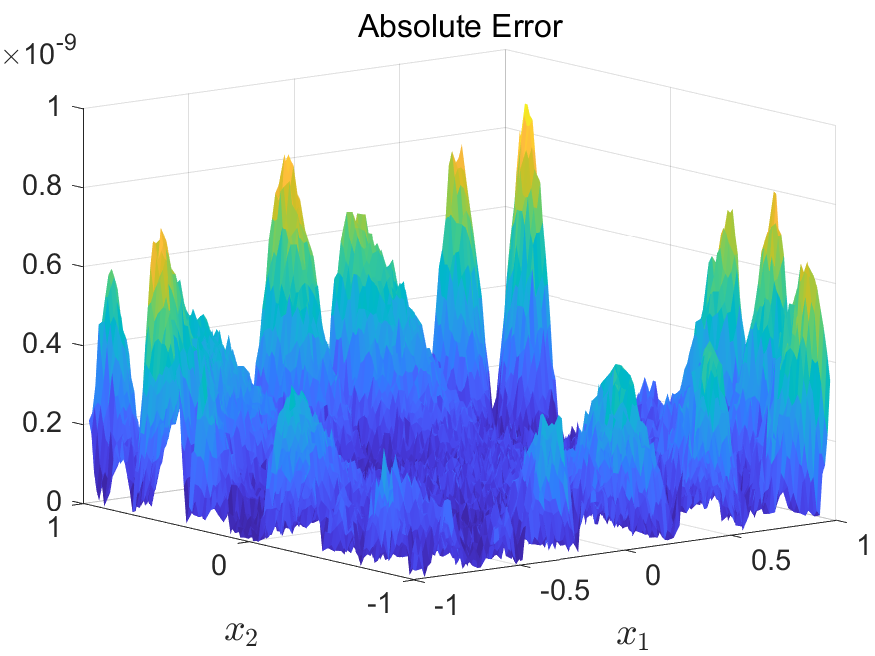}
	}
	\subfigure[GPIELM]{
		\label{2D_Dirichlet_E1:Perr2Gauss}
		\includegraphics[scale=0.325]{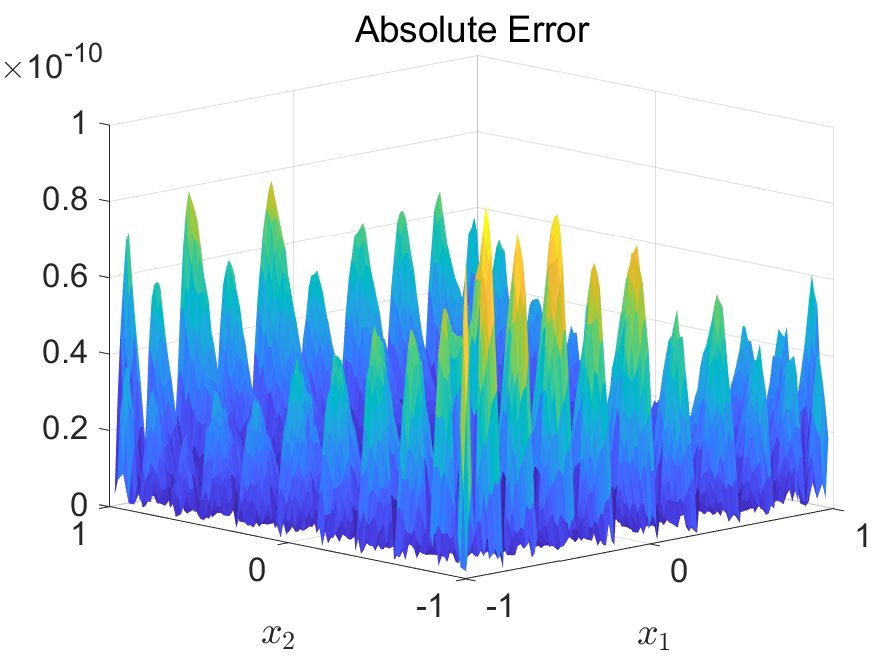}
	}
	\subfigure[TPIELM]{
		\label{2D_Dirichlet_E1:Perr2Tanh}
		\includegraphics[scale=0.35]{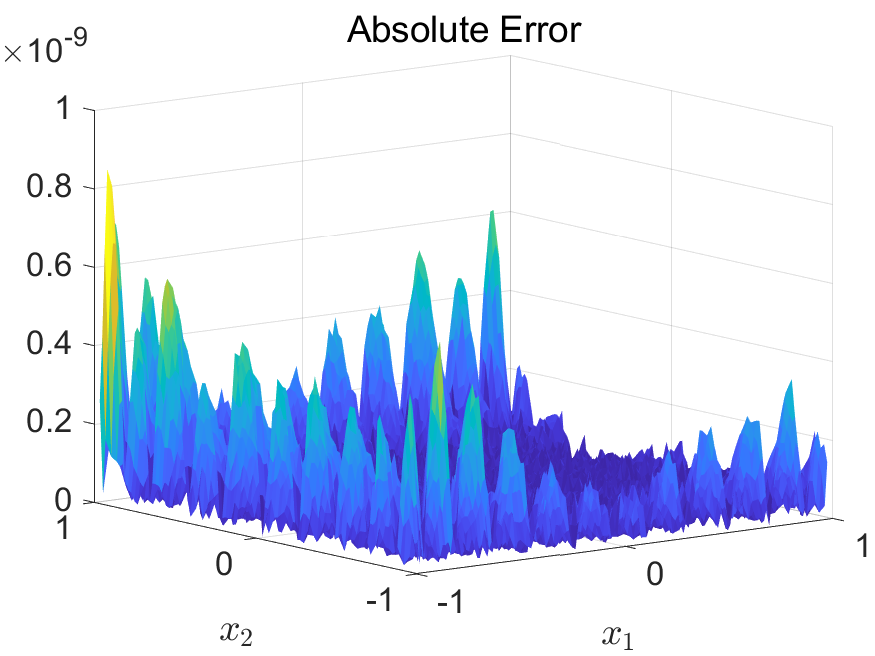}
	}
	\subfigure[FPIELM]{
		\label{2D_Dirichlet_E1:Perr2Sine}
		\includegraphics[scale=0.35]{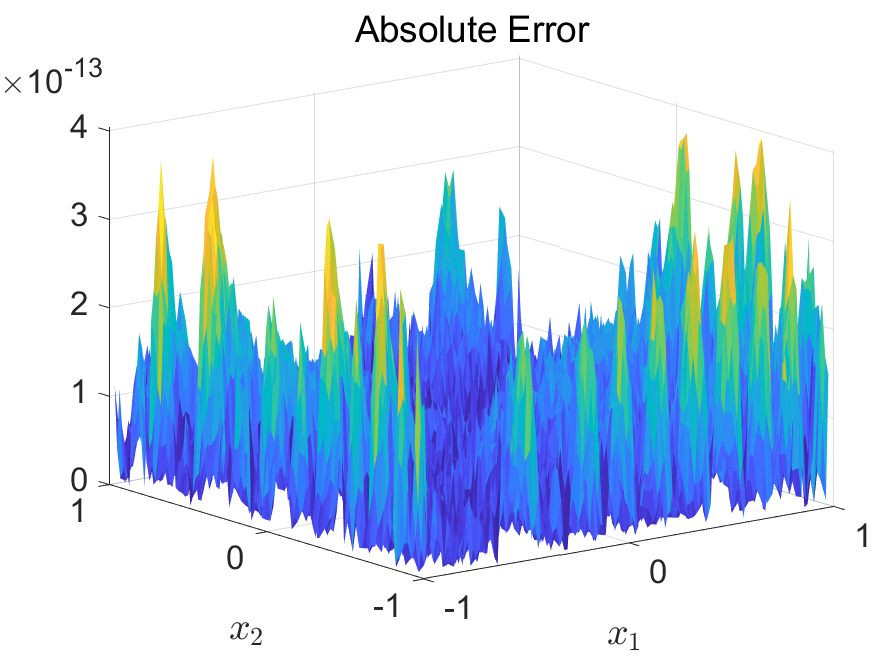}
	}
	\caption{Numerical results of four PIELM methods for solving Example~\ref{2D_Dirichlet_E1} on domain $[-1,1]\times[-1,1]$.}
	\label{Plot_2D_Dirichlet_E1}
\end{figure}

Figure~\ref{Plot_2D_Dirichlet_E1} shows the point-wise absolute error of the four PIELM models on the domain $[-1,1]\times[-1,1]$. Table~\ref{Table_2D_Dirichlet_E1} presents the optimal value of $\delta$, the REL, and the runtime for each PIELM method across the different domains. 

\begin{figure}[H]
	\centering
	\subfigure[REL VS the hidden node]{
		\label{2D_Dirichlet_E1:rels} 
		\includegraphics[scale=0.325]{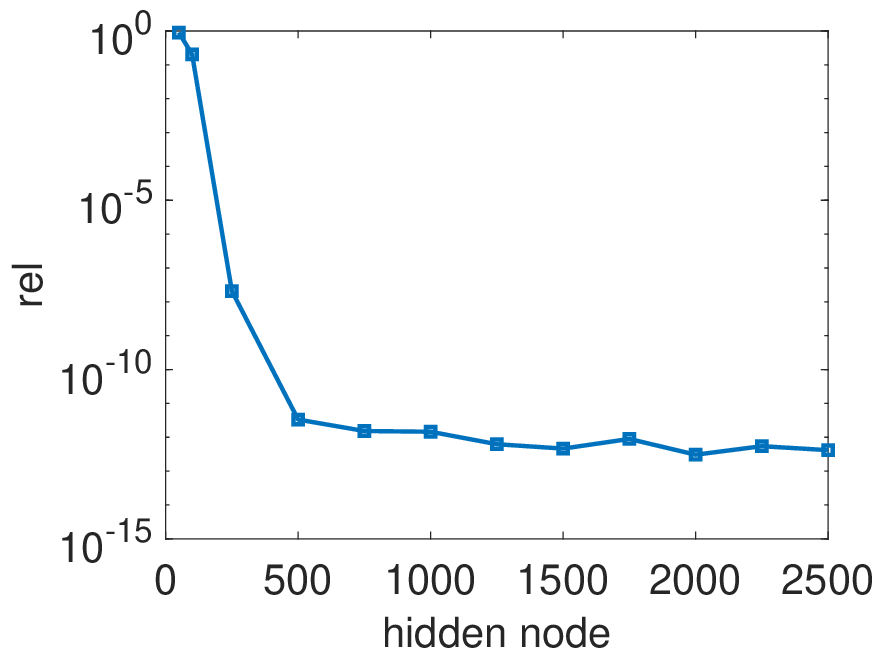}
	}
	\subfigure[Running time VS the hidden node]{
		\label{2D_Dirichlet_E1:times}
		\includegraphics[scale=0.325]{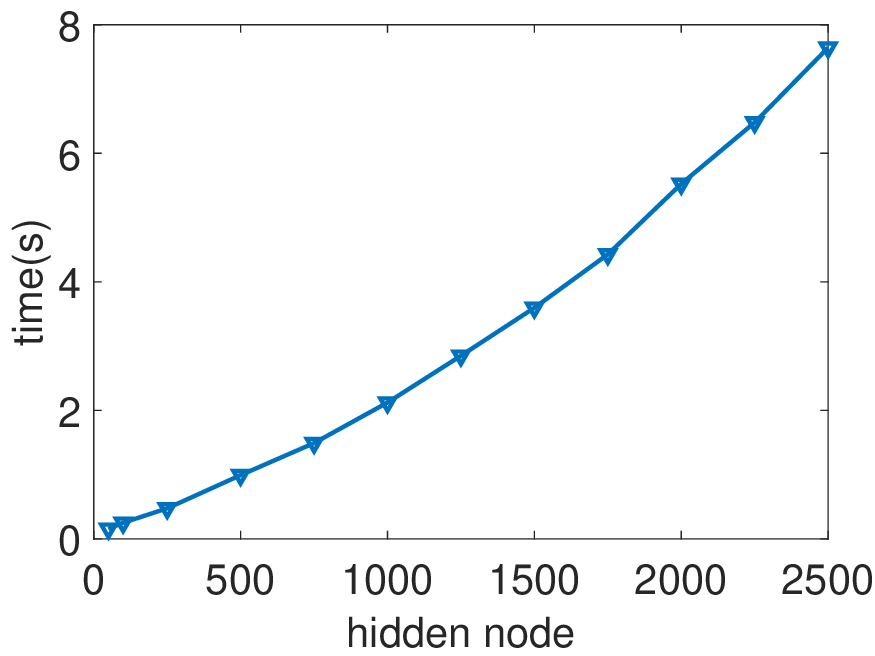}
	}
	\subfigure[REL VS the scale factor $\delta$]{
		\label{2D_Dirichlet_E1:sigma}
		\includegraphics[scale=0.325]{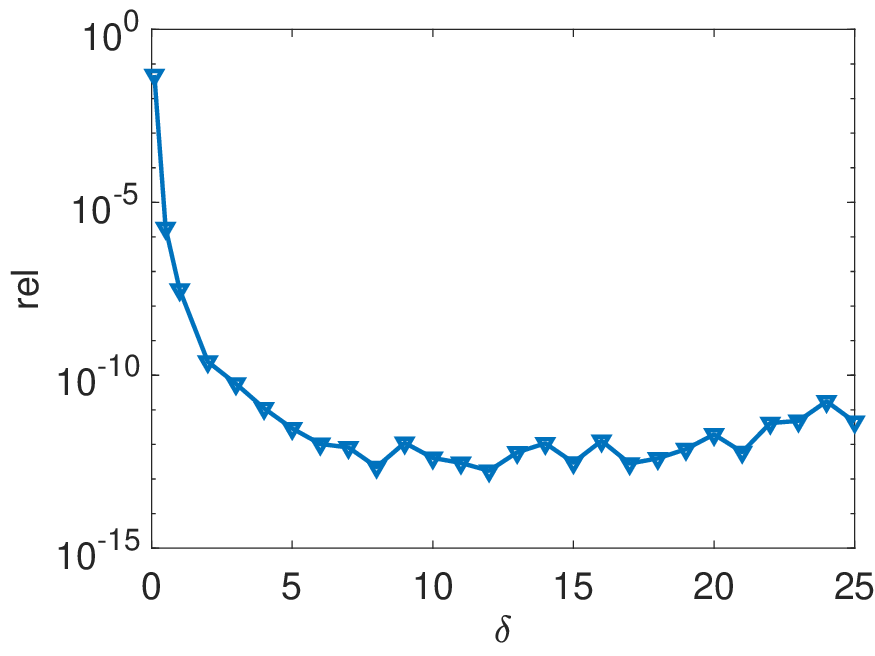}
	}
	\caption{REL VS the hidden node (left), running time VS the hidden node (middle) and REL VS the scale factor $\delta$ (right) for FPIELM method to Example \ref{2D_Dirichlet_E1} on domain$[-1,1]\times[-1,1]$.}
	\label{Plot_2D_Diri_E1_rel_time_Sigma}
\end{figure}

\begin{figure}[H]
	\centering
	\subfigure[]{
		\label{2D_Dirichlet_E1:1111}
		\includegraphics[scale=0.45]{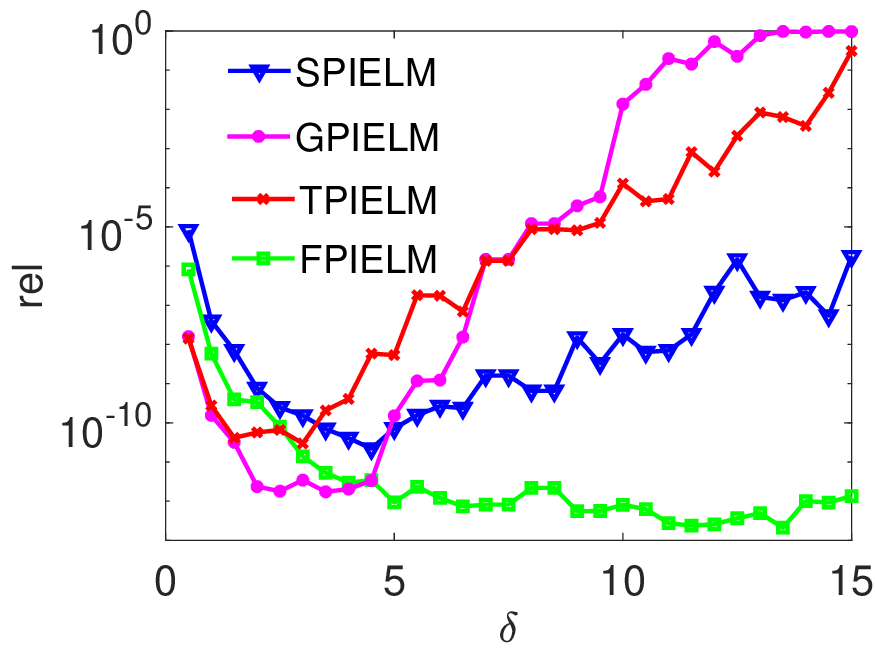}
	}
	\subfigure[]{
		\label{2D_Dirichlet_E1:0505}
		\includegraphics[scale=0.45]{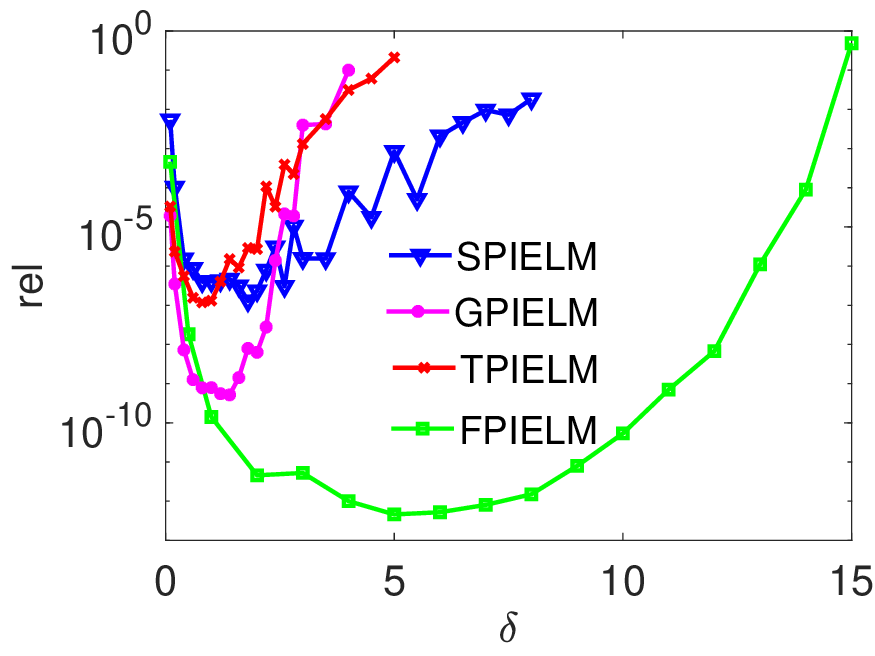}
	}
	\caption{REL VS the scale factor $\delta$ for the four PIELM models to solve Example \ref{2D_Dirichlet_E1} on regular domain $[-1,1]\times[-1,1]$(left) and $[0,5]\times[0,5]$(right).}
	\label{Plot_2D_Diri_E1_rel05051111}
\end{figure}

Next, we examine the impact of the scale factor $\delta$ and the number of hidden nodes on the FPIELM model in the domain $[-1,1]\times[-1,1]$, with results illustrated in Fig.~\ref{Plot_2D_Diri_E1_rel_time_Sigma}. Finally, comparison curves for REL as a function of the scale factor $\delta$ for the four PIELM methods on the domains $[-1,1]\times[-1,1]$ and $[0,5]\times[0,5]$ are shown in Fig.~\ref{Plot_2D_Diri_E1_rel05051111}.

\begin{table}[!ht]
	\centering
	\caption{REL, $\delta$ and runtime for four PIELM methods to Example \ref{2D_Dirichlet_E1} on various domain}
	\label{Table_2D_Dirichlet_E1}
	\begin{tabular}{|l|c|c|c|c|c|}
		\hline  
		Domain               &          & SPIELM               & GPIELM               &  TPIELM              &  FPIELM    \\    \hline
		&$\delta$  &6.0                   &1.5                   &1.4                   &8.0           \\
		$[-1,1]\times[-1,1]$ & REL      &$4.799\times 10^{-10}$&$4.059\times 10^{-11}$ &$1.396\times 10^{-10}$ &$2.223\times 10^{-13}$ \\ 
		& Time(s)  &8.203                      &3.939                      &3.265                      &2.875    \\  \hline
		&$\delta$  &2.0                   &0.8                  &1.4                   &5.0           \\ 
		$[0,5]\times[0,5]$   & REL      &$3.358\times 10^{-7}$&$5.795\times 10^{-10}$ &$1.819\times 10^{-7}$ &$2.573\times 10^{-13}$ \\ 
		& Time(s)  &6.676                 &3.081                      &3.129                &2.152    \\  \hline
		&$\delta$  &0.5                   &0.3                   &0.2                   &1.0           \\ 
		$[-4,6]\times[-3,7]$ & REL      &$4.029\times 10^{-7}$&$2.436\times 10^{-8}$ &$2.955\times 10^{-7}$ &$4.935\times 10^{-11}$ \\ 
		& Time(s)  &8.077                 &4.046                 &3.328                &3.031    \\  \hline
		&$\delta$  &0.1                   &0.2                   &0.12                  &0.8           \\ 
		$[5,15]\times[0,10]$ & REL      &$3.503\times 10^{-3}$&$1.187\times 10^{-4}$ &$9.331\times 10^{-4}$   &$6.981\times 10^{-8}$ \\
		& Time(s)  &8.138                & 3.891                     &3.251                      &2.843    \\  \hline
	\end{tabular}
\end{table}

Based on the point-wise absolute error distribution shown in Figures~\ref{2D_Dirichlet_E1:Perr2Simoid} -- \ref{2D_Dirichlet_E1:Perr2Sine} and the data in Table~\ref{Table_2D_Dirichlet_E1}, it can be concluded that the FPIELM method outperforms other types of PIELM methods. As observed in the curves in Figure~\ref{2D_Dirichlet_E1:rels}, the relative error does not decrease with an increasing number of hidden units, indicating that the performance of the FPIELM model stabilizes as the number of hidden neurons grows. The time-neuron curve in Figure~\ref{2D_Dirichlet_E1:times} demonstrates a superlinear increase in runtime with the addition of hidden neurons, though this increase is gradual. Furthermore, Figure~\ref{2D_Dirichlet_E1:sigma} illustrates that the FPIELM model remains stable and robust across variations in the scale factor $\sigma$ on the domain $[-1,1]\times[-1,1]$. The REL curves in Figure~\ref{Plot_2D_Diri_E1_rel05051111} confirm that FPIELM maintains stability and allows for a broader range of scale factor $\delta$ values across both unitized and non-unitized domains. In summary, our FPIELM method not only achieves high accuracy but also demonstrates efficiency and robustness. Additionally, it is noteworthy that the coefficient matrix in the final linear system is full rank, given that the number of samples exceeds the number of hidden nodes.

For comparison, we employ the finite difference method (FDM) with a five-point central difference scheme to solve the coupled scheme of \eqref{eq:biharmonic} on the domain $[0,5]\times[0,5]$. The FDM is implemented in MATLAB (version 12.0.0.341360, R2019a). For the FPIELM model, the number of hidden nodes is set to 1,000, and the scale factor $\delta$ for initializing weights and biases is set to 5.0. The numerical results in Table \ref{Comparision2FDM_FPIELM_Dirichlet_E1} indicate that the proposed FPIELM method maintains stability across various mesh grid sizes and significantly outperforms FDM in terms of accuracy, particularly on sparse grids.

\begin{table}[!ht]
	\centering
	\caption{REL and runtime of FDM and FPIELM for Example \ref{2D_Dirichlet_E1} on $[0,5]\times[0,5]$.}
	\label{Comparision2FDM_FPIELM_Dirichlet_E1}
	\begin{tabular}{|l|c|c|c|c|c|c|}
		\hline
		&number of mesh grid &$64\times 64$           &$128\times 128$        &$256\times 256$         &$512\times 512$     &$1024\times 1024$      \\  \hline
		FDM	&REL             &$1.65\times 10^{-3}$    &$4.14\times 10^{-4}$   &$1.04\times 10^{-4}$    &$2.59\times 10^{-5}$&$6.47\times 10^{-6}$   \\ 
		&Time(s)             &0.028                   &0.037                  &0.127                   &0.485               &2.05 \\ \hline
		FPIELM&REL            &$4.20\times 10^{-13}$  &$2.57\times 10^{-13}$   &$3.42\times 10^{-13}$   &$1.64\times 10^{-13}$&$2.43\times 10^{-13}$  \\  
		&Time(s)             &$1.47$                 &2.09                    &$4.71$                  &$16.67$  &$69.16$\\  \hline
	\end{tabular}
\end{table}

\begin{example}\label{2D_Dirichlet_E2}
We consider the biharmonic equation \eqref{eq:biharmonic} with Dirichlet boundary conditions in two-dimensional space and obtain its numerical solution within two hexagram-shaped domains, $\Omega$, derived from $[-\pi,\pi]\times[-\pi,\pi]$ and $[0,3\pi]\times[\pi,2\pi]$, respectively. An analytical solution is given by
\begin{equation*}
   u(x_1,x_2) = \sin(x_1)e^{\cos(x_2)},
\end{equation*}
such that
\begin{equation*}
	\begin{aligned}
	    f(x_1,x_2) &= \sin(x_1)e^{\cos(x_2)} - 2\sin(x_1)e^{\cos(x_2)}[\sin^2(x_2) - e^{\cos(x_2)}] \\
	    &\quad + \sin(x_1)e^{\cos(x_2)}[\sin^4(x_2) - 6\sin^2(x_2)\cos(x_2) + 3\cos^2(x_2) - 4\sin^2(x_2) + \cos(x_2)].
	\end{aligned}
\end{equation*}
The boundary constraint functions $g(x_1,x_2)$ and $h(x_1,x_2)$ on $\partial\Omega$ can be obtained by direct computation, but are omitted here for brevity.
\end{example}

\begin{figure}[H]
	\centering
	\subfigure[Analytical solution]{
		\label{2D_Dirichlet_E2:a} 
		\includegraphics[scale=0.325]{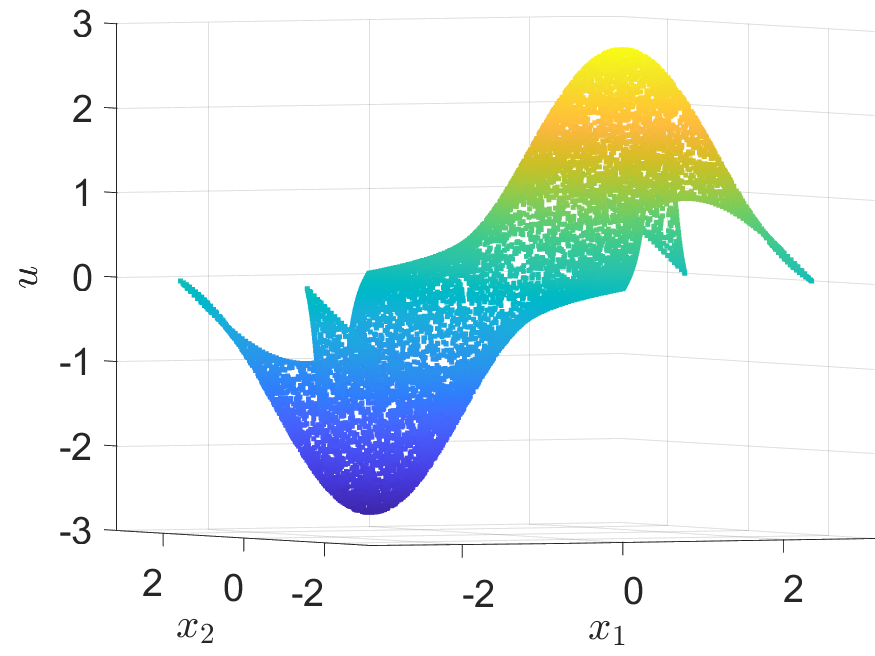}
	}
	\subfigure[Point-wise error of SPIELM]{
		\label{2D_Dirichlet_E2:Perr2Sigmoid}
		\includegraphics[scale=0.325]{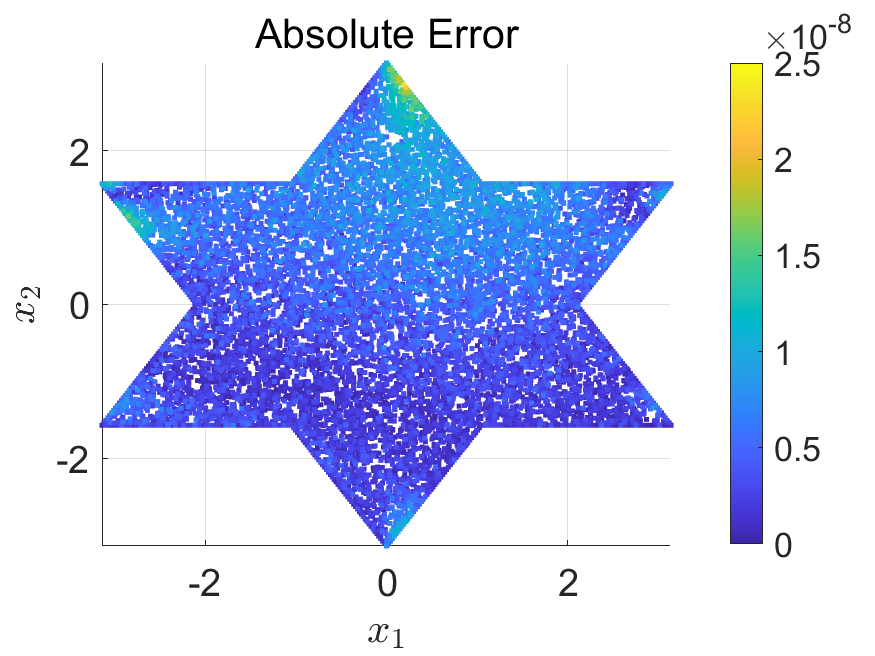}
	}
	\subfigure[Point-wise error of GPIELM]{
		\label{2D_Dirichlet_E2:Perr2Gauss}
		\includegraphics[scale=0.325]{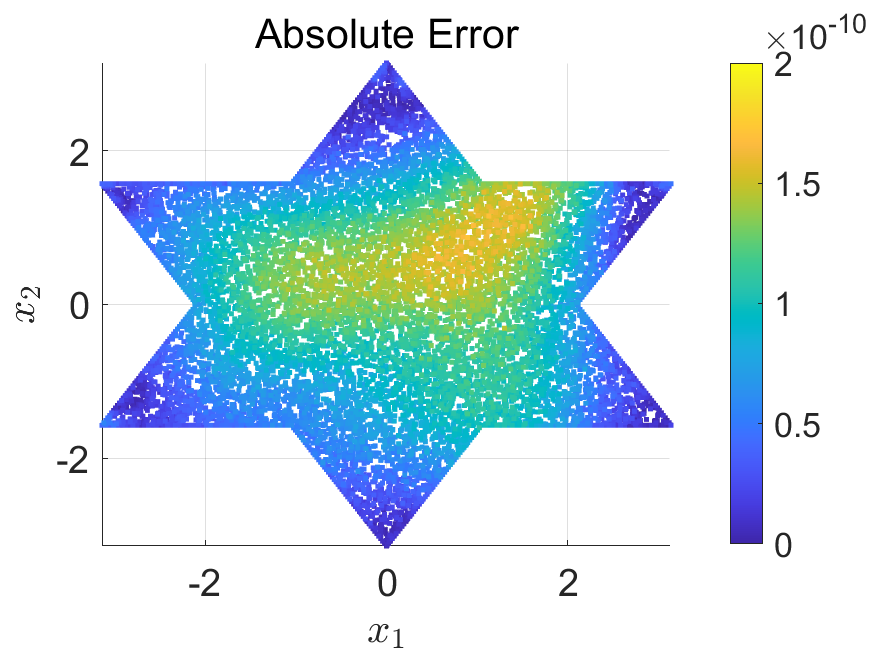}
	}
	\subfigure[Point-wise error of TPIELM]{
		\label{2D_Dirichlet_E2:Perr2Tanh}
		\includegraphics[scale=0.325]{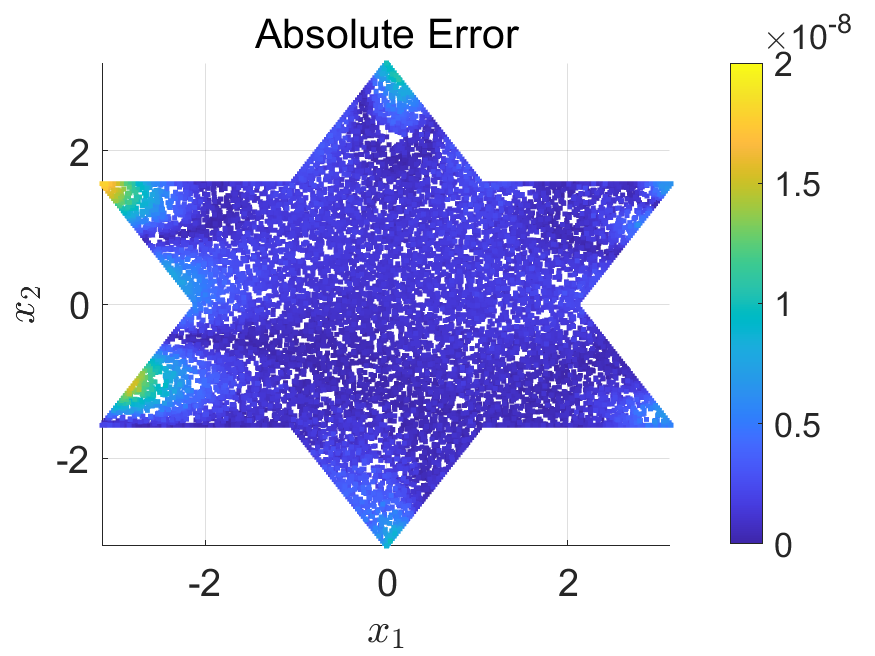}
	}
	\subfigure[Point-wise error of FPIELM]{
		\label{2D_Dirichlet_E2:Perr2SIN}
		\includegraphics[scale=0.325]{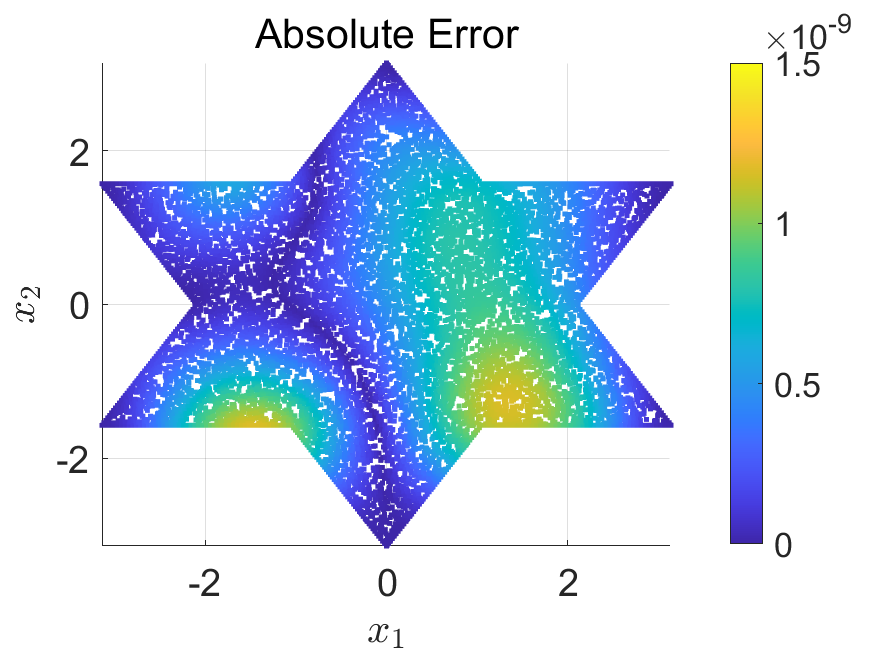}
	}
	\caption{Numerical results for four PIELM methods to Example~\ref{2D_Dirichlet_E2} on hexagram domian derived from $[-\pi,\pi]\times[-\pi,\pi]$.}
	\label{Plot_2D_Dirichlet_E2}
\end{figure}

In this example, the setup for the four PIELM methods is identical to that in Example~\ref{2D_Dirichlet_E1}. An approximate solution to \eqref{eq:biharmonic} is derived by applying the four methods on 19,214 collection points located within the hexagram domains $\Omega$. Numerical results are presented in Table~\ref{Table2D_Dirichlet_E2} for the two scenarios, and the point-wise absolute error for the four PIELM methods is shown in Figure~\ref{Plot_2D_Dirichlet_E2} over the hexagram domain derived from $[-\pi,\pi]\times[-\pi,\pi]$.

\begin{figure}[H]
    \centering
    \subfigure[REL VS the hidden node]{
        \label{2D_Dirichlet_E2:rels} 
        \includegraphics[scale=0.325]{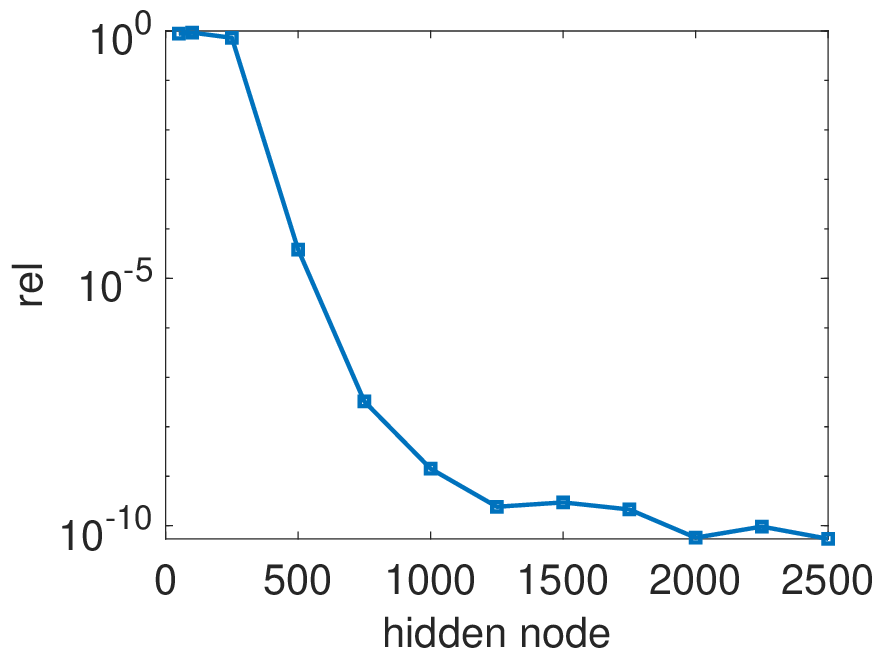}
    }
    \subfigure[Runtime VS the hidden node]{
        \label{2D_Dirichlet_E2:times}
        \includegraphics[scale=0.325]{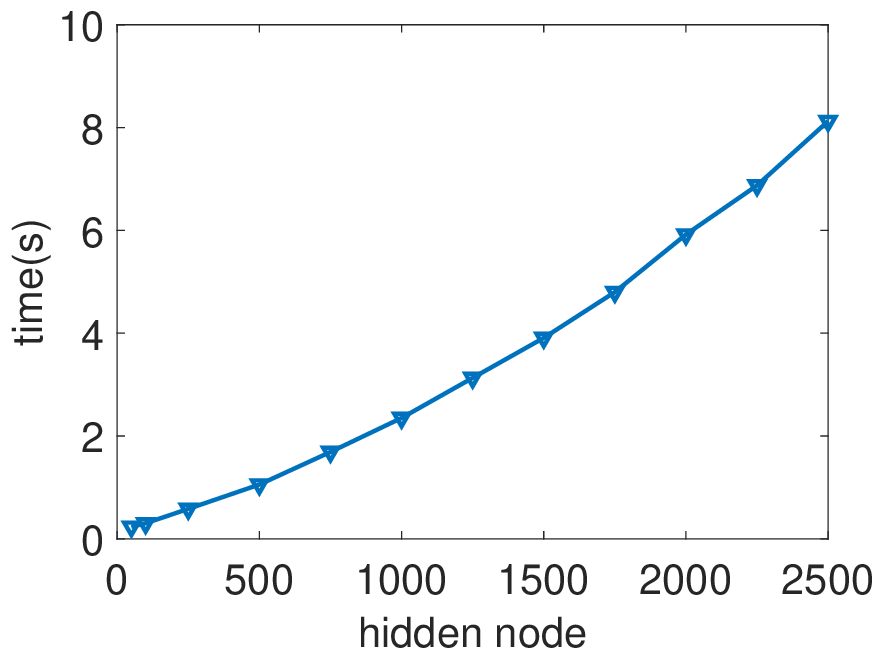}
    }
    \subfigure[REL VS the scale factor $\delta$]{
    \label{2D_Dirichlet_E2:sigma}
        \includegraphics[scale=0.325]{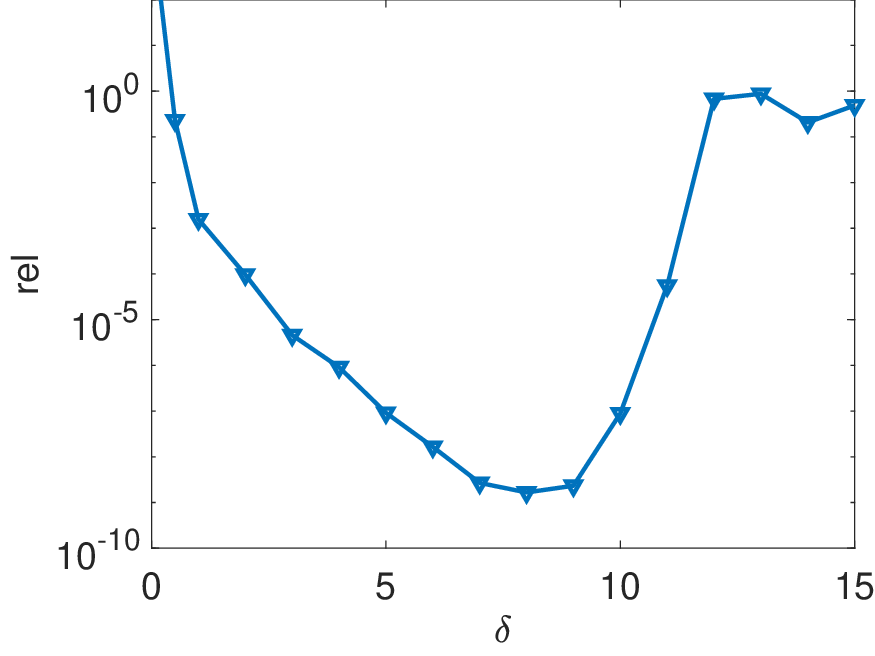}
    }
    \caption{REL VS the hidden node (left), runtime VS the hidden node (middle) and REL VS the scale factor $\delta$ (right) to FPIELM method to Example \ref{2D_Dirichlet_E2} on hexagram domain derived from $[-\pi,\pi]\times[-\pi,\pi]$.}
    \label{Plot_2D_rel_time_Sigma}
\end{figure}

\begin{figure}[H]
	\centering
	\subfigure[]{
		\includegraphics[scale=0.5]{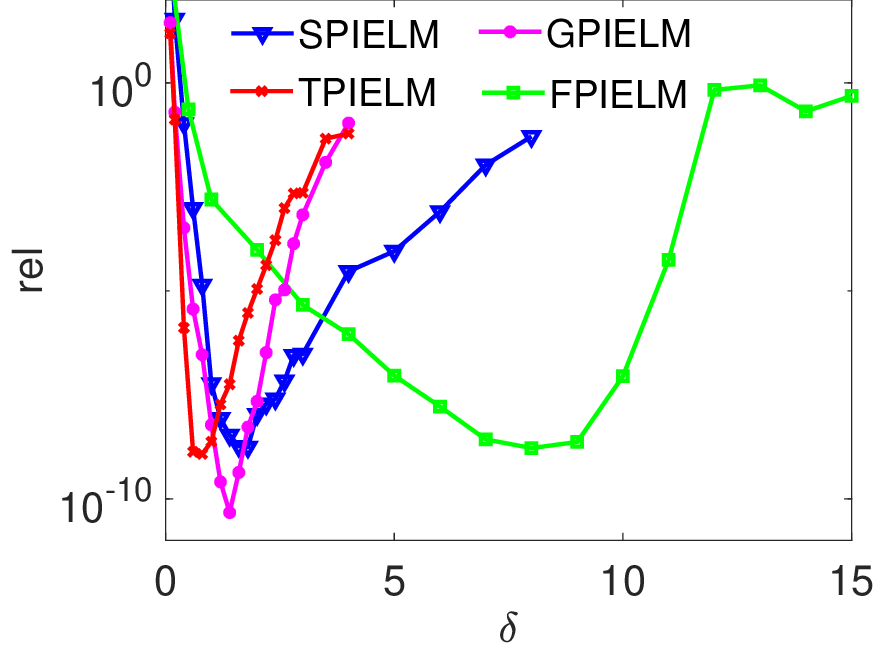}
	}
	\subfigure[]{
		\includegraphics[scale=0.475]{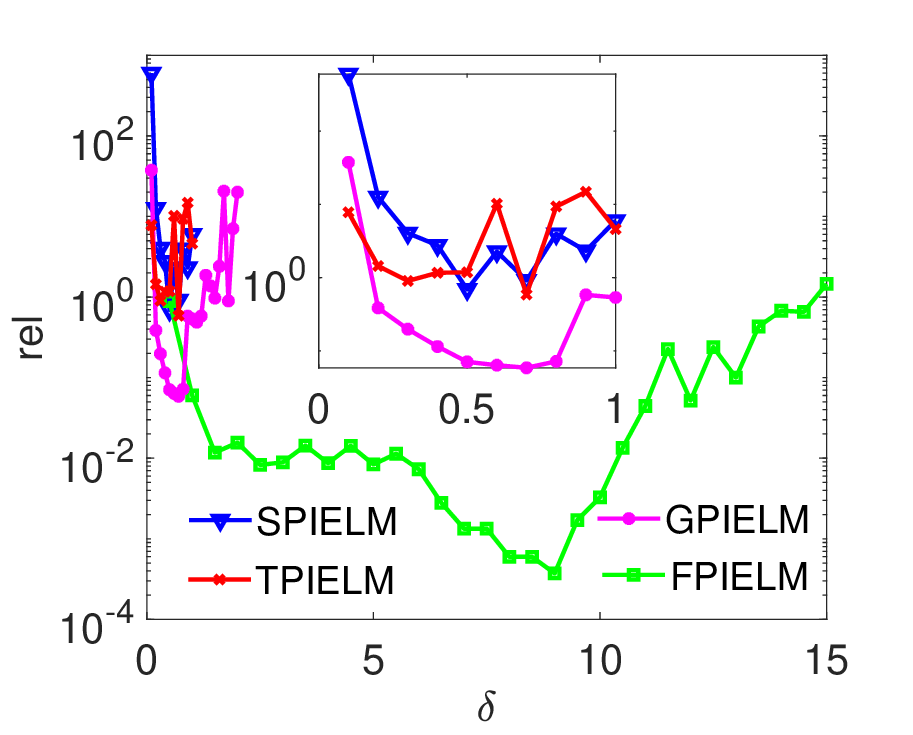}
	}
	\caption{REL VS the scale factor $\delta$ for four PIELM methods to Example \ref{2D_Dirichlet_E1} on hexagram domian derived from  $[-\pi,\pi]\times[-\pi,\pi]$(left) and $[0,3\pi]\times[\pi,2\pi]$(right), respectively.}
	\label{Plot_2D_Diri_E2_rel_0_pi}
\end{figure}

Additionally, Figure~\ref{Plot_2D_rel_time_Sigma} illustrates the variation in relative error and runtime concerning the number of hidden nodes, as well as the variation in relative error with the scale factor $\delta$ for FPIELM on the hexagram domain derived from $[-\pi,\pi]\times[-\pi,\pi]$. Finally, Figure~\ref{Plot_2D_Diri_E2_rel_0_pi} plots the REL comparison curves for the four PIELM methods, showing how REL varies with the scale factor $\delta$ across the hexagram domains derived from $[-\pi,\pi]\times[-\pi,\pi]$ and $[0,3\pi]\times[\pi,2\pi]$.

\begin{table}[!ht]
	\centering
	\caption{REL, $\delta$ and runtime for four PIELM methods to Example \ref{2D_Dirichlet_E2}.}
	\label{Table2D_Dirichlet_E2}
	\begin{tabular}{|l|c|c|c|c|c|}
		\hline  
		           Domain            &          & SPIELM               & GPIELM               &  TPIELM              &  FPIELM    \\    \hline
		                             &$\delta$  &1.2                   &1.4                   &0.6                   &8.5           \\ 
		$[-\pi,\pi]\times[-\pi,\pi]$ & REL      &$3.469\times 10^{-9}$ &$6.949\times 10^{-11}$&$1.453\times 10^{-9}$ &$3.876\times 10^{-10}$   \\  
		                             & Time(s)  &6.121                 &2.794                 &2.868                      &2.095    \\  \hline
		                             &$\delta$  &0.6                   &0.6                   &0.3                   &8.5           \\ 
		$[0,3\pi]\times[\pi,2\pi]$   & REL      &$0.792$               &$0.0147$              &$0.422$               &$7.866\times 10^{-5}$   \\ 
		                             & Time(s)  &4.879                 &2.541                 &2.651                 &2.011    \\  \hline
	\end{tabular}
\end{table}

Based on the data in Table~\ref{Table2D_Dirichlet_E2} and the point-wise absolute errors depicted in Figs.~\ref{2D_Dirichlet_E2:Perr2Sigmoid} -- \ref{2D_Dirichlet_E2:Perr2SIN}, we can infer that the FPIELM method outperforms the PIELM methods with sigmoid, tanh, and Gaussian activation functions. As observed in the curves in Fig.~\ref{2D_Navier_E1:rels}, the performance of the FPIELM method stabilizes as the number of hidden nodes increases. The time-neuron curve in Fig.~\ref{2D_Navier_E1:times} shows that the runtime exhibits a superlinear increase with the growing number of hidden nodes, although the increase remains gradual. Lastly, Figs.~\ref{2D_Dirichlet_E2:sigma} and \ref{Plot_2D_Diri_E2_rel_0_pi} illustrate that the FPIELM method remains stable and robust across variations in the scale factor $\sigma$. In summary, our FPIELM method demonstrates not only high accuracy but also efficiency and robustness.

\begin{example}\label{3D_Dirichlet_E1}

We now approximate the solution of the biharmonic equation \eqref{eq:biharmonic} in a three-dimensional cubic domain $\Omega = [1,3] \times [1,3] \times [1,3]$ containing one large hole (cyan) and eight smaller holes (red and blue), as shown in Fig.~\ref{fig:ThreeDim_Holes}. An analytical solution and its corresponding source term are given by
\begin{equation*}
u(x_1,x_2,x_3) = 50e^{-0.25(x_1 + x_2 + x_3)},
\end{equation*}
and
\begin{equation*}
f(x_1,x_2,x_3) = \frac{225}{128}e^{-0.25(x_1 + x_2 + x_3)},
\end{equation*}
respectively. The boundary functions $g(x_1,x_2,x_3)$ and $h(x_1,x_2,x_3)$ can be directly obtained from the exact solution; however, we omit them here for brevity.
\begin{figure}[H]
	\centering
	\subfigure[Domain of interest]{       
		\label{3D_Holes}
		\includegraphics[scale=0.425]{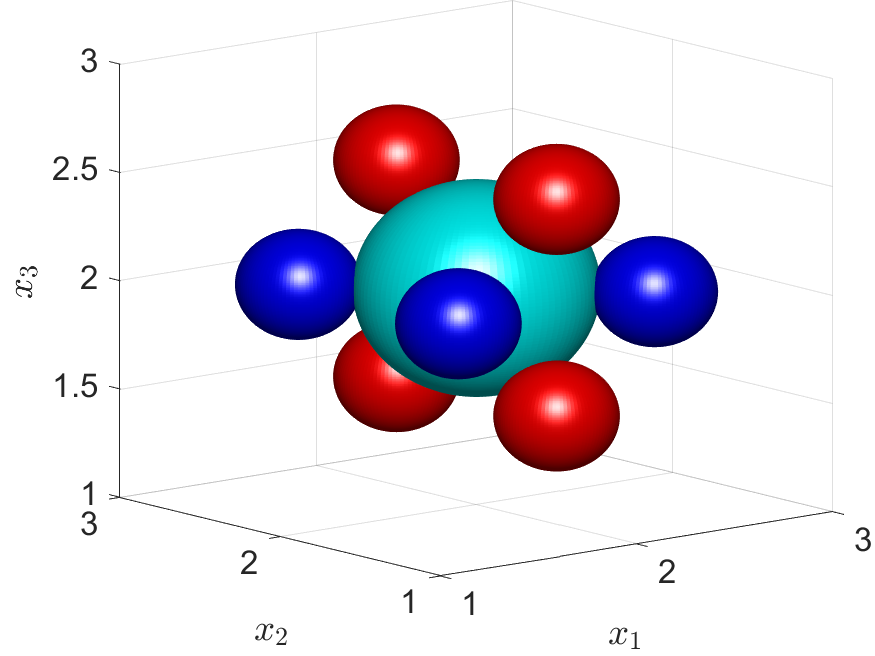}
	}
	\subfigure[Cut planes of domain]{       
		\label{Slice_3DHoles}
		\includegraphics[scale=0.425]{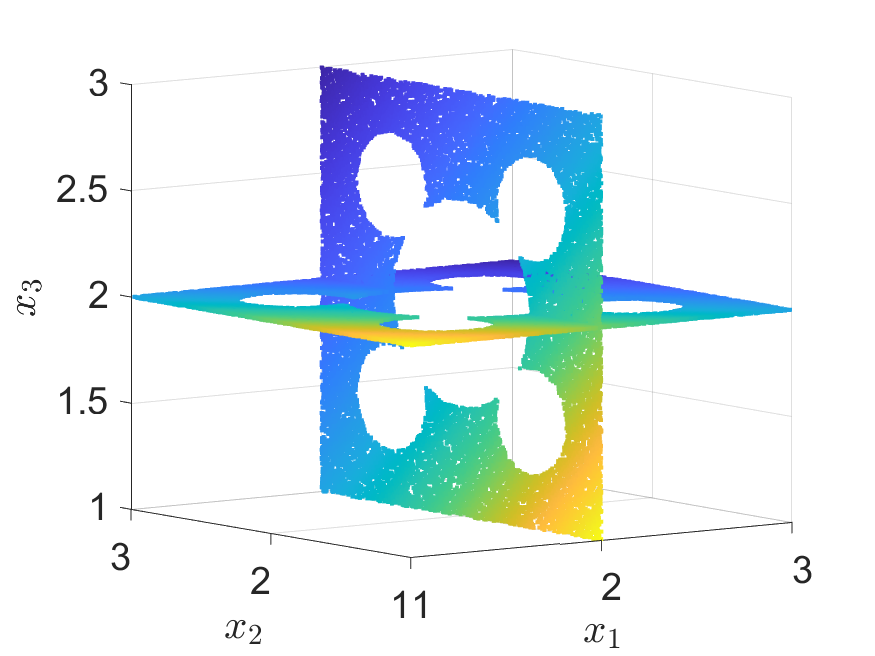}
	}
	\caption{The domain of interest with several holes and its cut planes.}
	\label{fig:ThreeDim_Holes}
\end{figure}

\end{example}

\begin{figure}[H]
	\centering
	\subfigure[Analytical solution]{
		\label{3D_Dirichlet_E1:a} 
		\includegraphics[scale=0.325]{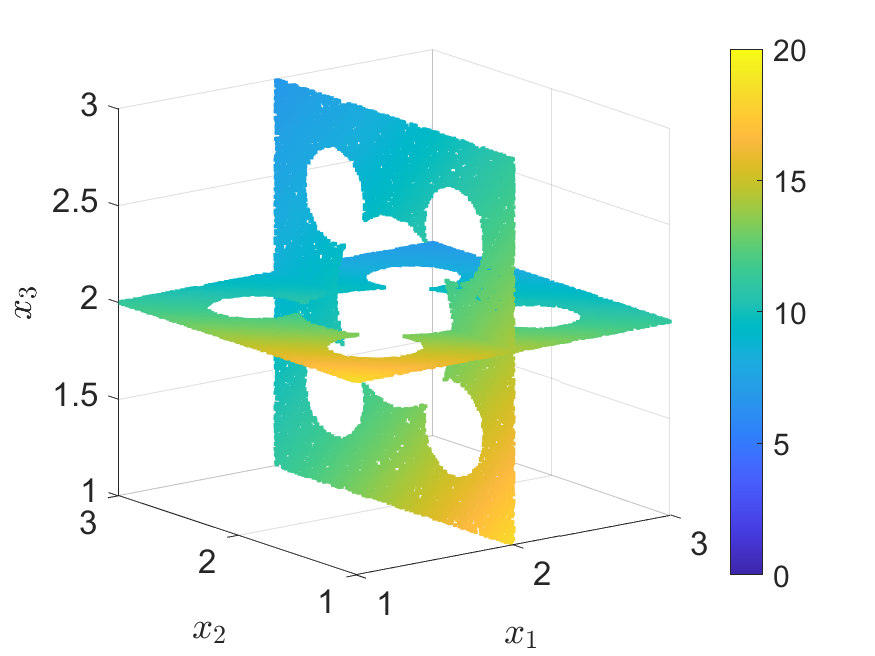}
	}
	\subfigure[Point-wise error of SPIELM]{
		\label{3D_Dirichlet_E1:Perr2Sigmoid}
		\includegraphics[scale=0.325]{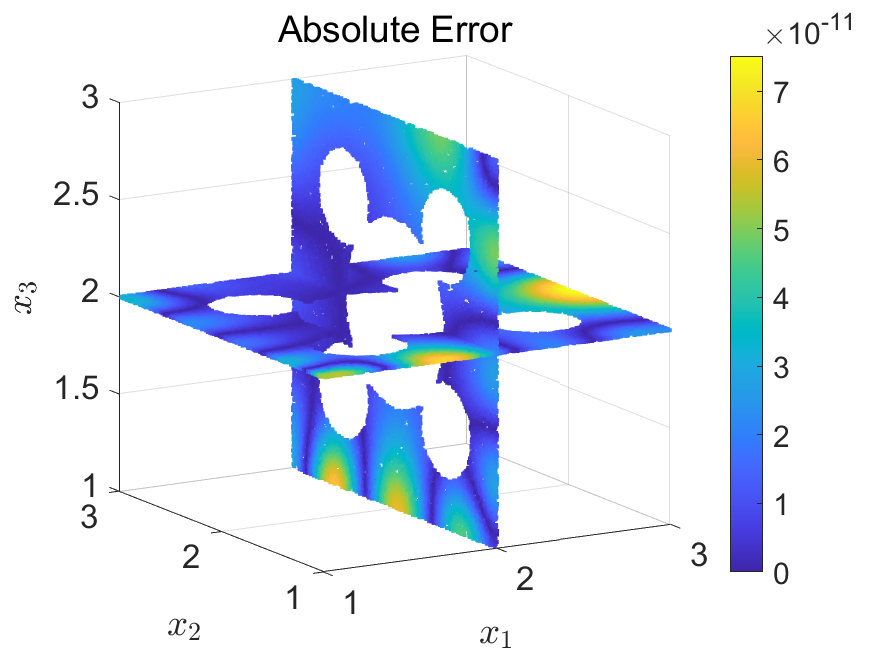}
	}
	\subfigure[Point-wise error of GPIELM]{
		\label{3D_Dirichlet_E1:Perr2Gauss}
		\includegraphics[scale=0.325]{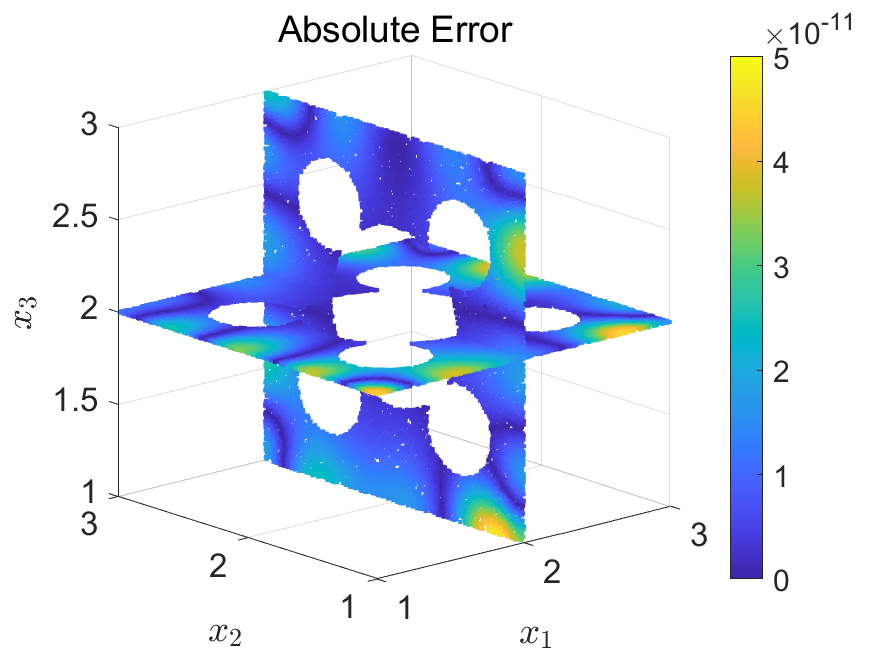}
	}
	\subfigure[Point-wise error of TPIELM]{
		\label{3D_Dirichlet_E1:Perr2Tanh}
		\includegraphics[scale=0.325]{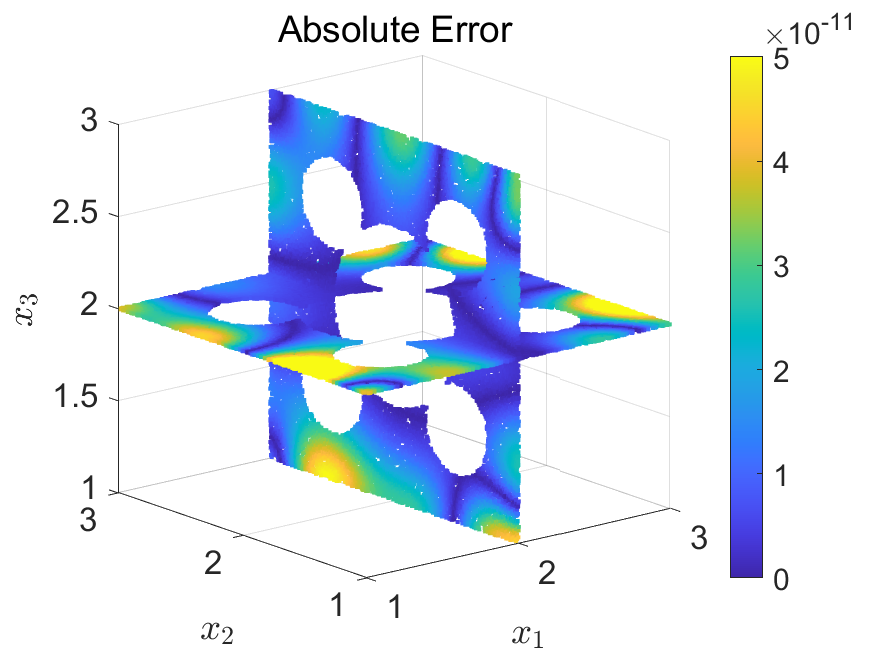}
	}
	\subfigure[Point-wise error of FPIELM]{
		\label{3D_Dirichlet_E1:Perr2SIN}
		\includegraphics[scale=0.325]{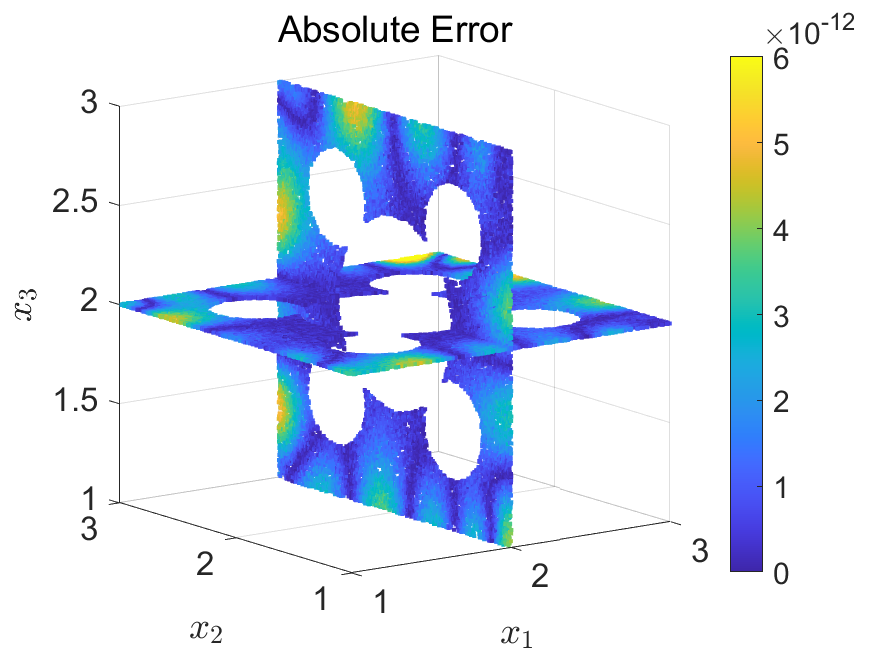}
	}
	\caption{Numerical results for four PIELM methods to Example~\ref{3D_Dirichlet_E1}.}
	\label{Plot_3D_E1}
\end{figure}

\begin{figure}[H]
	\centering
	\subfigure[REL VS the hidden node]{
		\label{3D_Dirichlet_E1:rels} 
		\includegraphics[scale=0.325]{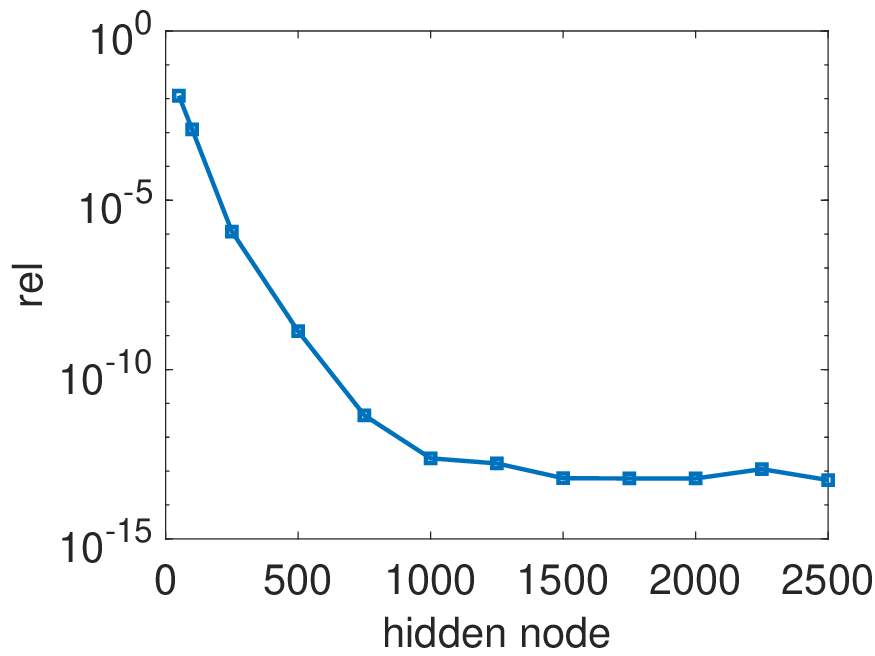}
	}
	\subfigure[Runtime VS the hidden node]{
		\label{3D_Dirichlet_E1:times}
		\includegraphics[scale=0.325]{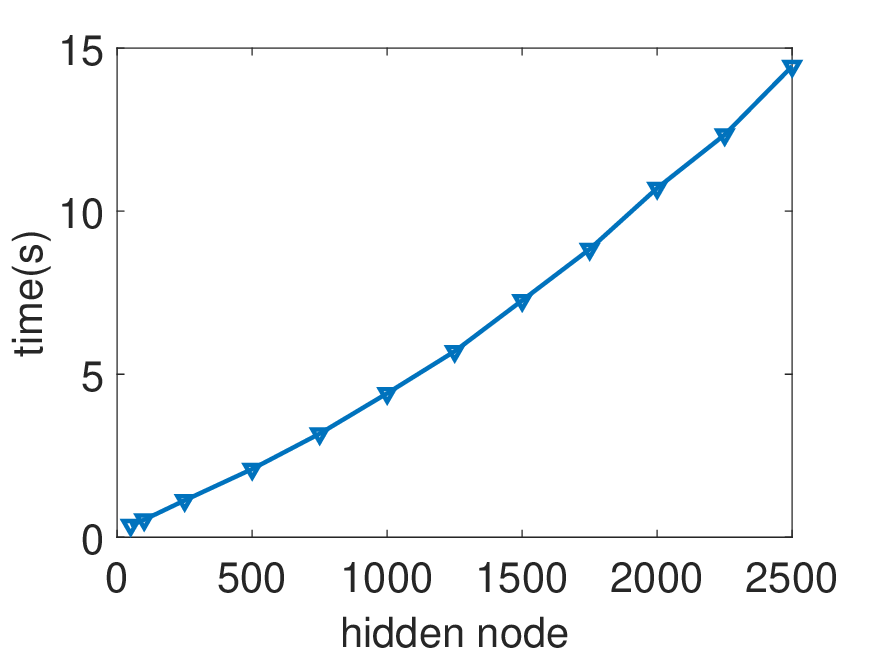}
	}
	\subfigure[REL VS the scale factor $\sigma$]{
		\label{3D_Dirichlet_E1:sigma}
		\includegraphics[scale=0.325]{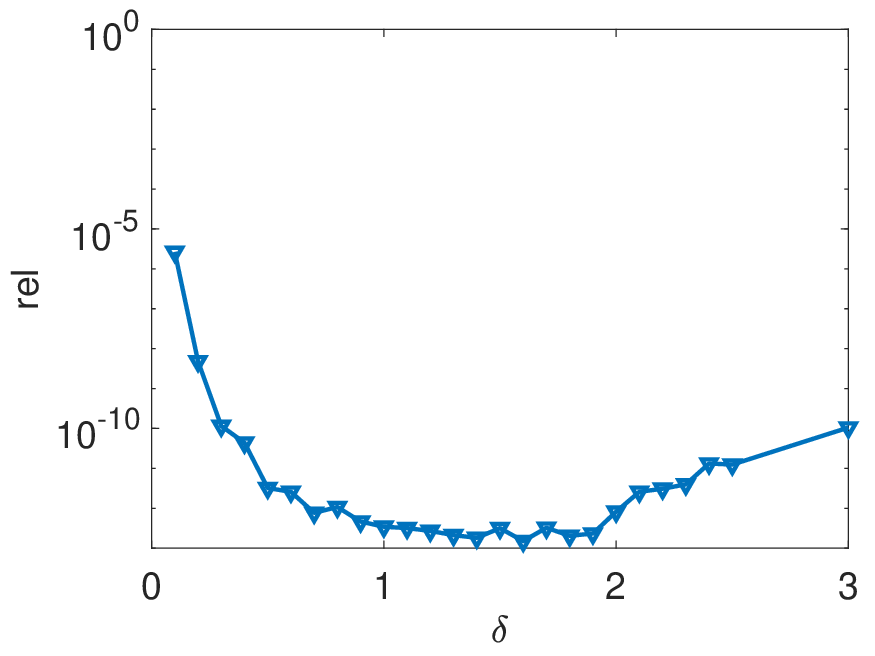}
	}
	\caption{REL VS the hidden node (left), runtime VS the hidden node (middle) and REL VS the scale factor $\sigma$ (right) for the PIELM method to Example \ref{3D_Dirichlet_E1}.}
	\label{Plot_3DE1_rel_time_Sigma}
\end{figure}

\begin{figure}[H]
	\centering 
	\includegraphics[scale=0.75]{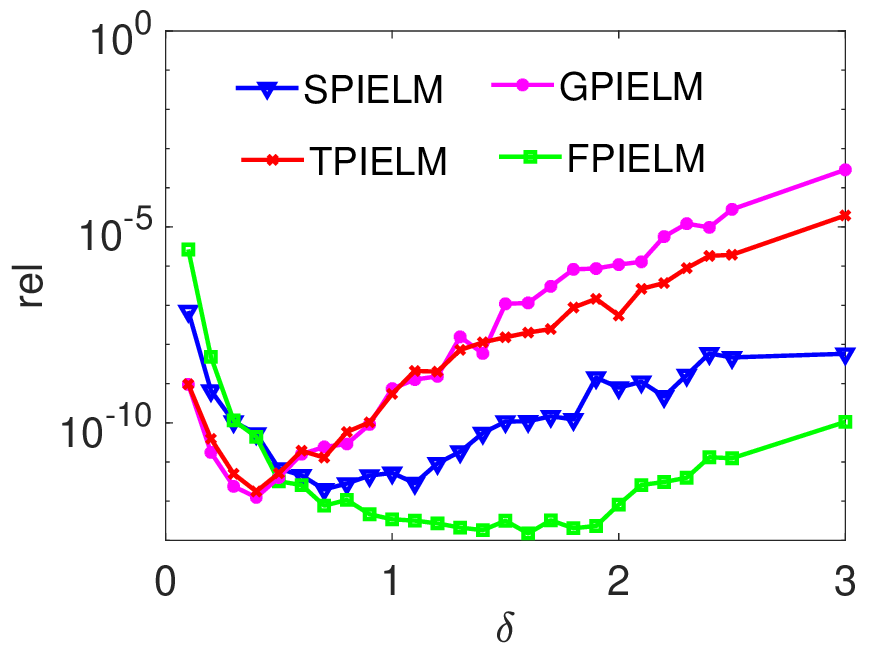}
	\caption{REL VS the scale factor $\delta$ for four PIELM models to Example \ref{3D_Dirichlet_E1}.}
	\label{Plot_3D_Diri_E2_rel_13}
\end{figure}

In this example, the number of hidden units for each of the four PIELM methods is set to 2,000. The testing set consists of 40,000 random points within the cut planes, excluding the holes, and includes planes parallel to the $xoy$ and $yoz$ coordinate planes. Additionally, 6,000 points randomly sampled from the six boundary surfaces are integrated into the testing set. The results are presented in Figs.~\ref{Plot_3D_E1} -- \ref{Plot_3D_Diri_E2_rel_13}.

Based on the point-wise error shown in Figs.~\ref{3D_Dirichlet_E1:Perr2Sigmoid} -- \ref{3D_Dirichlet_E1:Perr2SIN}, our proposed FPIELM method outperforms the other three methods in solving the biharmonic equation \eqref{eq:biharmonic} in three-dimensional space. The relative errors presented in Fig.~\ref{Plot_3D_Diri_E2_rel_13} indicate that the FPIELM method remains stable across various scale factors $\delta$ and is superior to the other methods. Furthermore, the curves of relative error and runtime in Fig.~\ref{Plot_3DE1_rel_time_Sigma} demonstrate the favorable stability and robustness of the FPIELM method as both the number of hidden nodes and the scale factor $\delta$ increase. It is also noteworthy that the coefficient matrix in the final linear system is full rank for this example, eliminating the need for adjustment through Tikhonov regularization.

\subsection{Numerical examples for Navier boundary}
\begin{example}\label{2D_Navier_E1} 
In this example, our goal is to solve the biharmonic equation \eqref{eq:biharmonic} with Navier boundary conditions in a regular rectangular domain $\Omega = [x_1^{\text{min}}, x_1^{\text{max}}] \times [x_2^{\text{min}}, x_2^{\text{max}}]$. An exact solution is given by
\begin{equation*}
u(x_1,x_2) = \sin(x_1^2 + x_2^2),
\end{equation*}
which naturally leads to the essential boundary function $g(x_1, x_2)$. The function $h(x_1, x_2)$ is defined as
\begin{equation*}
    h(x_1, x_2) = 2\cos(x_1^2 + x_2^2) - 4x_1^2\sin(x_1^2 + x_2^2) + 2\cos(x_1^2 + x_2^2) - 4x_2^2\sin(x_1^2 + x_2^2).
\end{equation*}
With careful calculations, the expression for the force term is determined as
\begin{equation*}
	\begin{aligned}
	   f(x_1, x_2) &= 16x_1^4\sin(x_1^2 + x_2^2) + 16x_2^4\sin(x_1^2 + x_2^2) - 64x_1^2\cos(x_1^2 + x_2^2) - 64x_2^2\cos(x_1^2 + x_2^2) \\
	   &\quad + 32x_1^2x_2^2\sin(x_1^2 + x_2^2) - 32\sin(x_1^2 + x_2^2).
	\end{aligned}
\end{equation*}
\end{example}

\begin{figure}[H]
	\centering
	\subfigure[Analytical solution]{
		\label{2D_Navier_E1:a} 
		\includegraphics[scale=0.325]{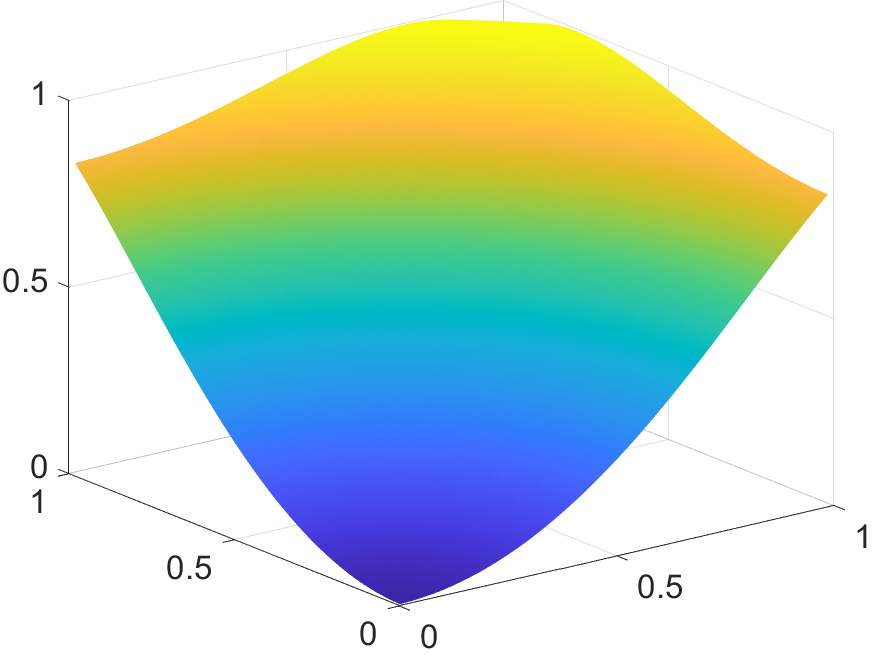}
	}
	\subfigure[Point-wise error of SPIELM]{
		\label{2D_Navier_E1:Perr2Sigmoid}
		\includegraphics[scale=0.325]{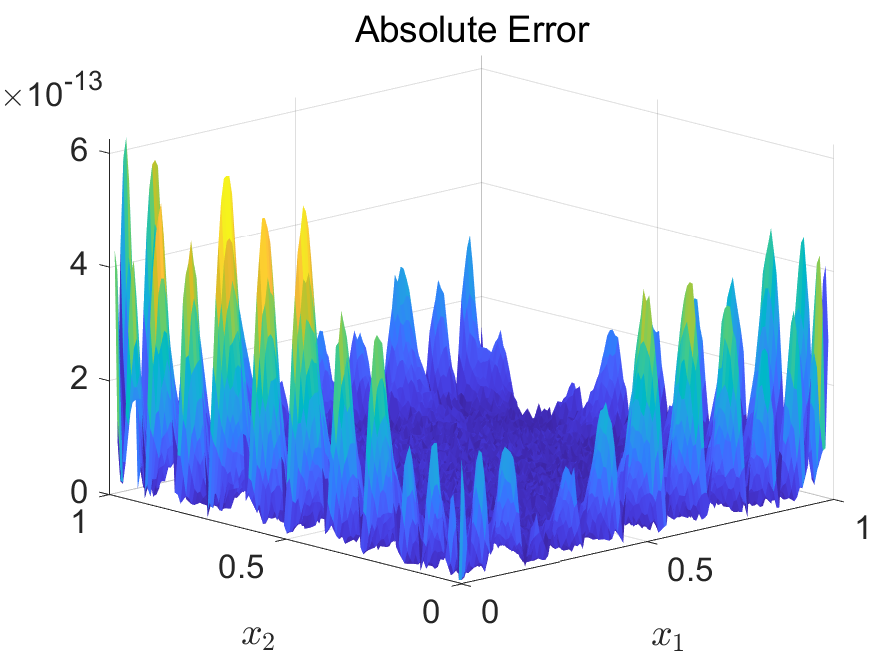}
	}
	\subfigure[Point-wise error of GPIELM]{
		\label{2D_Navier_E1:Perr2Gauss}
		\includegraphics[scale=0.325]{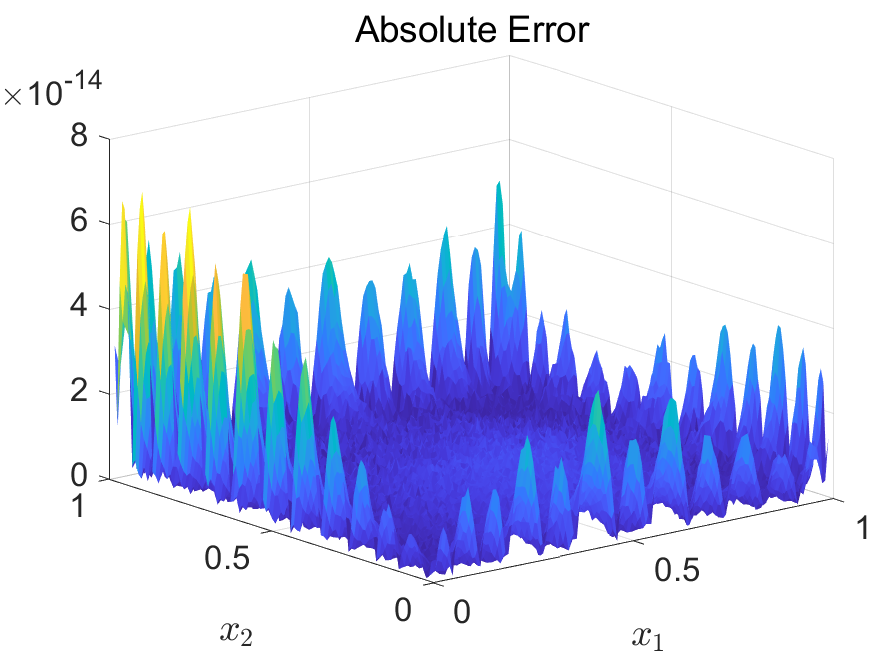}
	}
	\subfigure[Point-wise error of TPIELM]{
		\label{2D_Navier_E1:Perr2Tanh}
		\includegraphics[scale=0.325]{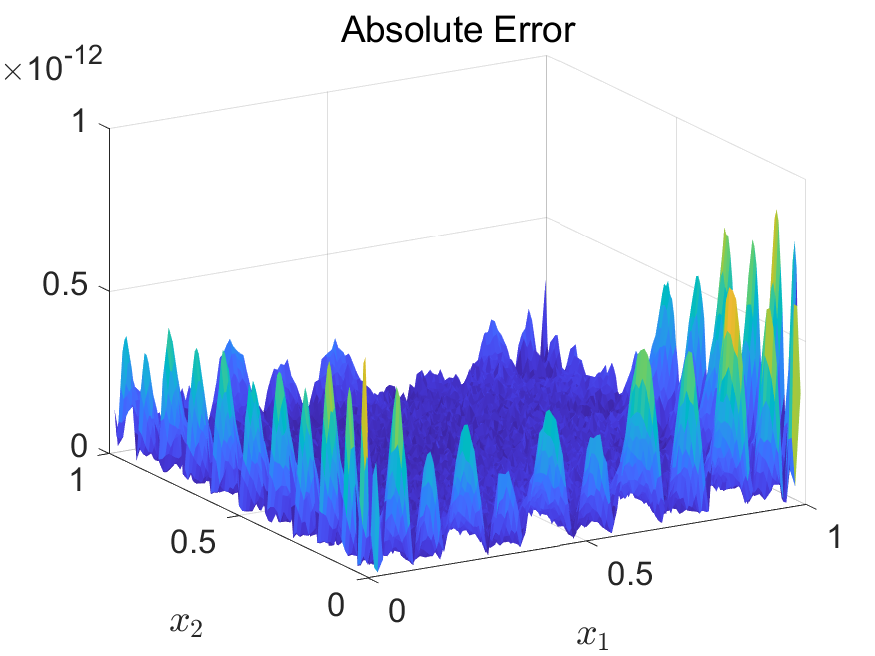}
	}
	\subfigure[Point-wise error of FPIELM]{
		\label{2D_Navier_E1:Perr2SIN}
		\includegraphics[scale=0.325]{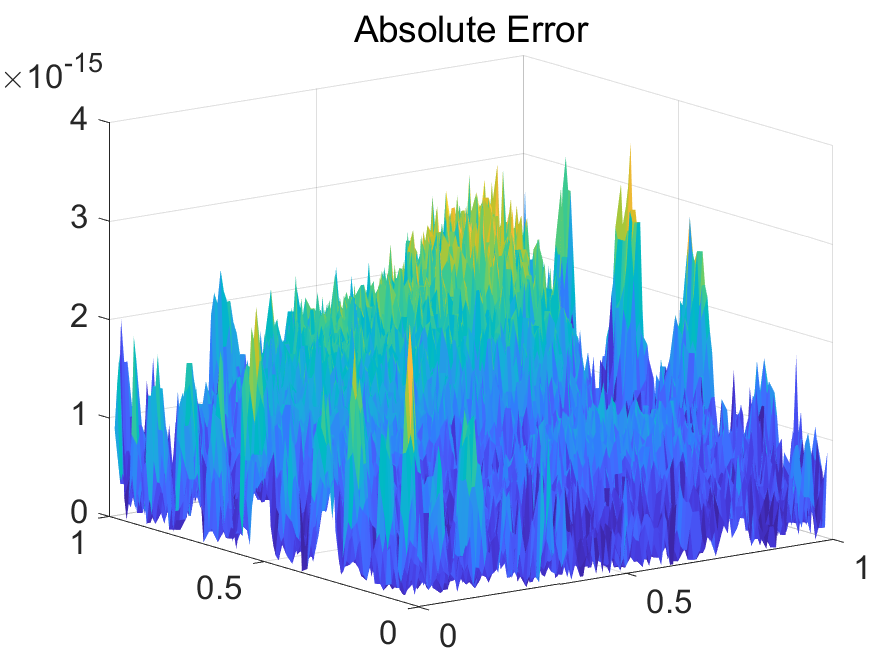}
	}
	\caption{Numerical results for four PIELM methods to Example~\ref{2D_Navier_E1} on domian $[0,1]\times[0,1]$.}
	\label{Plot_2D_Navier_E1}
\end{figure}

In this example, the setup for the four PIELM methods is identical to that in Example~\ref{2D_Dirichlet_E1}. An approximate solution of \eqref{eq:biharmonic} is obtained by applying the four methods on 16,384 grid points within the rectangular domains $[0,1] \times [0,1]$ and $[0,4] \times [0,4]$. The numerical results are listed in Table~\ref{Table2D_Navier_E1} for these two scenarios, and the point-wise absolute error for the four PIELM methods is shown in Fig.~\ref{Plot_2D_Navier_E1} on the domain $[0,1] \times [0,1]$. 

\begin{figure}[H]
	\centering
	\subfigure[REL VS the hidden node]{
		\label{2D_Navier_E1:rels} 
		\includegraphics[scale=0.325]{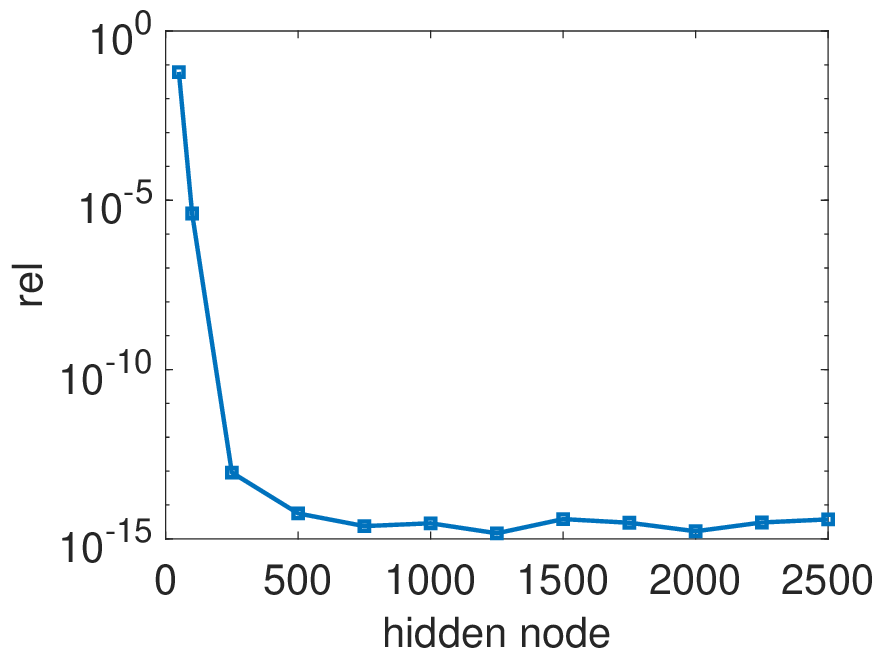}
	}
	\subfigure[Runtime VS the hidden node]{
		\label{2D_Navier_E1:times}
		\includegraphics[scale=0.325]{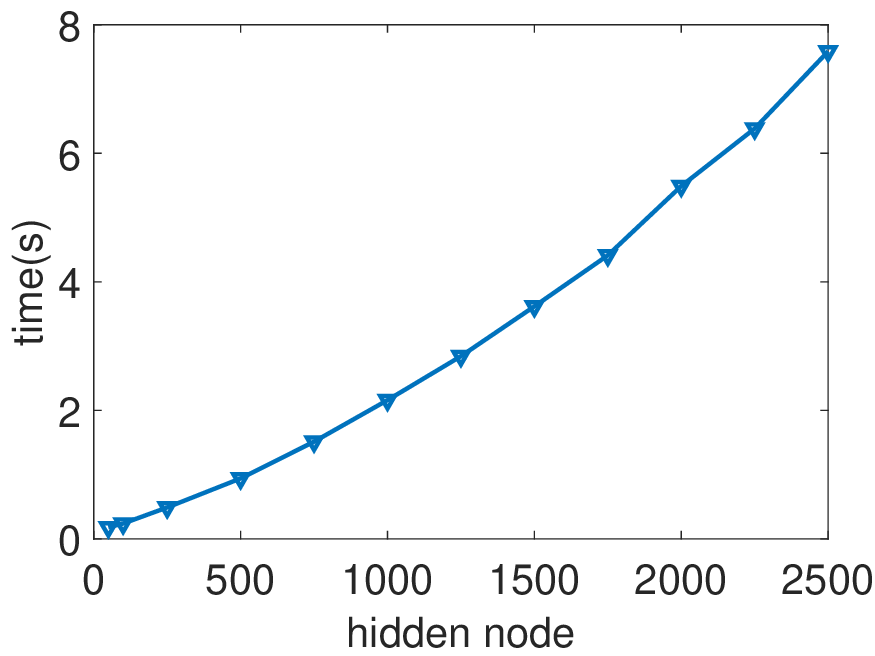}
	}
	\subfigure[REL VS the scale factor $\delta$]{
		\label{2D_Navier_E1:sigma}
		\includegraphics[scale=0.325]{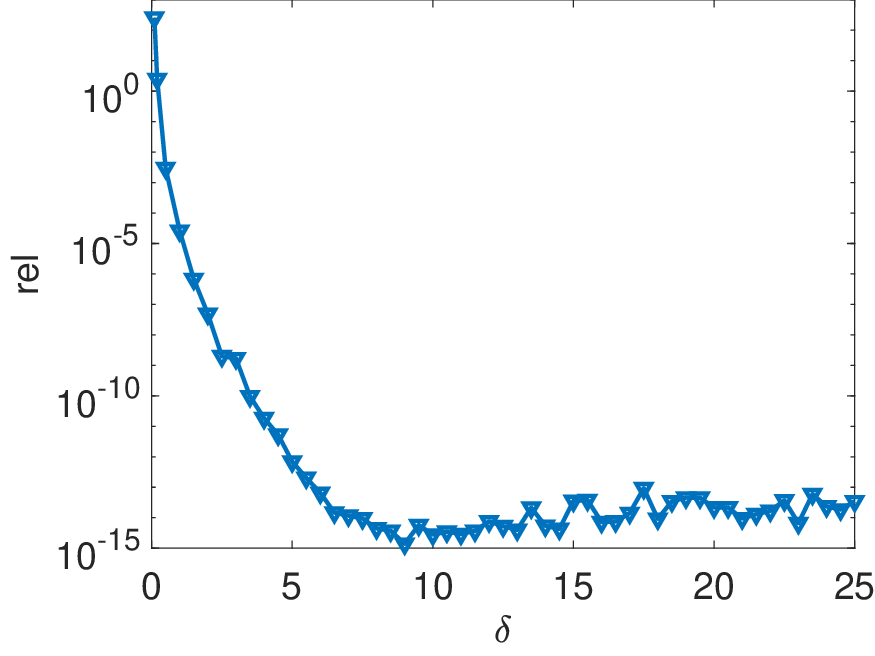}
	}
	\caption{REL VS the hidden nodes (left), runtime VS the hidden node (middle) and REL VS the scale factor $\delta$ (right) for FPIELM method to Example \ref{2D_Navier_E1} on  domain $[0,1]\times[0,1]$.}
	\label{Plot_Navier2D_E1_rel_time_Sigma}
\end{figure}

\begin{figure}[H]
	\centering 
	\subfigure[]{
		\includegraphics[scale=0.5]{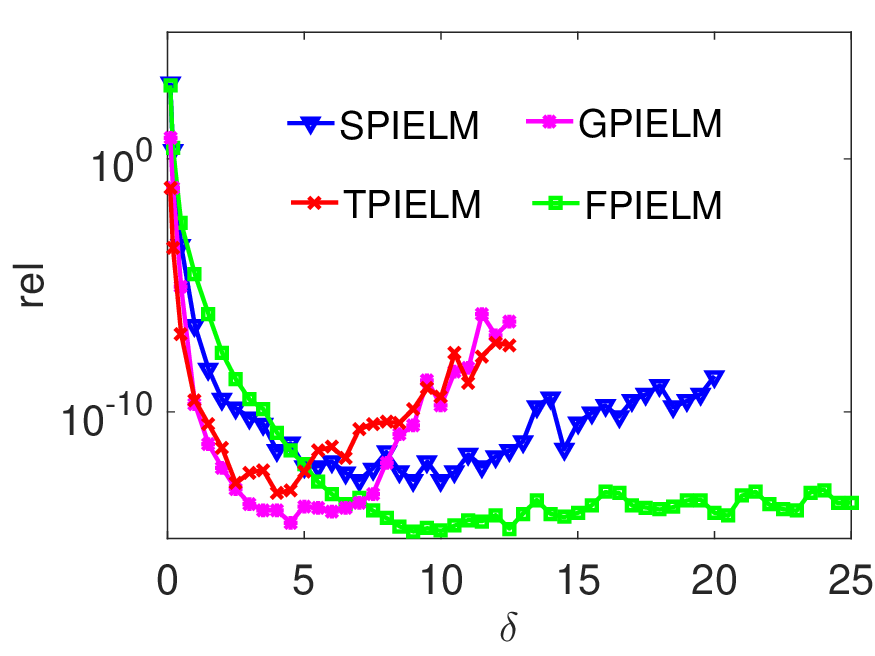}
	}
	\subfigure[]{
		\includegraphics[scale=0.5]{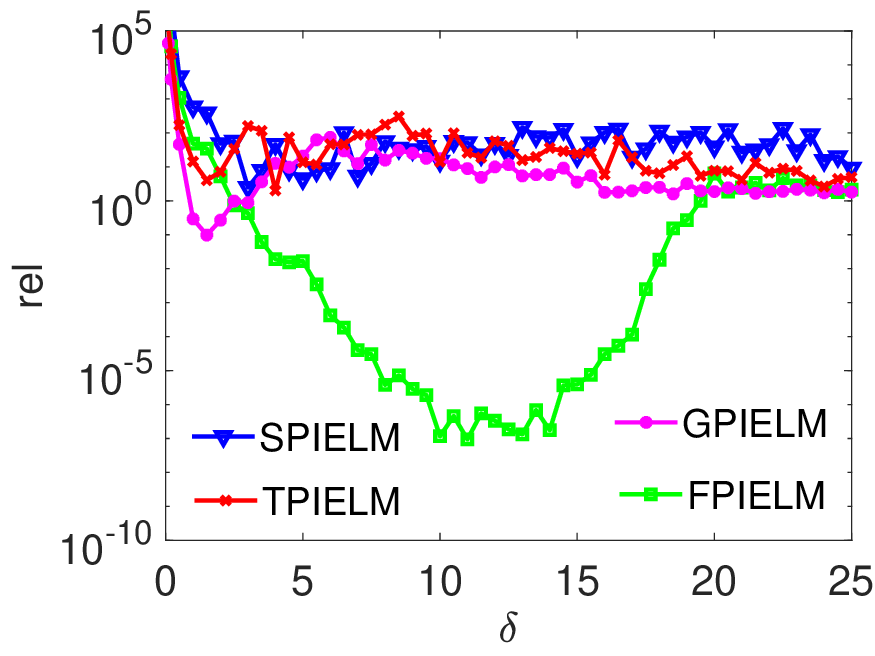}
	}
	\caption{REL VS the scale factor $\delta$ for the four PIELM methods to Example \ref{2D_Navier_E1} on regualr domains $[0,1]\times[0,1]$(left) and $[0,4]\times[0,4]$(right).}
	\label{Plot_Navier2D_E2_rel_04}
\end{figure}

Additionally, Fig.~\ref{Plot_Navier2D_E1_rel_time_Sigma} illustrates the variation in relative error and runtime concerning the number of hidden nodes, as well as the variation in relative error with the scale factor $\delta$ for FPIELM on the domain $[0,1] \times [0,1]$. Finally, Fig.~\ref{Plot_Navier2D_E2_rel_04} plots the comparison curves for REL as a function of the scale factor $\delta$ for the four PIELM methods on the rectangular domains $[0,1] \times [0,1]$ and $[0,4] \times [0,4]$, respectively.

\begin{table}[!ht]
	\centering
	\caption{REL, $\delta$ and runtime for four PIELM methods to Example \ref{2D_Navier_E1}.}
	\label{Table2D_Navier_E1}
	\begin{tabular}{|l|c|c|c|c|c|}
		\hline  
		Domain                &         & SPIELM               & GPIELM               &TPIELM                &  FPIELM               \\  \hline
	                          &$\delta$ &6.5                   &4.1                   & 3.2                  &9.0           \\ 
		   $[0,1]\times[0,1]$ &REL      &$1.222\times 10^{-13}$&$1.315\times 10^{-14}$&$1.552\times 10^{-13}$&$1.891\times 10^{-15}$   \\  
		     	              & Time(s) &11.958                &4.648                 &4.823                 &2.775    \\  \hline
		                      &$\delta$ &3.0                   &1.5                   &4.0                   &11.0           \\ 
           $[0,4]\times[0,4]$ & REL     &$2.471$               &$0.0979$              &$2.016$               &$9.418\times 10^{-8}$   \\  
                              & Time(s) &6.795                 &3.133                 &3.121                 &2.178    \\  \hline
	\end{tabular}
\end{table}

Based on the point-wise errors depicted in Figs.~\ref{2D_Navier_E1:Perr2Sigmoid} -- \ref{2D_Navier_E1:Perr2SIN} and the relative errors presented in Fig.~\ref{Plot_Navier2D_E1_rel_time_Sigma}, our proposed FPIELM method outperforms the other three methods for solving the biharmonic equation \eqref{eq:biharmonic} with Navier boundary conditions in both unitized and non-unitized domains. The FPIELM method remains stable across various scale factors $\delta$, and its runtime increases linearly with the number of hidden nodes. Additionally, the curves of relative error and runtime in Fig.~\ref{Plot_Navier2D_E2_rel_04} highlight the favorable stability and robustness of the FPIELM method as both the number of hidden nodes and the scale factor $\delta$ increase.

\begin{example}\label{2D_Navier_E2}
We consider the biharmonic equation \eqref{eq:biharmonic} on two porous domains $\Omega$, derived from the two-dimensional domains $[-1, 1] \times [-\pi, \pi]$ and $[0, 4] \times [0, 4\pi]$, respectively. The exact solution is given by
\begin{equation*}
u(x_1, x_2) = e^{x_1} \sin(x_2),
\end{equation*}
such that $f(x_1, x_2) = 0$ and $k(x_1, x_2) = 0$. The boundary conditions $g(x_1, x_2)$ can be easily obtained through direct calculation.
\end{example}

\begin{figure}[H]
	\centering
	\subfigure[Analytical solution]{
		\label{2D_Navier_E2:a} 
		\includegraphics[scale=0.325]{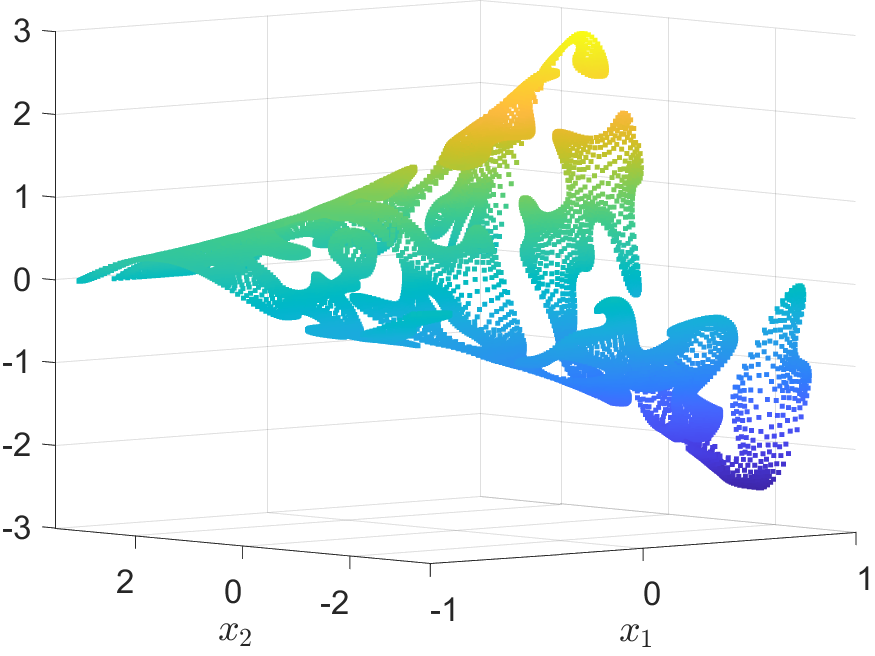}
	}
	\subfigure[Point-wise error (SIgmoid)]{
		\label{2D_Navier_E2:Perr2Tanh}
		\includegraphics[scale=0.325]{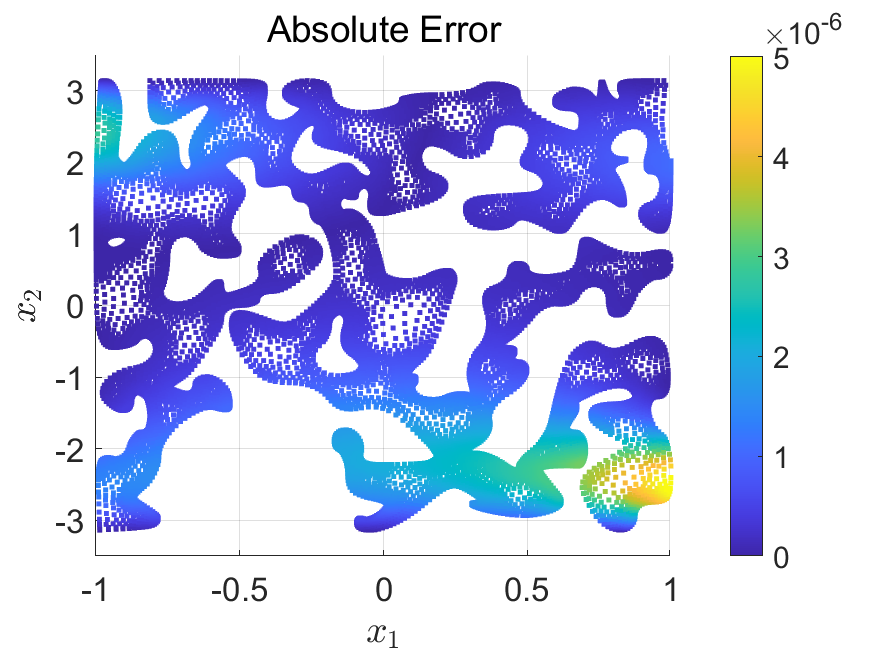}
	}
	\subfigure[Point-wise error (Gaussian)]{
		\label{2D_Navier_E2:Perr2Gauss}
		\includegraphics[scale=0.325]{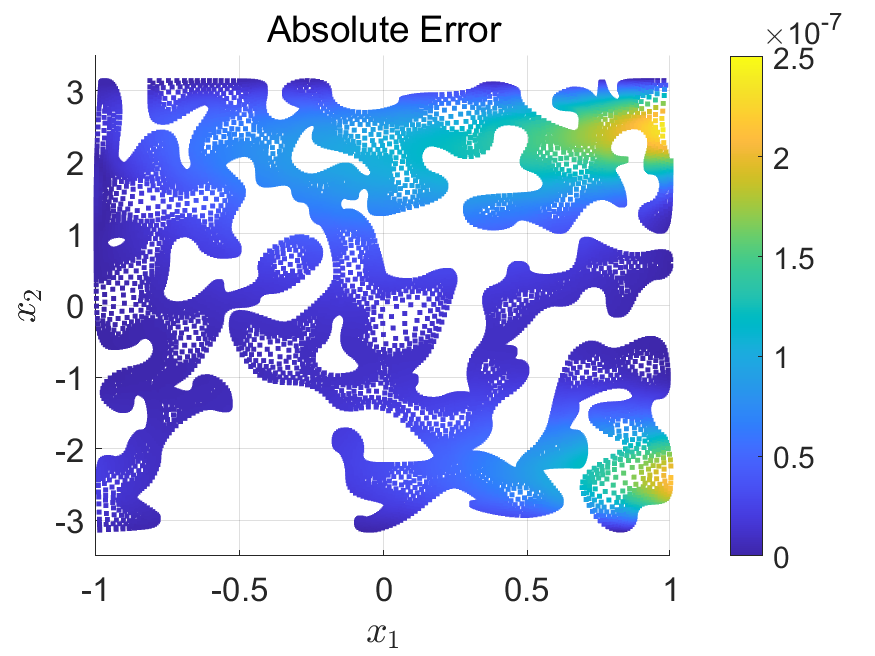}
	}
	\subfigure[Point-wise error for TPIELM]{
		\label{2D_Navier_E2:Perr2SIN}
		\includegraphics[scale=0.325]{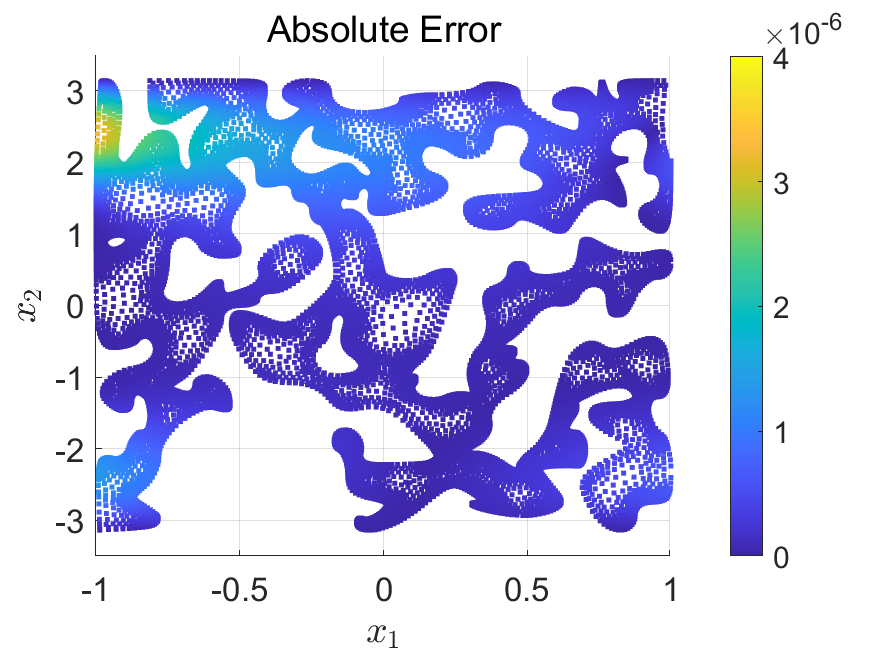}
	}
	\subfigure[Point-wise error for FPIELM]{
		\label{2D_Navier_E2:Perr2SINCOS}
		\includegraphics[scale=0.325]{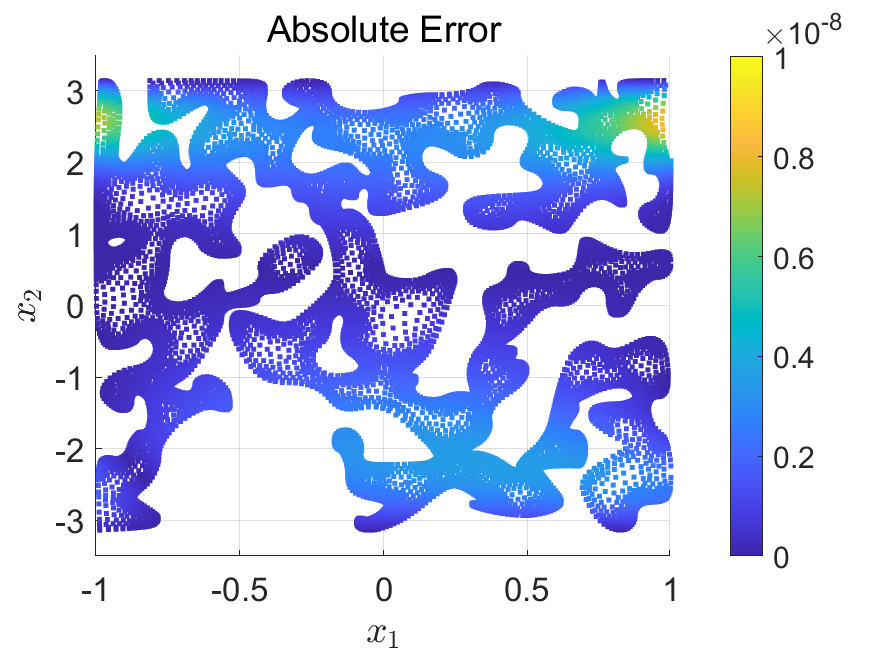}
	}
	\caption{Numerical results for four PIELM methods to Example~\ref{2D_Navier_E2} on porous domian derived from $[-1,1]\times[-\pi,\pi]$.}
	\label{Plot_2D_Navier_E2}
\end{figure}

\begin{figure}[H]
	\centering
	\subfigure[REL VS the hidden nodes]{
		\label{2D_Navier_E2:rels} 
		\includegraphics[scale=0.325]{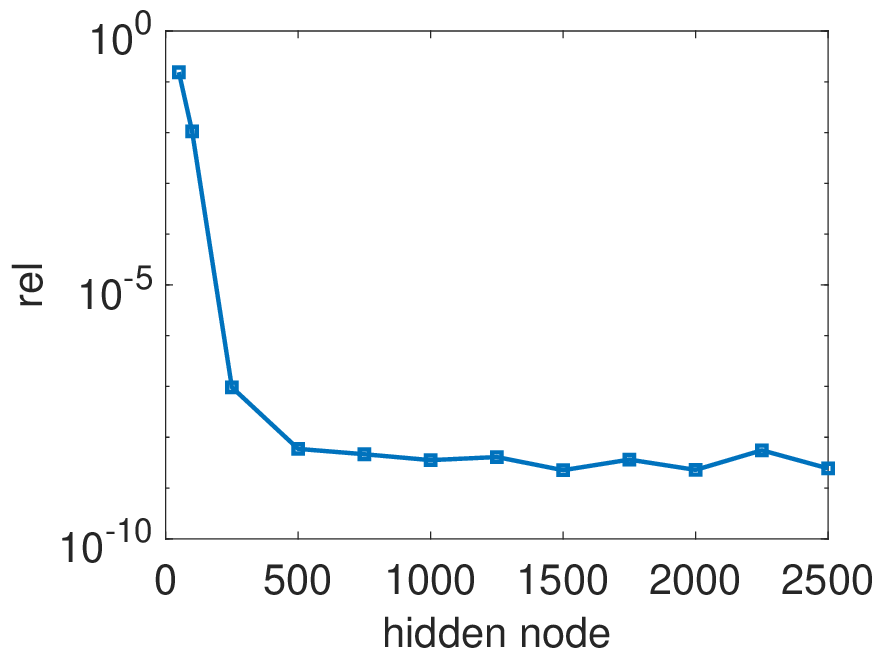}
	}
	\subfigure[Runtime VS the hidden nodes]{
		\label{2D_Navier_E2:times}
		\includegraphics[scale=0.325]{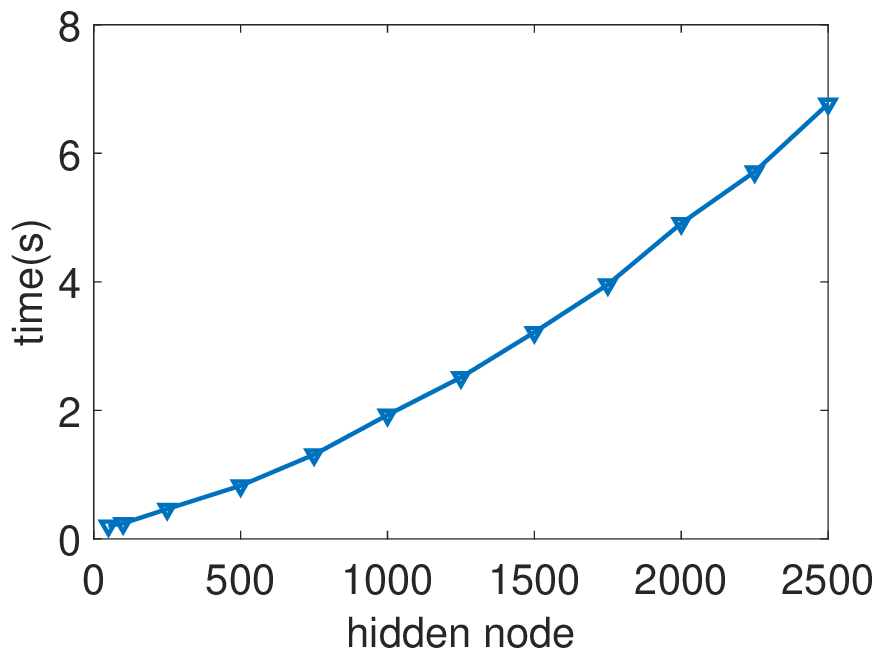}
	}
	\subfigure[REL VS the scale factor $\sigma$]{
		\label{2D_Navier_E2:sigma}
		\includegraphics[scale=0.325]{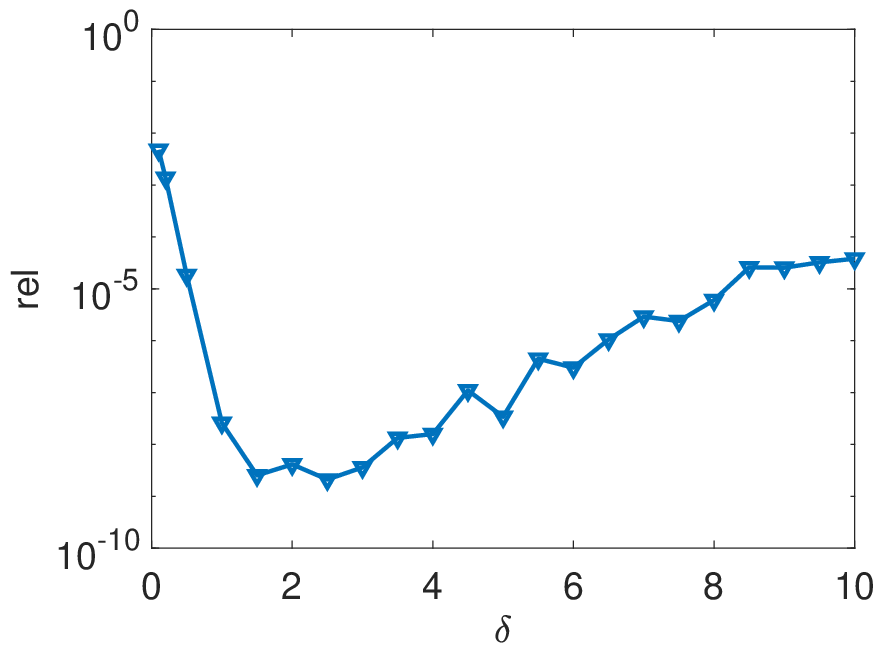}
	}
	\caption{REL VS the hidden nodes (left), runtime VS the hidden nodes (middle) and REL VS the scale factor $\sigma$ (right) to FPIELM method for Example \ref{2D_Navier_E2}  on the porous domain.}
	\label{Plot_2D_Navier_E2_rel_time_Sigma}
\end{figure}

\begin{figure}[H]
	\centering
	\subfigure[Analytical solution]{
		\label{Navier2D_E2:REL_Delta_11_PiPi} 
		\includegraphics[scale=0.5]{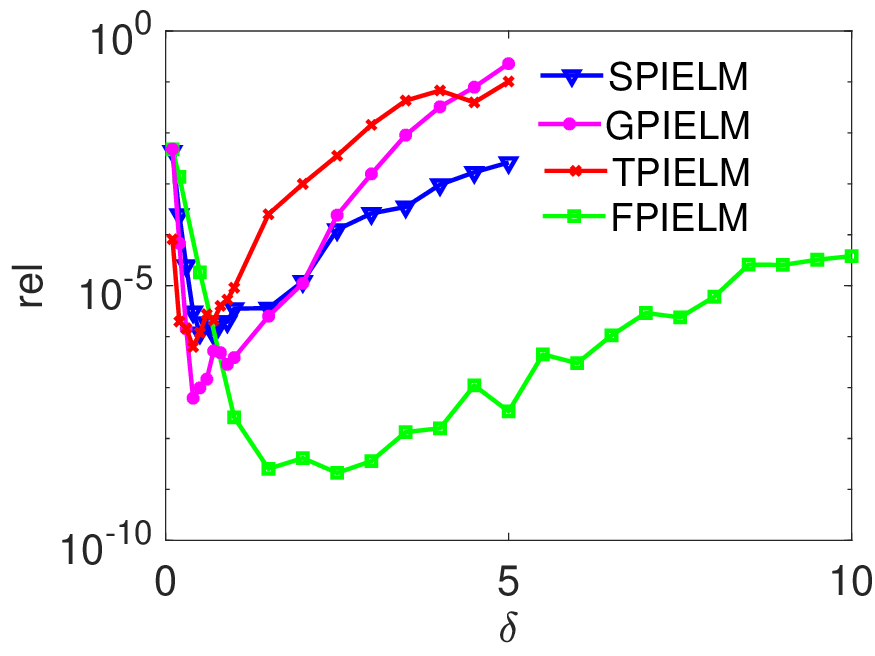}
	}
	\subfigure[Point-wise error (SIgmoid)]{
		\label{Navier2D_E2:REL_Delta_04_04Pi}
		\includegraphics[scale=0.5]{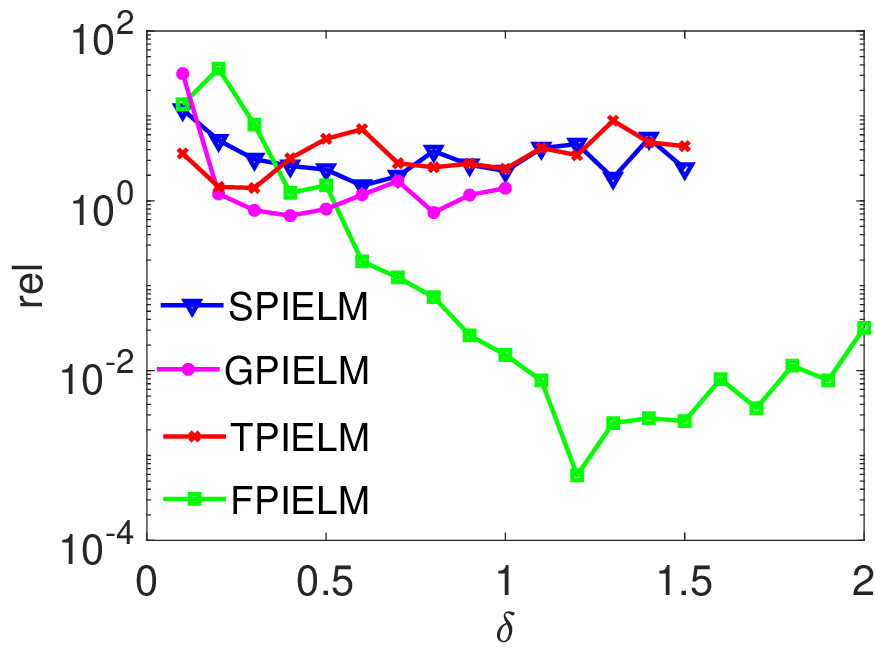}
	}
	\caption{REL VS the scale factor $\delta$ to the four PIELM models for solving Example \ref{2D_Navier_E2} on porous domain from $[-1,1]\times[-\pi,\pi]$(left) and $[0,4]\times[0,4\pi]$(right).}
	\label{rel_delta_11PiPi_0404Pi}
\end{figure}

In this example, the setup for the four PIELM methods is identical to that in Example~\ref{2D_Navier_E1}. An approximation of the solution to the biharmonic equation \eqref{eq:biharmonic} is obtained using the four methods on 20,000 points randomly sampled from the hexagonal domain $\Omega$. The numerical results are presented in Figs.~\ref{Plot_2D_Navier_E2} -- \ref{rel_delta_11PiPi_0404Pi} and in Table~\ref{Table2D_Navier_E2}.

\begin{table}[!ht]
	\centering
	\caption{REL, $\delta$ and runtime for four PIELM methods to Example \ref{2D_Navier_E2}.}
	\label{Table2D_Navier_E2}
	\begin{tabular}{|l|c|c|c|c|c|}
		\hline  
		Domain               &         & SPIELM               & GPIELM               &TPIELM                &  FPIELM               \\  \hline
		                     &$\delta$ &0.7                   &0.4                   & 0.4                  &2.5           \\ 
	$[-1,1]\times[-\pi,\pi]$ &REL      &$1.013\times 10^{-6}$&$6.185\times 10^{-8}$&$6.286\times 10^{-7}$&$2.072\times 10^{-9}$   \\  
		                     & Time(s) &5.902                 &2.788                 &2.715                 &1.924    \\  \hline
		                     &$\delta$ &0.6                   &0.4                   &0.3                  &1.2           \\ 
	$[0,4]\times[0,4\pi]$    & REL     &$1.508$               &$0.668$               &$1.417$               &$5.581\times 10^{-4}$   \\  
		                     & Time(s) &6.095                 &2.748                 &2.842                 &1.941   \\  \hline
	\end{tabular}
\end{table}

Based on the results in Fig.~\ref{Plot_2D_Navier_E2} and Table~\ref{Table2D_Navier_E2}, the FPIELM method demonstrates a remarkable ability to achieve a satisfactory approximation of the solution to \eqref{eq:biharmonic} with Navier boundary conditions, even within an irregular domain in two-dimensional space. Its performance surpasses that of the other three methods. As shown in Figs.~\ref{2D_Navier_E2:rels} and \ref{2D_Navier_E2:times}, the FPIELM method remains stable as the number of hidden nodes varies, with computational complexity increasing only linearly as hidden nodes increase. Regarding the scale factor $\sigma$, the FPIELM method shows insensitivity across a broad range, while the other three methods are sensitive in the unitized domain and fail to converge on the non-unitized domain. Additionally, the coefficient matrix in the final linear systems for this example is full rank.

\begin{example}\label{3D_Navier_E1}
We consider the biharmonic equation \eqref{eq:biharmonic} in three-dimensional space and obtain its numerical solution within a closed domain $\Omega$, bounded by two spherical surfaces with radii $r_1 = 0.2$ and $r_2 = 1.0$, respectively.

An exact solution is given by
\begin{equation*}
    u(x_1, x_2, x_3) = \sin(\pi x_1) \sin(\pi x_2) \sin(\pi x_3),
\end{equation*}
which naturally determines the boundary conditions $g(x_1, x_2, x_3)$ on $\partial \Omega$. The second-order boundary function is given by \( k(x_1, x_2, x_3) = -3\pi^2 \sin(\pi x_1) \sin(\pi x_2) \sin(\pi x_3) \). Through careful calculation, the source term is determined as \( f(x_1, x_2, x_3) = 9\pi^4 \sin(\pi x_1) \sin(\pi x_2) \sin(\pi x_3) \).
\end{example}

\begin{figure}[H]
	\centering
	\subfigure[Domain]{
		\label{3D_Navier_E1:Domain} 
		\includegraphics[scale=0.335]{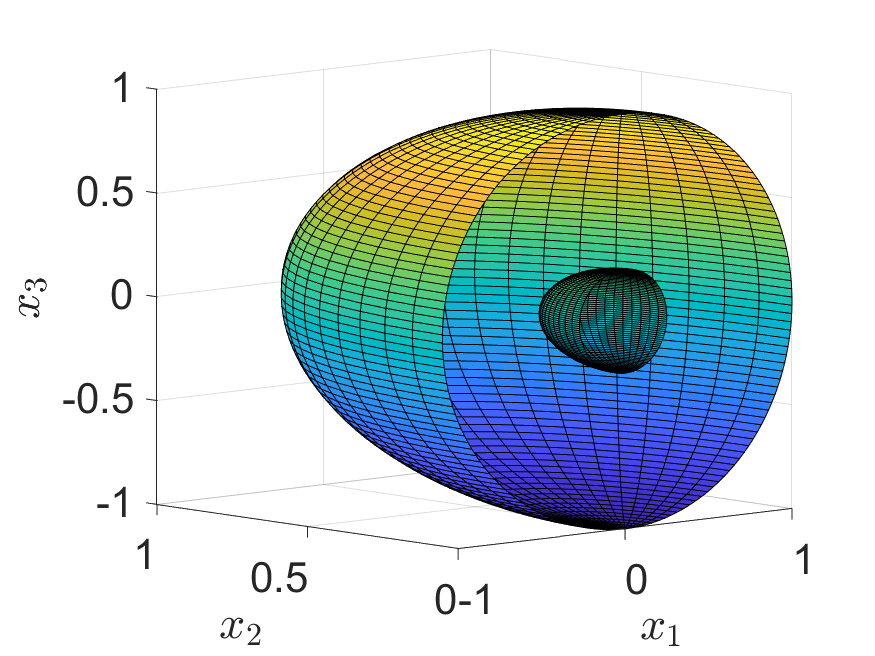}
	}
	\subfigure[Analytical solution]{
		\label{3D_Navier_E1:Solu} 
		\includegraphics[scale=0.335]{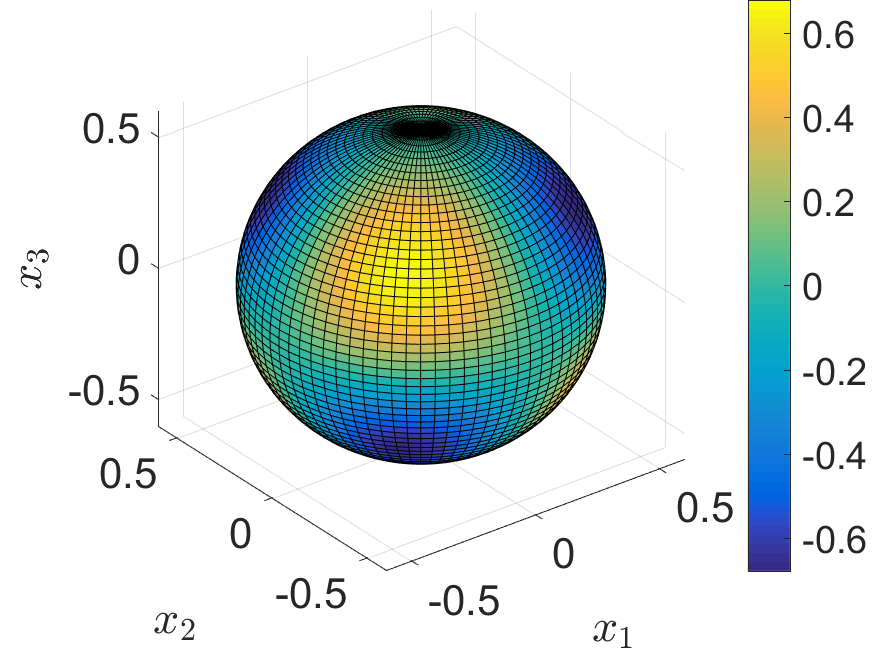}
	}
	\subfigure[Point-wise error of SPIELM]{
		\label{3D_Navier_E1:Perr2Sigmoid}
		\includegraphics[scale=0.325]{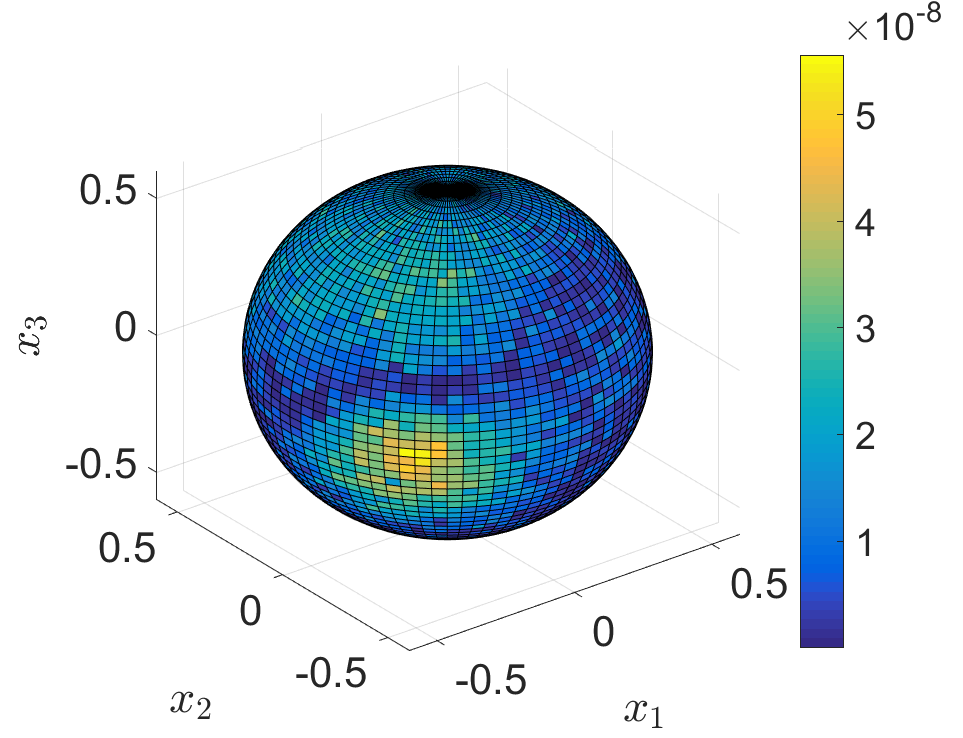}
	}
	\subfigure[Point-wise error of GPIELM]{
		\label{3D_Navier_E1:Perr2Gauss}
		\includegraphics[scale=0.32]{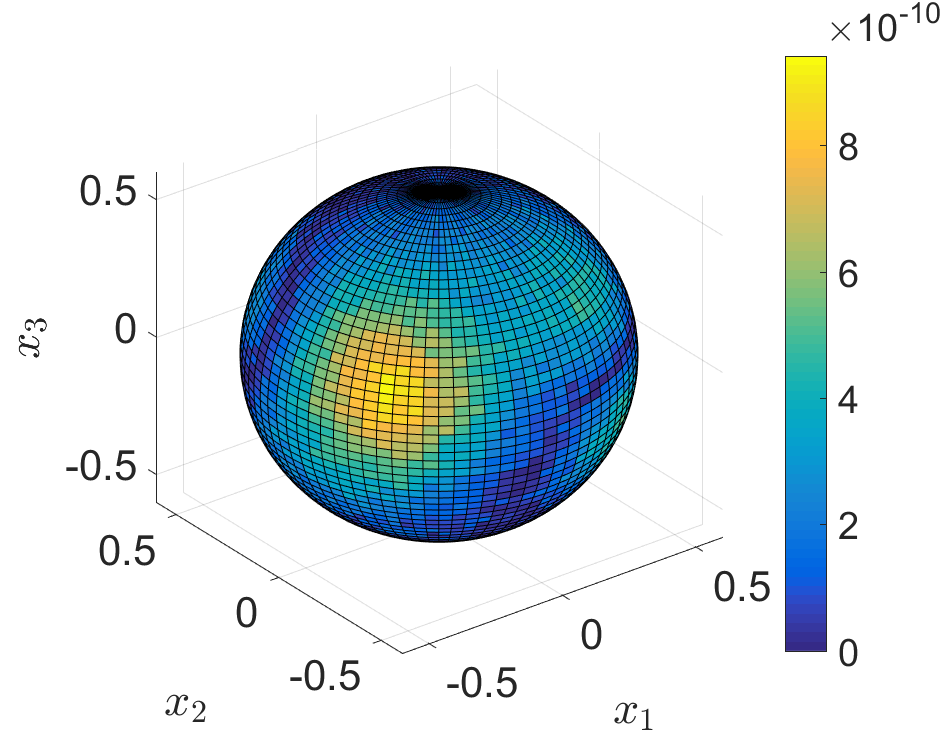}
	}
	\subfigure[Point-wise error of TPIELM]{
		\label{3D_Navier_E1:Perr2Tanh}
		\includegraphics[scale=0.32]{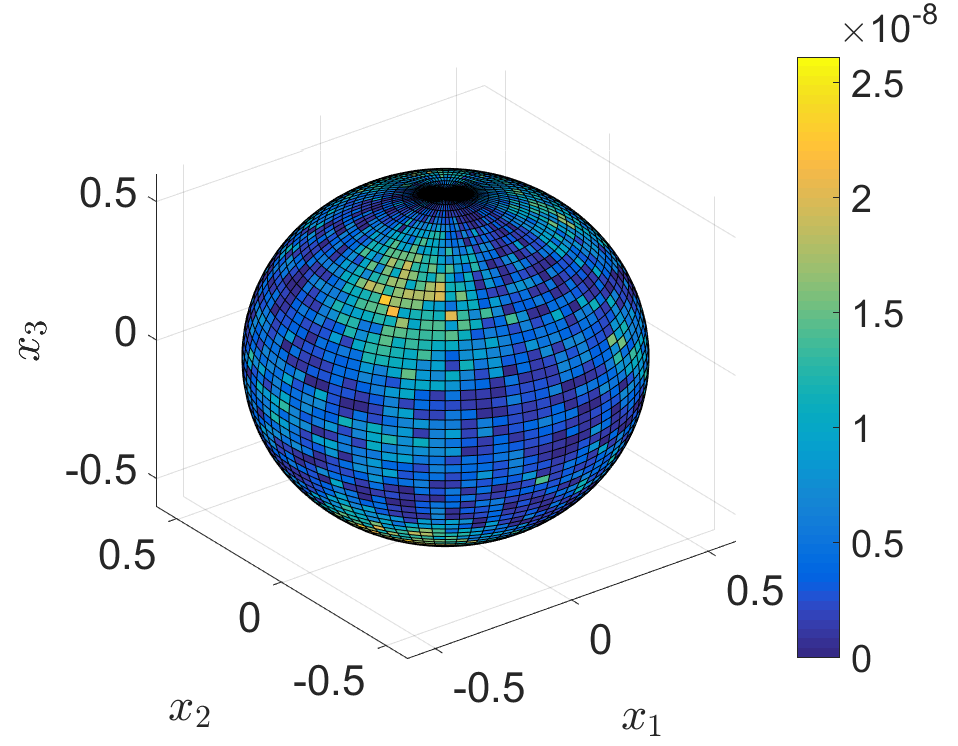}
	}
	\subfigure[Point-wise error of FPIELM]{
		\label{3D_Navier_E1:Perr2SIN}
		\includegraphics[scale=0.32]{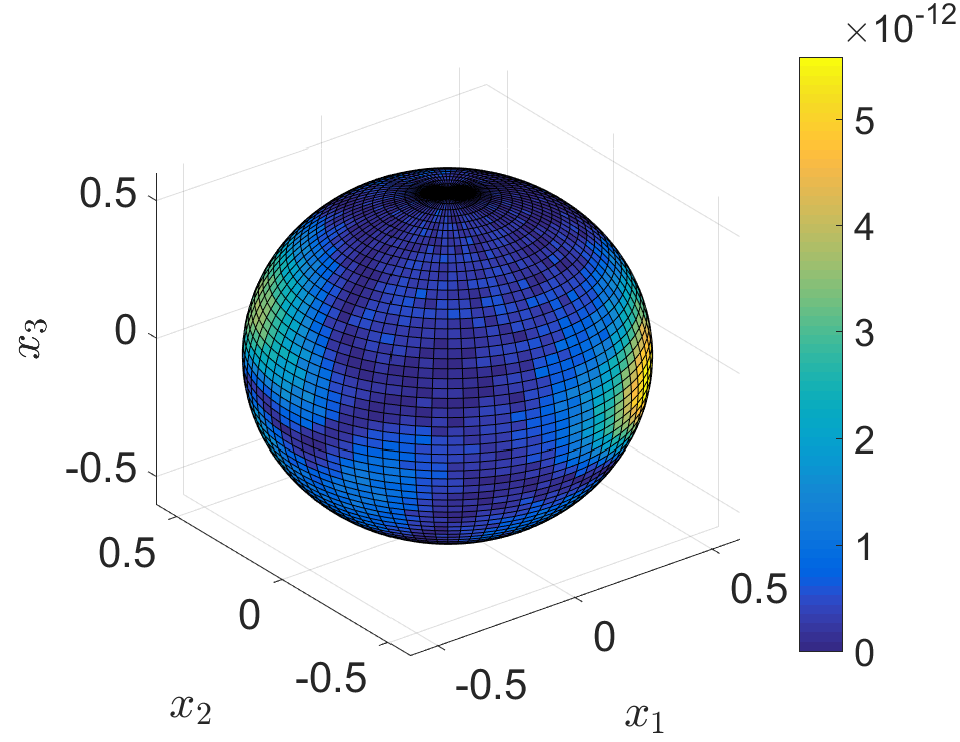}
	}
	\caption{Numerical results for four PIELM methods to Example~\ref{3D_Dirichlet_E1} on the cut planes.}
	\label{Plot_3D_Navier_E1}
\end{figure}

\begin{figure}[H]
	\centering
	\subfigure[REL VS the hidden nodes]{
		\label{3D_Navier_E1:rels} 
		\includegraphics[scale=0.325]{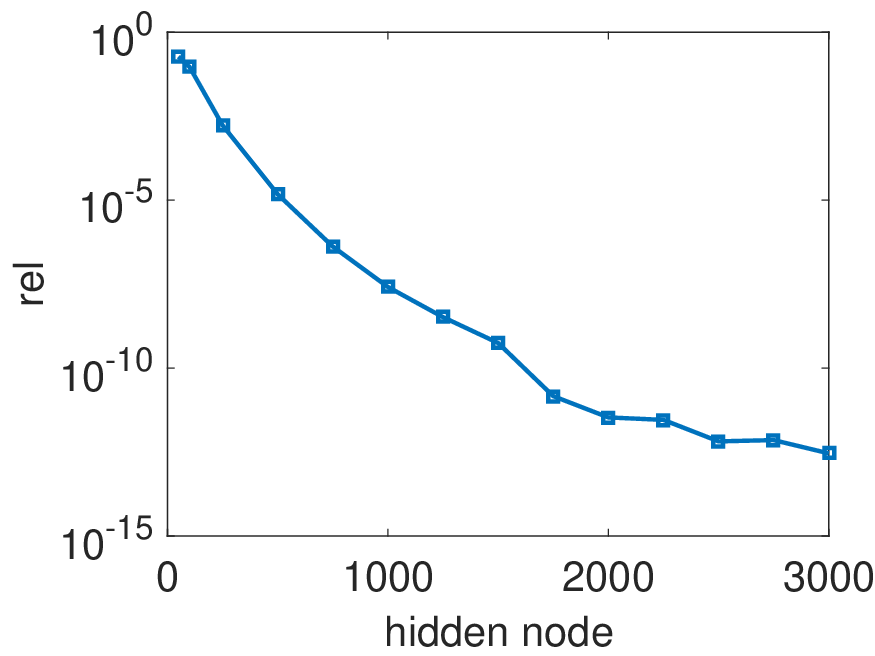}
	}
	\subfigure[Runtime VS the hidden nodes]{
		\label{3D_Navier_E1:times}
		\includegraphics[scale=0.325]{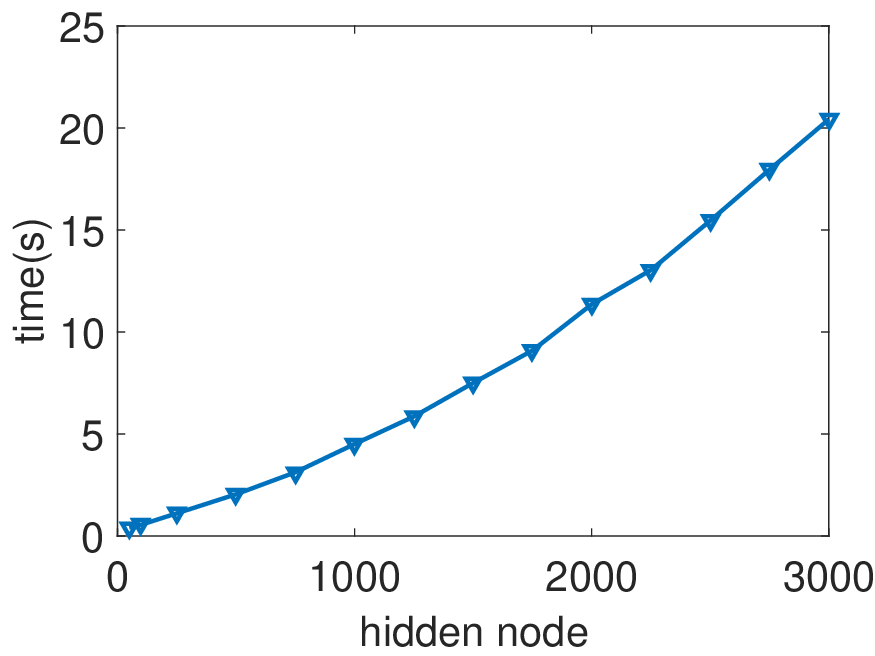}
	}
	\subfigure[REL VS the scale factor $\sigma$]{
		\label{3D_Navier_E1:sigma}
		\includegraphics[scale=0.325]{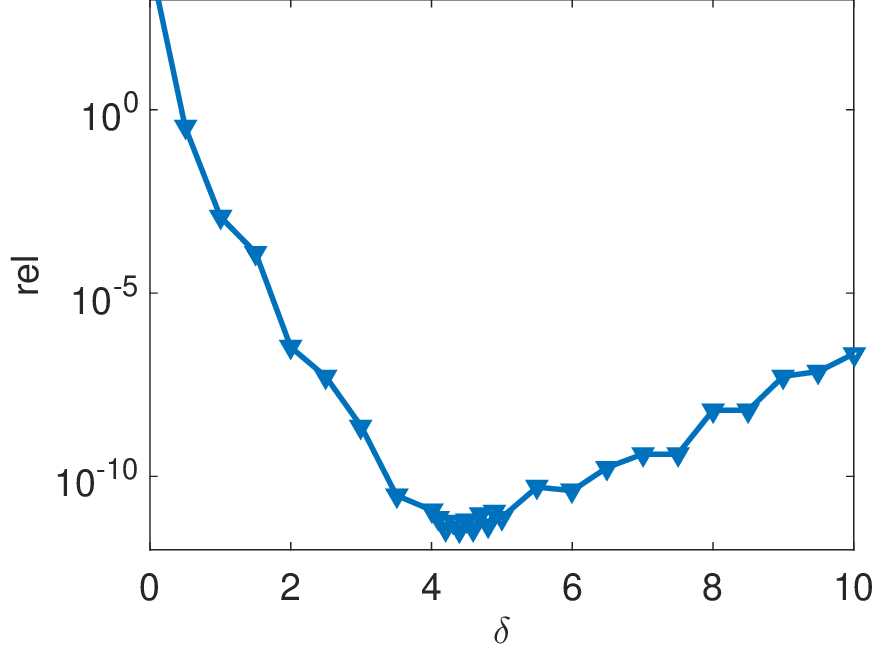}
	}
	\caption{REL VS the hidden nodes (left), runtime VS the hidden nodes (middle), and REL VS the scale factor $\sigma$ (right) to PIELM models for Example \ref{3D_Navier_E1}  on the cut planes.}
	\label{Plot_3D_Navier_E1_rel_time_Sigma}
\end{figure}

\begin{figure}[H]
	\centering
	\includegraphics[scale=0.7]{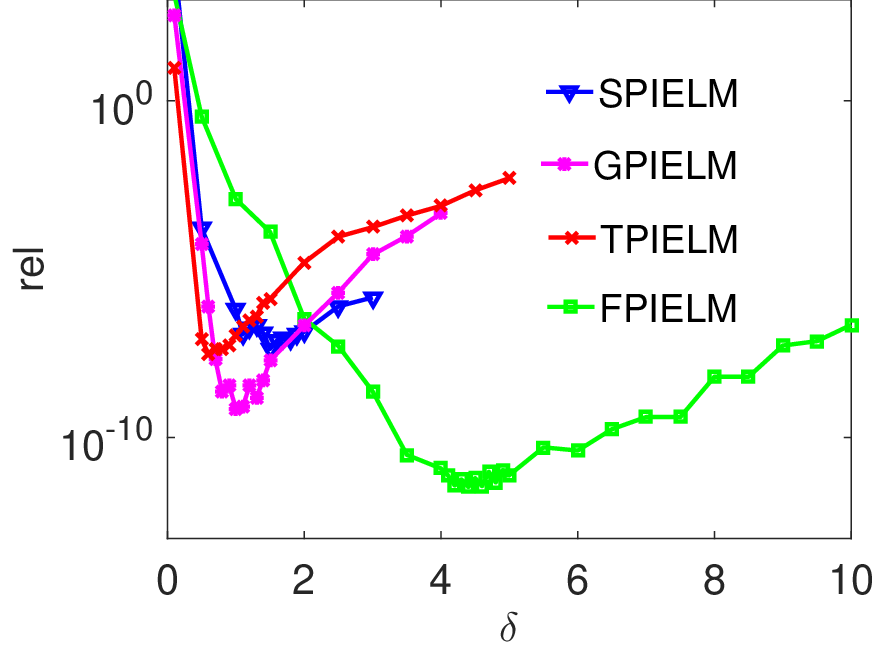}
	\caption{REL VS the scale factor $\delta$ to the four PIELM models for solving Example \ref{3D_Navier_E1} on sphere domain.}
	\label{Navier3D_E1:REL_Delta_11_PiPi} 
\end{figure}

In this example, we obtain the numerical solution of \eqref{eq:biharmonic} using the CELM model with 2,000 hidden units in $\Omega$. The other settings of the CELM model are the same as in Example~\ref{2D_Navier_E1}. The testing set consists of 3,000 mesh grid points on the spherical surface with radius $r = 1.2$. In addition to the mesh grid points, the training set includes 40,000 randomly sampled points from the domain $\Omega$ and 20,000 random points equally sampled from the inner and outer spherical surfaces. The numerical results are shown in Figs.~\ref{Plot_3D_Navier_E1} -- \ref{Navier3D_E1:REL_Delta_11_PiPi}.


For the spherical cavity in three-dimensional space, the point-wise errors in Figs.~\ref{3D_Navier_E1:Perr2Sigmoid} -- \ref{3D_Navier_E1:Perr2SIN} indicate that our proposed FPIELM method produces lower errors than those achieved using sigmoid, hyperbolic tangent, and Gaussian functions in solving the biharmonic equation \eqref{eq:biharmonic}. The curves for relative error and runtime in Figs.~\ref{Plot_3D_Navier_E1_rel_time_Sigma} and \ref{Navier3D_E1:REL_Delta_11_PiPi} demonstrate the excellent stability and robustness of the FPIELM method across varying scale factors $\sigma$, with a notably broad effective range for the scale factor $\delta$. Furthermore, the runtime of FPIELM increases linearly as the number of hidden nodes gradually grows.

\section{Conclusion}\label{sec:06}

We have introduced a Fourier-induced extreme learning machine (FPIELM) to solve biharmonic problems, utilizing the Fourier expansion of a given function within a simple ELM network. Similar to other mesh-free approaches like PINN and radial basis function models, FPIELM operates independently of a mesh and does not require an initial guess. Unlike the PINN framework, FPIELM is an efficient solution for approximating \eqref{eq:biharmonic} without needing parameter updates via iterative, gradient descent-based methods. Additionally, we examined the effects of different activation functions—tanh, Gaussian, Sigmoid, and Sine—on the performance of the PIELM model. The computational results strongly suggest that the proposed method achieves high accuracy and efficiency in solving \eqref{eq:biharmonic} across both regular and irregular domains.

It is also worth noting that random projection neural networks have been applied not only to forward problems in steady-state and time-dependent scenarios~\cite{fabiani2023parsimonious,dong2021local} but have also shown effectiveness in addressing inverse problems~\cite{galaris2022numerical, dong2023method}. These networks outperform traditional solvers~\cite{fabiani2023parsimonious} and other machine learning architectures~\cite{galaris2022numerical}, excelling in both numerical accuracy and computational efficiency. Given these advantages, future work will focus on extending our neural network models to address other PDE-related problems.

\section*{Contributions}
Xi'an Li: Conceptualization, Methodology, Investigation, Validation, Software and Writing - Original Draft; 
Jinran Wu: Conceptualization, Methodology, Validation, Writing - review $\&$ editing;  
Yujia Huang: Conceptualization, Investigation, Writing - review $\&$ editing; 
Zhe Ding: Conceptualization, Writing - review $\&$ editing; 
Xin Tai: Writing - Review \& Editing, Funding acquisition and Project administration.; 
Liang Liu: Writing - Review \& Editing and Project administration; 
You-Gan Wang: Writing - Review \& Editing and Project administration. 
All authors have read and agreed to the published version of the manuscript.

\section*{Declaration of interests}
All authors declare that they have no known competing financial interests or personal relationships that could have appeared to influence the work reported in this paper.

\section*{Acknowledgements}
This study was supported by the National Key R$\&$D Program of China (No. 2023YFF0717300).

\bibliographystyle{model1-num-names}
\bibliography{References}

\end{document}